\documentclass{amsart}
\usepackage[dvips]{epsfig}
\usepackage{graphicx}
\usepackage{latexsym}
\usepackage{amsmath}
\usepackage{amsthm}
\usepackage{amssymb}
\usepackage{mathrsfs}
\RequirePackage[colorlinks=true,citecolor=blue,urlcolor=blue]{hyperref}

\usepackage{caption}
\usepackage{subcaption}

\usepackage{eucal}

\usepackage{xcolor,soul}

\usepackage{mathtools}

 \usepackage[applemac]{inputenc}

\newtheorem{thm}{Theorem}[section]
\newtheorem{theo}[thm]{Theorem}

\newtheorem{lem}[thm]{Lemma}
\newtheorem{cor}[thm]{Corollary}
\newtheorem{defi}[thm]{Definition}
\newtheorem{hyp}[thm]{Assumption}

\theoremstyle{remark}
\newtheorem{ex}[thm]{Example}

\newtheorem{rem}[thm]{Remark}

\newcommand{\inP} {\underset{n \rightarrow \infty}{\overset{\mathbb{P}}{\xrightarrow{\hspace*{0.75cm}}}} }

\newcommand{\vip}{\vskip.2cm}
\newcommand{\COMMENTAIRE}[1]{}

\newcommand{\field}[1]{\mathbb{#1}}
\newcommand{\EE}{\field{E}}

\newcommand{\GG}{\field{G}}

\newcommand{\NN}{\field{N}}

\newcommand{\RR}{\field{R}} 
\newcommand{\TT}{\field{T}}

\newcommand{\Bb}{{\mathcal B}}

\newcommand{\Hh}{{\mathcal H}}

\newcommand{\Nn}{{\mathcal N}}
\newcommand{\Pp}{{\mathcal P}}

\newcommand{\Qq}{{\mathcal Q}}

\def \ep {\varepsilon}

\newcommand{\rd}{{\rm d}}

\newcommand{\bff}{{\bf f}}
\newcommand{\bF}{{\mathfrak f}}

\newcommand{\bH}{{\mathfrak h}}

\newcommand{\cb}{{\mathcal B}}

\newcommand{\ce}{{\mathcal E}}

\newcommand{\cf}{{\mathcal F}}

\newcommand{\ch}{{\mathcal H}}

\newcommand{\cn}{{\mathcal N}}

\newcommand{\cp}{{\mathcal P}}
\newcommand{\cq}{{\mathcal Q}}
\newcommand{\crr}{{\mathcal R}}

\newcommand{\cs}{{\mathscr S}}

\newcommand{\A}{{\mathbb A}}
\newcommand{\C}{{\mathbb C}}

\newcommand{\E}{{\mathbb E}}

\newcommand{\G}{\mathbb{G}}
\newcommand{\N}{{\mathbb N}}
\renewcommand{\P}{{\mathbb P}}
\newcommand{\p}{{\hat p}}

\newcommand{\R}{{\mathbb R}}

\newcommand{\T}{\mathbb{T}}

\newcommand{\ind}{{\bf 1}}

\newcommand{\Var}{{\mathrm{Var}}}

\newcommand{\sot}{\otimes_{\rm sym}}
\newcommand{\ssub}{\Sigma^{\rm sub}}

\newcommand{\scrit}{\Sigma^{\rm crit}}

\newcommand{\norm}[1]{\mathop{\parallel\! #1 \! \parallel}\nolimits}
\newcommand{\normm}[1]{\mathop{\parallel\! #1 \! \parallel}\nolimits_{L^{2}(\mu)}}

\newcommand{\inv}[1]{\mathop{\frac{1}{ #1}}\nolimits}
\newcommand{\expp}[1]{\mathop {\mathrm{e}^{ #1}}}

\newcommand{\reff}[1]{(\ref{#1})}

\begin{document}

\title[CLT for BMC]{Central limit theorem for bifurcating Markov chains under
  $L^{2}$-ergodic conditions}

\author{S. Val\`ere Bitseki Penda and Jean-Fran\c cois Delmas}

\address{S. Val\`ere Bitseki Penda, IMB, CNRS-UMR 5584, Universit\'e Bourgogne Franche-Comt\'e, 9 avenue Alain Savary, 21078 Dijon Cedex, France.}

\email{simeon-valere.bitseki-penda@u-bourgogne.fr}

\address{Jean-Fran\c cois Delmas, CERMICS, Ecole des Ponts,  France.}

\email{jean-francois.delmas@enpc.fr}

\begin{abstract}
  Bifurcating Markov  chains (BMC) are  Markov chains indexed by  a full
  binary tree representing  the evolution of a trait  along a population
  where each  individual has two  children.  We provide a  central limit
  theorem for additive functionals of BMC under $L^2$-ergodic conditions
  with three different regimes.  This  completes the pointwise approach developed
  in a  previous work. As application,  we study the elementary  case of
  symmetric  bifurcating  autoregressive   process,  which  justify  the
  non-trivial hypothesis considered on the kernel transition of the BMC.
  We illustrate  in this  example the phase  transition observed  in the
  fluctuations.
\end{abstract}

\maketitle






\section{Introduction}
Bifurcating Markov chains (BMC) are a class of stochastic processes indexed by
regular binary tree and which satisfy the branching Markov property (see
below for a precise definition).  This model represents the evolution of
a trait along  a population where each individual has  two children.
We refer to \cite{BD} for references on this subject. 
The recent study of BMC models was motivated by the understanding of the
cell division mechanism (where the trait of an individual is given by its
growth rate).  The first model  of BMC, named  ``symmetric'' bifurcating
auto-regressive process  (BAR), see  Section \ref{sec:sym-bar}  for more
details in  a Gaussian  framework, were introduced  by Cowan  \& Staudte
\cite{CS86} in  order to analyze cell  lineage data. In \cite{Guyon}, Guyon has studied
``asymmetric'' BAR in order to prove
statistical  evidence  of aging  in  Escherichia  Coli. 
\medskip

In this paper, our objective is to establish a central limit theorem for
additive  functionals of  BMC.   This  will be  done  for  the class  of
functions which  belong to  $L^{4}(\mu)$, where  $\mu$ is  the invariant
probability measure  associated to the associated Markov chain given by
the    genealogical evolution of
an individual  taken at random  in the population.  This  paper complete
the  pointwise  approach   developed  in  \cite{BD}  in   a  very  close
framework. Let us emphasize that the $L^2$-approach is an important step
toward  the  kernel  approximation  of   the  densities  of  the  kernel
transition of the BMC and  the invariant probability measure $\mu$ which
will be  developed in a companion  paper. The main contribution  of this
paper, with  respect to  \cite{BD}, is the  derivation of  a non-trivial
hypothesis on the kernel transition given in Assumption \ref{hyp:Q} (i).
More precisely let the random variable  $(X,Y,Z)$ model the trait of the
mother, $X$, and the traits of its  two children $Y$ and $Z$. Notice, we
do not  assume that conditionally on  $X$, the random variables  $Y$ and
$Z$ are  independent nor  have the same  distribution. In  this setting,
$\mu$  is the  distribution of  an individual  picked at  random in  the
stationary regime. From an ergodic point of view, it would be natural to
assume  some $L^2(\mu)$  continuity in  the sense  that for  some finite
constant $M$ and all functions $f$ and $g$:
\[
  \E_{X \sim \mu} [f(Y)^2 g(Z)^2]
  \leq M\,  \E_{Y\sim \mu} [f(Y)^2]\,\, \E_{Z\sim \mu} [f(Z)^2],
\]
where  $\E_{W \sim\mu}$   means  that  the  random   variable  $W$  has
distribution $\mu$. However,  this condition is not always  true even in
the simplest  case of the symmetric  BAR model, see comments  in Remarks
\ref{rem:hypPL2}    and   the    detailed    computation   in    Section
\ref{sec:sym-bar-M}.  This  motivate   the  introduction  of  Assumption
\ref{hyp:Q} (i), which  allows to recover the results  from \cite{BD} in
the context of the $L^2$ approach,  and in particular the three regimes:
sub-critical, critical and super-critical  regime. Since the results are
similar  and the  proofs  follows  the same  steps,  we  only provide  a
detailed proof in the sub-critical case.  To finish, let us mention that
the numerical  study on the  symmetric BAR, see  Section \ref{sec:simul}
illustrates the phase transitions for  the fluctuations. We also provide
an example  where the asymptotic variance  in the critical regime  is 0;
this happens when  the considered function is orthogonal  to the second
eigenspace of the associated Markov chain. 
\medskip

The paper is organized as  follows. In Section \ref{sec:BMC}, we present
the  model and  give  the assumptions:  we introduce  the  BMC model  in
Section  \ref{sec:BMC-sub},  we give  the  assumptions  under which  our
results will be stated in Section \ref{sec:L2-approach} and we give some
useful   notations    in   Section   \ref{sec:more-not}.    In   Section
\ref{sec:results}, we state  our main results: the  sub-critical case in
Section   \ref{sec:sub-critical},   the   critical   case   in   Section
\ref{sec:critical}    and   the    super-critical   case    in   Section
\ref{sec:main-res-super}. In  Section \ref{sec:sym-bar-M}, we  study the
special case of symmetric BAR process.

The  proof of  the results  in the  sub-critical case  given in  Section
\ref{sec:proof-sub-cor-L2}, which  are in the same  spirit of \cite{BD},
rely  essentially on  explicit second  moments computations  and precise
upper bounds  of fourth moments  for BMC  which are recalled  in Section
\ref{sec:moment}. The  proof of the results  in the critical case  is an
adaptation of the sub-critical space in the same spirit as in \cite{BD};
the interested reader can find the details in \cite{BD-arxiv}. The proof
of the results in the super-critical  case does not involve the original
Assumption \ref{hyp:Q} (i);  it not reproduced here as it  is very close
to its counter-part in \cite{BD}.

\section{Models and assumptions}\label{sec:BMC}

\subsection{Bifurcating Markov chain: the model}
\label{sec:BMC-sub}
We   denote   by   $\N$   the   set   of   non-negative   integers   and
$\N^*= \N  \setminus \{0\}$. If $(E,  \ce)$ is a measurable  space, then
$\cb(E)$ (resp. $\cb_b(E)$, resp.  $\cb_+(E)$) denotes the set of (resp.
bounded, resp.   non-negative) $\R$-valued measurable  functions defined
on       $E$.         For       $f\in       \cb(E)$,        we       set
$\norm{f}_\infty  =\sup\{|f(x)|, \,  x\in  E\}$.  For  a finite  measure
$\lambda$   on   $(E,\ce)$   and    $f\in   \cb(E)$   we   shall   write
$\langle \lambda,f  \rangle$ for  $\int f(x) \,  \rd\lambda(x)$ whenever
this integral is well defined.  For  $p\geq 1$ and $f\in \cb(E)$, we set
$\|f\|_{L^{p}(\lambda)}  =\langle  \lambda,|f|^p \rangle^{1/p}$  and  we
define                             the                             space
$L^{p}(\lambda)=   \left\{f  \in   \Bb(E);\,  \|f\|_{L^{p}(\lambda)}   <
  +\infty\right\}$  of $p\text{-}$integrable  functions with  respect to
$\lambda$.  For $n\in \N^*$, the product space $E^n$ is endowed with the
product $\sigma$-field $\ce^{\otimes n}$.

Let  $(S, \cs)$   be  a  measurable  space. 
Let $Q$ be a   
probability kernel   on $S \times \cs$, that is:
$Q(\cdot  , A)$  is measurable  for all  $A\in \cs$,  and $Q(x,
\cdot)$ is  a probability measure on $(S,\cs)$ for all $x \in
S$. For any $f\in \cb_b(S)$,   we set for $x\in S$:
\begin{equation}
   \label{eq:Qf}
(Qf)(x)=\int_{S} f(y)\; Q(x,\rd y).
\end{equation}
We define $(Qf)$, or simply $Qf$, for $f\in \cb(S)$ as soon as the
integral \reff{eq:Qf} is well defined, and we have $\cq f\in \cb(S)$. For $n\in \N$, we denote by $Q^n$  the
$n$-th iterate of $Q$ defined by $Q^0=I_d$, the identity map on
$\cb(S)$, and $Q^{n+1}f=Q^n(Qf)$ for $f\in \cb_b(S)$.  

Let $P$ be a   
probability kernel   on $S \times \cs^{\otimes 2}$, that is:
$P(\cdot  , A)$  is measurable  for all  $A\in \cs^{\otimes 2}$,  and $P(x,
\cdot)$ is  a probability measure on $(S^2,\cs^{\otimes 2})$ for all $x \in
S$. For any $g\in \cb_b(S^3)$ and $h\in \cb_b(S^2)$,   we set for $x\in S$:
\begin{equation}
   \label{eq:Pg}
(Pg)(x)=\int_{S^2} g(x,y,z)\; P(x,\rd y,\rd z)
\quad\text{and}\quad
(Ph)(x)=\int_{S^2} h(y,z)\; P(x,\rd y,\rd z).
\end{equation}
We define $(Pg)$ (resp. $(Ph)$), or simply $Pg$ for $g\in \cb(S^3)$
(resp. $Ph$ for $h\in \cb(S^2)$), as soon as the corresponding 
integral \reff{eq:Pg} is well defined, and we have  that $Pg$ and
$Ph$ belong to $\cb(S)$.
\medskip 

We  now introduce  some notations  related to  the regular  binary tree.
 We   set   $\T_0=\G_0=\{\emptyset\}$,
$\G_k=\{0,1\}^k$  and $\T_k  =  \bigcup _{0  \leq r  \leq  k} \G_r$  for
$k\in  \N^*$, and  $\T  =  \bigcup _{r\in  \N}  \G_r$.   The set  $\G_k$
corresponds to the  $k$-th generation, $\T_k$ to the tree  up to the $k$-th
generation, and $\T$ the complete binary  tree. For $i\in \T$, we denote
by $|i|$ the generation of $i$ ($|i|=k$  if and only if $i\in \G_k$) and
$iA=\{ij; j\in A\}$  for $A\subset \T$, where $ij$  is the concatenation
of   the  two   sequences  $i,j\in   \T$,  with   the  convention   that
$\emptyset i=i\emptyset=i$.

We recall the definition of bifurcating Markov chain  from
\cite{Guyon}. 
\begin{defi}
  We say  a stochastic process indexed  by $\T$, $X=(X_i,  i\in \T)$, is
  a bifurcating Markov chain (BMC) on a measurable space $(S, \cs)$ with
  initial probability distribution  $\nu$ on $(S, \cs)$ and probability
  kernel $\cp$ on $S\times \cs^{\otimes 2}$ if:
\begin{itemize}
\item[-] (Initial  distribution.) The  random variable  $X_\emptyset$ is
  distributed as $\nu$.
   \item[-] (Branching Markov property.) For  a sequence   $(g_i, i\in
     \T)$ of functions belonging to $\cb_b(S^3)$, we have for all $k\geq 0$,
\[
\E\Big[\prod_{i\in \G_k} g_i(X_i,X_{i0},X_{i1}) |\sigma(X_j; j\in \T_k)\Big] 
=\prod_{i\in \G_k} \cp g_i(X_{i}).
\]
\end{itemize}
\end{defi}

Let $X=(X_i,  i\in \T)$ be a BMC  on a measurable space $(S, \cs)$ with
  initial probability distribution  $\nu$ and probability
  kernel $\cp$. 
We    define     three
probability kernels $P_0, P_1$ and $\cq$ on $S\times \cs$ by:
\[
P_0(x,A)=\cp(x, A\times S), \quad
P_1(x,A)=\cp(x, S\times A) \quad\text{for
$(x,A)\in S\times  \cs$, and}\quad
\cq=\inv{2}(P_0+P_1).
\] 
Notice  that  $P_0$ (resp.   $P_1$)  is  the  restriction of  the  first
(resp. second) marginal of $\cp$ to $S$.  Following \cite{Guyon}, we
introduce an  auxiliary Markov  chain $Y=(Y_n, n\in  \N) $  on $(S,\cs)$
with  $Y_0$ distributed  as $X_\emptyset$  and transition  kernel $\cq$.
The  distribution of  $Y_n$ corresponds  to the  distribution of  $X_I$,
where $I$  is chosen independently from  $X$ and uniformly at  random in
generation  $\G_n$.    We  shall   write  $\E_x$   when  $X_\emptyset=x$
(\textit{i.e.}  the initial  distribution  $\nu$ is  the  Dirac mass  at
$x\in S$).  \medskip

We end this section with a useful inequality and the Gaussian BAR model.

\begin{rem}
   \label{rem:ineq-simple}
By convention, for  $f,g\in  \cb(S)$,   we  define  the  function
$f\otimes  g\in \cb(S^2)$  by $(f\otimes g)(x,y)=f(x)g(y)$ for  $x,y\in S$ and
introduce the notations:
\[
f\sot g= \inv{2}(f\otimes g + g\otimes f) \quad\text{and}\quad f\otimes ^2= f\otimes f.
\]
Notice that 
$\cp(g\sot \ind)=\cq(g)$ for  $g\in \cb_+(S)$.
 For $f \in \cb_+(S)$, as $ f\otimes f \leq f^2 \sot \ind $, we get:
\begin{equation}\label{eq:majo-pfxf}
\cp(f\otimes^2)=\cp(f\otimes f) \leq  \cp(f^2\sot \ind) = \cq\left(f^2\right). 
\end{equation}
\end{rem}

\begin{ex}[Gaussian bifurcating autoregressive process] 
\label{ex:valide} 
We will consider the real-valued Gaussian bifurcating autoregressive
process (BAR) $X=(X_{u},u\in\TT)$ where for all $ u \in \T$: 
\[
\begin{cases}
X_{u0} = a_{0} X_{u} + b_{0} + \ep_{u0},  \\
 X_{u1} = a_{1} X_{u} + b_{1} + \ep_{u1},
\end{cases}
\]
with $a_{0}, a_{1} \in (-1,1)$, $b_{0}, b_{1} \in \RR$ and
$((\ep_{u0},\ep_{u1}),\, u \in \T)$  an
independent  sequence of bivariate Gaussian  $\Nn(0,\Gamma)$ random vectors
independent of $X_{\emptyset}$ with covariance matrix, with $\sigma>0$
and $\rho\in \R$ such that $|\rho|\leq  \sigma^2$:
\begin{equation*}
\Gamma = \begin{pmatrix} \sigma^{2} \quad \rho \\ \rho \quad \sigma^{2} \end{pmatrix}.
\end{equation*}
Then the process $X=(X_{u},u\in\TT)$ is a BMC with transition
probability $\Pp$ given by:
\[
\Pp(x,dy,dz) = \frac{1}{2\pi \sqrt{\sigma^{4} - \rho^{2}}}
\, \exp\left( -\frac{\sigma^{2}}{2(\sigma^{4} -
  \rho^{2})}\,  g(x,y,z) \right)\, dydz,
\]
with 
\[
g(x,y,z)=(y-a_{0}x-b_{0})^{2}  - 2\rho\sigma^{-2}(y - a_{0}x -
b_{0})(z - a_{1}x - b_{1}) + (z - a_{1}x - b_{1})^{2}.
\]
The transition kernel $\Qq$ of the auxiliary Markov  chain is defined by:
\begin{equation*}
\Qq(x,dy) = \frac{1}{2\sqrt{2\pi\sigma^2}}\left(\expp{-(y - a_{0}x
    - b_{0})^{2}/2\sigma^{2}}  + \expp{-(y - a_{1}x -
    b_{1})^{2}/2\sigma^{2}}\right)\, dy.
\end{equation*}
\end{ex}

\subsection{Assumptions}
\label{sec:L2-approach}
We assume that $\mu$ is an invariant probability measure for $\cq$. 

\medskip

We state first some regularity assumptions on
the kernels $\cp$ and $\cq$ and the invariant measure $\mu$ we will use
later on. Notice first that by Cauchy-Schwartz we have for $f,g\in L^4(\mu)$:
\[
|\cp ( f\otimes
g)| ^2\leq  \cp (f^2 \otimes 1) \,\cp(1\otimes g^2)
\leq 4 \cq(f^2)\, \cq(g^2),
\]
so that, as $\mu$ is an invariant measure of $\cq$:
\begin{equation}
   \label{eq:mp2}
\normm{ \cp ( f\otimes
g)} \leq 2 \normm{\cq(f^2)}^{1/2} \,\normm{\cq (g^2)}^{1/2} \leq 2
\norm{f}_{L^4(\mu)}\, \norm{g}_{L^4(\mu)},
\end{equation}
and similarly for $f,g\in L^2(\mu)$:
\begin{equation}
   \label{eq:mp}
\langle \mu,  \cp ( f\otimes
g)\rangle \leq 
2\normm{f}\, \normm{g}.
\end{equation}
We shall in fact assume that $\cp$ (in fact only its symmetrized
version)  is in a sense an $L^2(\mu)$ operator,
see also Remark \ref{rem:hypPL2} below.

\begin{hyp}\label{hyp:Q}
There exists an invariant probability measure, $\mu$, for the Markov transition kernel $\cq$.
\begin{itemize}
\item[(i)]
There exists a finite constant $M$ such that for all $f,g,h \in L^{2}(\mu)$:
\begin{align}
\label{eq:P-L2-1}
\normm{\cp (\cq f \sot \cq g)}
&\leq M  \normm{f}  \normm{g} ,\\
\label{eq:P-L2-3}
\normm{\cp \left(\cp(\cq f \sot \cq g)\sot \cq h\right)}
&\leq M  \normm{f}  \normm{g} \normm{h} ,\\
\label{eq:P-L2-2}
\normm{\cp ( f \sot  \cq g)}
&\leq M  \norm{f}_{L^{4}(\mu)}   \normm{g}.
\end{align}
\item[(ii)]  There exists $k_0\in \N$, such that the  probability measure $\nu \cq^{k_0}$ has a bounded density, say $\nu_0$, with respect to $\mu$. That is:
\begin{equation*}
\nu \cq^{k_0}(dy) = \nu_0(y) \mu(y) \, dy
\quad\text{and} \quad
\norm{\nu_0}_\infty <+\infty .
\end{equation*}
\end{itemize}
\end{hyp}

\begin{rem}\label{rem:hypPL2}
Let $\mu$ be an invariant probability measure of $\cq$.  If  there exists a finite constant $M$ such that for all $f,g \in
L^{2}(\mu)$:
\begin{equation}
\label{eq:P-L2-0}
\normm{ \cp ( f\otimes g)} \leq  M \normm{f}\normm{g},
\end{equation}
then we deduce that \reff{eq:P-L2-1}, \reff{eq:P-L2-3} and \reff{eq:P-L2-2}
hold. Condition \reff{eq:P-L2-0} is much more natural and simpler than
the latter ones, and it allows to give shorter
proofs. However
Condition \reff{eq:P-L2-0} appears to be too strong even in the simplest
case of the symmetric BAR model developed in Example
\ref{ex:valide} with $a_0=a_1$ and $b_0=b_1$. Let $a$ denote the common
value of $a_0$ and $a_1$. In fact, according to the value of $a\in (-1, 1)$ in the symmetric BAR model,  there exists
$k_1\in \N$ such that for all $f,g \in
L^{2}(\mu)$
\begin{equation}
\label{eq:P-L2-11}
\normm{\cp (\cq^{k_1} f \otimes \cq^{k_1} g)}
\leq M  \normm{f}  \normm{g} ,
\end{equation}
with $k_1$  increasing with $|a|$. Since Assumption \ref{hyp:Q} (i) is only necessary for the asymptotic normality in the case
$|a|\in [0, 1/\sqrt{2}]$ (corresponding to the sub-critical and critical
regime), it will be enough to consider $k_1=1$ (but not
sufficient  to   consider  $k_1=0$).   For  this   reason,  we  consider
\reff{eq:P-L2-1},   that   is  \eqref{eq:P-L2-11}  with
$k_1=1$.    A   similar   remark    holds   for   \reff{eq:P-L2-3}   and
\reff{eq:P-L2-2}. In a sense Condition  \reff{eq:P-L2-11} (and similar
extensions of \eqref{eq:P-L2-3} and  \eqref{eq:P-L2-2})
is in the same
spirit  as item  (ii) of  Assumption \ref{hyp:Q}:  ones use  iterates of
$\cq$ to get smoothness on the kernel $\cp$ and the initial distribution
$\nu$.
\end{rem}

\begin{rem}\label{rem:hypP}
  Let  $\mu$ be  an  invariant probability measure  of $\cq$  and  assume that  the
  transition kernel $\cp$ has a density, denoted by $p$, with respect to
  the       measure       $\mu^{\otimes        2}$,       that       is:
  $\cp(x, dy, dz)  = p(x,y,z) \, \mu(dy)\mu(dz)$ for all  $x\in S$. Then
  the  transition kernel  $\cq$  has  a density,  denoted  by $q$,  with
  respect  to $\mu$,  that  is:  $\cq(x, dy)=  q(x,y)  \mu(dy)$ for  all
  $x   \in   S$   with   $q(x,y)=2^{-1}   \int_S   (p(x,y,z)+p(x,z,y))\,
  \mu(dz)$. We set:
\begin{equation}
   \label{eq:def-bH}
\bH(x) = \left(\int_{S} q(x,y)^{2}\,  \mu(dy)\right)^{1/2}.
\end{equation}
Assume that:
\begin{align}
  \label{eq:Hil-SchP-1}
\normm{\cp(\bH\otimes^2)} &<+\infty ,\\
  \label{eq:Hil-SchP-2}
\normm{\cp(\cp(\bH\otimes^2)\sot \bH)} &<+\infty ,
\end{align}
and that  there exists a finite
constant $C$ such that for all $f\in L^4(\mu)$:
\begin{equation}
  \label{eq:Hil-SchP-3}
\normm{\cp(f\sot \bH)} \leq C \,\norm{f}_{L^4(\mu)}.
\end{equation}
Since $|\cq f| \leq \normm{f} \bH$, we deduce that
\eqref{eq:Hil-SchP-1}, \eqref{eq:Hil-SchP-2} and \eqref{eq:Hil-SchP-3}
 imply respectively  \reff{eq:P-L2-1}, 
\eqref{eq:P-L2-3} and \eqref{eq:P-L2-2}. 
\end{rem}

We consider the following ergodic properties of $\cq$, which in
particular  implies that $\mu$ is indeed the unique invariant
probability measure for $\cq$.   We refer to
\cite{dmps:mc} Section 22 for a detailed account on
$L^2(\mu)$-ergodicity (and in particular Definition 22.2.2 on exponentially convergent  
Markov kernel). 

\begin{hyp}\label{hyp:Q1}
The Markov kernel $\cq$ has an (unique) invariant probability measure $\mu$, and
$\cq$ is $L^2(\mu)$ exponentially convergent, that is there exists
$\alpha \in (0,1)$ and $M$ finite such that for all $f \in L^{2}(\mu)$:
\begin{equation}
\label{eq:L2-erg}
\normm{\cq^{n}f - \langle \mu, f\rangle}\leq M \alpha^{n} \normm{f} \quad \text{for all $n \in \NN$}.
\end{equation}
\end{hyp}

We consider the stronger ergodic property based on a second spectral gap. (Notice in particular that
Assumption \ref{hyp:Q2} implies Assumption \ref{hyp:Q1}.)
\begin{hyp}\label{hyp:Q2}
The Markov kernel $\cq$ has an (unique)  invariant probability measure $\mu$, and there exists $\alpha \in (0,1)$,   a finite  non-empty set
  $J$ of indices, distinct complex eigenvalues $\{\alpha_j, \, j\in J\}$
  of  the  operator  $\cq$ with  $|\alpha_j|=\alpha$,  non-zero  complex
  projectors $\{\crr_j, \, j\in J\}$  defined on $\C L^2(\mu)$, the $\C$-vector
  space    spanned by        $L^2(\mu)$,       such        that
  $\crr_j\circ \crr_{j'}=\crr_{j'}\circ  \crr_{j}=0$ for all  $j\neq j'$
  (so that  $\sum_{j\in J} \crr_j$ is  also a projector defined  on $\C L^2(\mu)$)
  and a positive  sequence $(\beta_n, n\in \N)$ converging  to $0$ such
  that  for  all  $f\in  L^2(\mu)$,  with
  $\theta_j=\alpha_j/\alpha$:
\begin{equation}\label{eq:L2-erg-crit}
\|\Qq^{n}f - \langle \mu, f\rangle - \alpha^{n}\sum_{j\in J} \theta_{j}^{n} \crr_{j}(f)\|_{L^{2}(\mu)} \leq \beta_{n}\alpha^{n} \|f\|_{L^{2}(\mu)} \quad \text{for all $n \in \NN$}.
\end{equation}
\end{hyp}

Assumptions \ref{hyp:Q1} and \ref{hyp:Q2} stated in an $L^2$ framework
corresponds to  \cite[Assumptions 2.4 and 2.6]{BD} stated in a pointwise
framework. The structural Assumption \ref{hyp:Q} on the transition
kernel $\cp$ replace the structural \cite[Assumptions 2.2]{BD} on the
set of considered functions. 

\begin{rem}\label{rem:HS-0}
Assume that $\cq$ has a density $q$ with respect to an invariant
probability  measure $\mu$
such that $\bH\in L^2(\mu)$, where $\bH$ is defined in
\eqref{eq:def-bH}, 
that is:
\begin{equation*}\label{eq:Hilbert-SchQ}
\int_{S^2} q(x,y)^2 \, \mu(dx)\mu(dy) < +\infty. 
\end{equation*}
Then the operator $\Qq$ is  a non-negative Hilbert-Schmidt operator (and
then a  compact operator) on $L^2(\mu)$.  It is well known  that in this
case, except for the possible value  $0$, the spectrum of $\Qq$ is equal
to the set $\sigma_{p}(\Qq)$  of eigenvalues of $\Qq$; $\sigma_{p}(\Qq)$
is a countable set with $0$  as the only possible accumulation point and
for  all $\lambda  \in \sigma_{p}(\Qq)\setminus  \{0\}$, the  eigenspace
associated  to $\lambda$  is finite-dimensional  (we refer  for e.g.  to
\cite[chap. 4]{b:iotis} for more details).
In particular, if 1 is the only eigenvalue of $\cq$ with modulus 1
  and if it has multiplicity 1 (that is the corresponding eigenspace is
  reduced to the constant functions),
then  Assumptions \ref{hyp:Q1} and \ref{hyp:Q2} also hold. Let us
mention that  $q(x,y)>0$  $\mu(dx)\otimes \mu(dy)$-a.s. is a standard
condition which  implies that  1 is the only eigenvalue of $\cq$ with modulus 1
  and that it  has multiplicity 1, see for example \cite{BR95}.
\end{rem}

\medskip

\subsection{Notations for average of different functions over different generations}
\label{sec:more-not}

Let  $X=(X_u, u\in  \T)$ be  a BMC  on $(S,  \mathcal{S})$ with  initial
 probability distribution $\nu$, and probability kernel $\cp$.  Recall $\cq$ is the induced Markov kernel. We shall assume that $\mu$ is an invariant probability measure of $\cq$. For a finite set $A\subset \T$ and a function $f\in \cb(S)$, we set:
\[
M_A(f)=\sum_{i\in A} f(X_i).
\]
We shall be interested in the cases $A=\G_n$ (the  $n$-th generation)
and $A=\T_n$ (the tree up to the $n$-th generation). We recall  from
\cite[Theorem~11 and Corollary~15]{Guyon}  that under geometric
ergodicity assumption, we have for $f$ a continuous bounded real-valued
function defined on $S$, the following convergence in $L^2(\mu)$ (resp. a.s.):
 \begin{equation}\label{eq:lfgn-G}
 \lim_{n\rightarrow\infty } |\G_n|^{-1} M_{\G_n}(f)=\langle \mu, f
 \rangle
 \quad\text{and}\quad
 \lim_{n\rightarrow\infty } |\T_n|^{-1} M_{\T_n}(f)=\langle \mu, f
 \rangle.
 \end{equation}
Using Lemma \ref{lem:L2MG} and the Borel-Cantelli Theorem, one can prove
that we also have \eqref{eq:lfgn-G} with the $L^2(\mu)$ and
a.s. convergences under Assumptions \ref{hyp:Q}-(ii) and \ref{hyp:Q1}. 

We shall now consider the corresponding fluctuations.  We will use
frequently the following notation: 
\[
\boxed{\tilde f= f - \langle \mu,
    f \rangle\quad \text{for $f\in L^1(\mu)$.}}
\]

Recall that for $f\in L^1(\mu)$, we set $\tilde f= f - \langle \mu,
    f \rangle$. In order to study the asymptotics of $ M_{\GG_{n-\ell }}(\tilde
f)$,  we shall consider the contribution of the descendants of the
individual $i\in \T_{n-\ell}$ for  $n\geq \ell\geq 0$:
\begin{equation}
   \label{eq:def-Nnil}
N^\ell_{n,i}(f)=|\G_n|^{-1/2} M_{i\G_{n-|i|-\ell}}(\tilde f), 
\end{equation}
where  $i\G_{n-|i|-\ell}=\{ij, \, j\in
\G_{n-|i|-\ell}\}\subset \G_{n-\ell}$. For  all $k\in \N$ such that $n\geq k+\ell$, we have:
\[
M_{\G_{n-\ell}}(\tilde f)=\sqrt{|\G_n|}\,\, \sum_{i\in \G_k}
N^\ell_{n,i}(f)=
\sqrt{|\G_n|}\, \,N_{n, \emptyset}^\ell(f).
\]
Let $\bF=(f_\ell, \ell\in \N)$ be a sequence of elements of
$L^1(\mu)$. We set
for $n\in \N$ and $i\in \T_n$:
\begin{equation}
   \label{eq:def-NiF}
N_{n,i}(\bF)=\sum_{\ell=0}^{n-|i|} N_{n,i}^\ell(f_\ell) 
=|\G_n|^{-1/2 }\sum_{\ell=0}^{n-|i|}
M_{i\G_{n-|i|-\ell}}(\tilde f_\ell).
\end{equation}
We deduce that $ \sum_{i\in \G_k} N_{n,i}(\bF)=|\G_n|^{-1/2 }\sum_{\ell=0}^{n-k}
M_{\G_{n-\ell}}(\tilde f_\ell)$
which gives for $k=0$:
\begin{equation}
   \label{eq:def-NOf1}
\boxed{N_{n, \emptyset}(\bF)= |\G_n|^{-1/2 }\sum_{\ell=0}^{n}
  M_{\G_{n-\ell}}(\tilde f_{\ell}).}
\end{equation}
The notation $N_{n, \emptyset}$ means  that we consider the average from
the root $\emptyset$  to the $n$-th generation.

\begin{rem}
   \label{rem:simpleN0n}
 We shall consider in particular the following two simple cases. 
Let      $f\in     L^1(\mu)$      and     consider      the     sequence
$\bF=(f_\ell,  \ell\in \N)$.   If  $f_0=f$ and $f_\ell=0$ for  $\ell
\in \N^*$,  then we get:
\[
N_{n, \emptyset}(\bF)= |\G_n|^{-1/2} M_{\G_n}(\tilde f).
\]
  If $f_\ell=f$ for $\ell \in \N$, then we shall write  $\bff=(f,f,
  \ldots)$, and we get, as  $|\T_n|=2^{n+1} - 1 $ and $|\G_n|= 2^n$:
\[
N_{n, \emptyset}(\bff)= |\G_n|^{-1/2} M_{\T_n}(\tilde f)
=  \sqrt{2 - 2^{-n}}\,\, |\T_n|^{-1/2} M_{\T_n}(\tilde f).
\]
Thus, we will deduce the
fluctuations of $M_{\T_n}(f)$ and $M_{\G_n}(f)$ from the asymptotics of $N_{n,
  \emptyset}(\bF)$. 
\end{rem}

\medskip 

Because of condition (ii) in  Assumption \ref{hyp:Q} which roughly state
that after  $k_0$ generations,  the distribution  of the  induced Markov
chain is  absolutely continuous  with respect  to the  invariant measure
$\mu$, it  is better to consider  only generations $k\geq k_0$  for some
$k_0\in  \N$ and  thus  remove  the first  $k_0-1$  generations in  the
quantity $N_{n,\emptyset}(\bF) $ defined in \reff{eq:def-NOf1}. 

To study the asymptotics of $N_{n, \emptyset}(\bF)$, it is convenient to
write for $n\geq k\geq 1$:
\begin{equation}
   \label{eq:nof-D}
N_{n, \emptyset}(\bF)= |\G_n|^{-1/2} \sum_{r=0}^{k-1} M_{\G_r}(\tilde
f_{n-r}) + \sum_{i\in \G_k} N_{n,i}(\bF).
\end{equation}
If $\bff=(f,f, \ldots)$ is the infinite sequence of the same function
$f$, this becomes:
\begin{equation*}
   \label{eq:Nnff}
N_{n, \emptyset}(\bff)= |\G_n|^{-1/2} M_{\T_n}(\tilde f)= |\G_n|^{-1/2}
M_{\T_{k-1}}(\tilde f)+ \ \sum_{i\in
  \G_k} N_{n,i}(\bff).
\end{equation*}

\section{Main results}\label{sec:results}
 
\subsection{The sub-critical case: $2\alpha^2<1$}\label{sec:sub-critical}
We shall consider, when well defined, for a sequence $\bF=(f_\ell,
\ell\in  \N)$ of measurable real-valued
functions defined on $S$, the quantities:
\begin{equation}\label{eq:ssub}
\ssub(\bF)=\ssub_1(\bF)+ 2 \ssub_2(\bF), 
\end{equation}
where:
\begin{align}
   \label{eq:S1}
   \ssub_1(\bF)
&=\sum_{\ell\geq 0} 
2^{-\ell} \, \langle \mu,   \tilde f_\ell^ 2\rangle
+ \sum_{\ell\geq 0, \, k\geq 0} 
2^{k-\ell} \, \langle \mu, \cp\left((\cq^k \tilde
    f_\ell) \otimes^2\right)\rangle,\\
   \label{eq:S2}
   \ssub_2(\bF)
&=\sum_{0\leq \ell< k}
2^{-\ell} \langle \mu,   \tilde f_k \cq^{k-\ell} \tilde f_\ell\rangle
+ \sum_{\substack{0\leq \ell< k\\ r\geq 0}} 2^{r-\ell} \langle \mu,
\cp\left( \cq^r \tilde
    f_k \sot
\cq^{k-\ell+r} \tilde f_\ell  \right)\rangle.
\end{align}
\medskip

The proof of the next result is detailed in Section
\ref{sec:proof-sub-cor-L2}. 
\begin{theo}\label{cor:subcritical-L2}
Let $X$ be a BMC with kernel $\cp$ and initial distribution $\nu$ such that Assumptions  \ref{hyp:Q} and \ref{hyp:Q1} are in  force with  $\alpha\in (0,  1/\sqrt{2})$.  We  have the  following convergence      in      distribution       for      all      sequence $\bF=(f_\ell,   \ell\in   \N)$   bounded  in   $L^4(\mu)$   (that   is $\sup_{\ell\in \N} \norm{f_\ell}_{L^4(\mu)}<+\infty $):
\[
N_{n, \emptyset}(\bF) \; \xrightarrow[n\rightarrow \infty
]{\text{(d)}} \; G, 
\]
where $G$ is centered Gaussian   random variable with 
variance $\ssub(\bF)$ given by \reff{eq:ssub} which is well defined and finite. 
\end{theo}

Notice   that  the   variance  $\ssub(\bF)$   already  appears   in  the
sub-critical  pointwise  approach  case,   see  \cite[(15)  and  Theorem
3.1]{BD}. Then, arguing similarly as in \cite[Section 3.1]{BD}, we deduce that if
Assumptions \ref{hyp:Q} and \ref{hyp:Q1} are in force
with $\alpha\in  (0, 1/\sqrt{2})$, then  for $f\in L^{4}(\mu)$,  we have
the following convergence in distribution:
\begin{equation}
  \label{eq:subcriticalGT}
|\GG_{n}|^{-1/2}M_{\GG_{n}}(\tilde{f}) \; \xrightarrow[n\rightarrow \infty ]{\text{(d)}} \; G_{1}
  \quad \text{and} \quad
  |\TT_{n}|^{-1/2}M_{\TT_{n}}(\tilde{f}) \; \xrightarrow[n\rightarrow
  \infty ]{\text{(d)}} \; G_{2},  
\end{equation}
where $G_{1}$ and $G_2$ are  centered Gaussian  random variables with
respective  
variances   $\ssub_{\GG}(f)=\ssub(\bF)$, with  $\bF = (f,0,0, \ldots)$,
and $\ssub_{\TT}(f)=\ssub(\bff)/2$  with $ \bff=(f,f, \ldots)$,
given in  \cite[Corollary 3.3]{BD}  which are  well defined
and finite.

\subsection{The critical case: $2\alpha^2=1$}\label{sec:critical}
In the critical  case $\alpha=1/\sqrt{2}$, we shall  denote by
$\crr_j$  the projector on the eigen-space  associated to the eigenvalue
$\alpha_j$ with $\alpha_j=\theta_j \alpha$, $|\theta_j|=1$ and for $j$
in the finite set of indices $J$. Since $\cq$  is a real operator, we get that if $\alpha_j$
is a non real eigenvalue, so is $\overline  \alpha_j$.   We  shall
denote  by  $\overline  \crr_j$  the projector  associated  to
$\overline  \alpha_j$.   Recall  that the sequence $(\beta_n, n\in \N)$
in Assumption \ref{hyp:Q2} is non-increasing and bounded from above by 1. For all measurable real-valued function $f$ defined on $S$, we set,
when this is well defined:
\begin{equation}\label{eq:fhatcrit-S}
\boxed{ \hat{f} = \tilde{f} - \sum_{j \in J} \crr_{j}(f)
\quad\text{with}\quad
\tilde f= f- \langle \mu, f \rangle.}
\end{equation}
We shall consider, when well defined, for a sequence $\bF=(f_\ell,
\ell\in  \N)$ of measurable real-valued
functions defined on $S$, the quantities:
\begin{equation}\label{eq:scrit}
\scrit (\bF)=\scrit _1(\bF)+ 2\scrit _2(\bF), 
\end{equation}
where:
\begin{align}
\label{eq:S1-crit}
\scrit_1(\bF) 
&= \sum_{k\geq 0} 2^{-k} \langle \mu, \cp f_{k, k}^*
\rangle
= \sum_{k\geq 0} 2^{-k} \sum_{j\in J} \langle \mu,
\cp(\crr_{j}(f_k) \sot \overline{\crr}_{j}(f_k)) \rangle,\\
\label{eq:S2-crit}
\scrit_2(\bF) 
&=   \sum_{0\leq \ell <k} 2^{-(k+\ell)/2} \langle \mu, \cp
f_{k, \ell}^*  \rangle,
\end{align}
with, for $k, \ell\in \N$:
\begin{equation*}\label{eq:def-fkl*}
f_{k, \ell}^*= \sum_{j\in J}\,\, \theta_j^{\ell -k}\,  \crr_j
(f_k)\sot   \overline \crr_j (f_\ell).  
\end{equation*}
Notice  that $f_{k, \ell}^*=f_{\ell, k}^*$ and that $f_{k, \ell}^*$ is
real-valued as  $\overline {\theta_j^{\ell -k}\,  \crr_j (f_k)\otimes
\overline \crr_j (f_\ell)}= \theta_{j'}^{\ell -k}\,  \crr_{j'} (f_k)\otimes
\overline \crr_{j'} (f_\ell)$ for $j'$ such that $\alpha_{j'}=\overline
\alpha _j$ and thus $\crr_{j'}=\overline \crr_j$. 
\medskip

The technical proof of the next result is omitted as it is an adaptation of the
proof of Theorem \ref{cor:subcritical-L2} in the sub-critical space in
the same spirit as
\cite[Theorem~3.4]{BD} (critical case)  is an adaptation of the
proof of \cite[Theorem~3.1]{BD} (sub-critical case). The interested
reader can find the details in \cite{BD-arxiv}. 

\begin{theo}\label{cor:critical-L2}
Let $X$ be a BMC with kernel $\cp$ and initial distribution $\nu$ such
that Assumptions  \ref{hyp:Q} (with $k_0\in \N$), \ref{hyp:Q1} and
\ref{hyp:Q2} are in  force with  $\alpha = 1/\sqrt{2}$. We  have the
following convergence in distribution for all sequence $\bF=(f_\ell,
\ell\in   \N)$ bounded  in   $L^4(\mu)$   (that is $\sup_{\ell\in \N}
\norm{f_\ell}_{L^4(\mu)}<+\infty $): 
\[
n^{-1/2}N_{n, \emptyset}(\bF) \; \xrightarrow[n\rightarrow \infty
]{\text{(d)}} \; G, 
\]
where $G$ is centered Gaussian   random variable with 
variance $\scrit(\bF)$ given by \reff{eq:scrit}, which is well defined
and finite. 
\end{theo}

Notice   that  the   variance  $\scrit(\bF)$   already  appears   in  the
critical  pointwise  approach  case,   see  \cite[(20)  and  Theorem
3.4]{BD}. Then, arguing similarly as in \cite[Section 3.2]{BD}, we deduce that if
 Assumptions  \ref{hyp:Q} (with $k_0\in \N$), \ref{hyp:Q1} and
\ref{hyp:Q2} are in  force with  $\alpha = 1/\sqrt{2}$, then  for $f\in L^{4}(\mu)$,  we have
the following convergence in distribution:
\begin{equation}
  \label{eq:criticalGT} 
(n|\GG_{n}|)^{-1/2}M_{\GG_{n}}(\tilde{f}) \; \xrightarrow[n\rightarrow
\infty ]{\text{(d)}} \; G_{1},
\quad \text{and} \quad
  (n|\TT_{n}|)^{-1/2}M_{\TT_{n}}(\tilde{f}) \; \xrightarrow[n\rightarrow
  \infty ]{\text{(d)}} \; G_{2},  
\end{equation}
where $G_{1}$ and $G_2$ are  centered Gaussian  random variables with
respective  
variances   $\scrit_{\GG}(f)=\scrit(\bF)$, with  $\bF = (f,0,0, \ldots)$,
and $\scrit_{\TT}(f)=\scrit(\bff)/2$  with $ \bff=(f,f, \ldots)$,
given in  \cite[Corollary 3.6]{BD}  which are  well defined
and finite.

\subsection{The super-critical case $2\alpha^2>1$}
\label{sec:main-res-super}

We  consider the  super-critical  case $\alpha\in  (1/\sqrt{2}, 1)$.
This case is very similar to the super-critical case in the pointwise
approach, see 
\cite[Section~3.3]{BD}. So we only mention  the most interesting
results without proof.  The interested
reader can find the details in \cite{BD-arxiv}.
\medskip

We shall assume that Assumptions \ref{hyp:Q} (ii) and \ref{hyp:Q2} hold.
In  particular we  do not  assume Assumption  \ref{hyp:Q2} (i).   Recall
\reff{eq:L2-erg-crit}           with           the           eigenvalues
$\{\alpha_j=\theta_j \alpha, j\in J\} $  of $\cq$, with modulus equal to
$\alpha$ (\textit{i.e.}   $|\theta_j|=1$) and the projector  $\crr_j$ on
the eigen-space  associated to  eigenvalue $\alpha_j$.  Recall  that the
sequence $(\beta_n, n\in \N)$ in  Assumption \ref{hyp:Q2} can (and will)
be chosen non-increasing and bounded from above by 1.  We shall consider
the     filtration      $\ch=(\ch_n,     n\in     )$      defined     by
$\ch_{n}   =  \sigma(X_{i},i\in\TT_{n})$.    The  next   lemma  exhibits
martingales related to the projector $\crr_j$.

\begin{lem}\label{lem:martingale}
Let $X$ be a BMC with kernel $\cp$ and initial distribution
$\nu$ such  that Assumptions \ref{hyp:Q} (ii) and \ref{hyp:Q2} are in
force with
$\alpha\in (1/\sqrt{2}, 1)$ in
\reff{eq:L2-erg-crit}. Then, for all $j\in J$ and $f\in L^{2}(\mu)$, the sequence
$M_{j}(f)=(M_{n,j}(f), n\in \N)$, with
\begin{equation*}
M_{n,j}(f) = (2 \alpha_j) ^{-n} \, M_{\G_n} (\crr_j(f)), 
\end{equation*}
is a $\ch$-martingale   which converges a.s. and in $L^{2}(\nu)$ to a  random
variable, say  $M_{\infty,j}(f)$.
\end{lem}

The next result corresponds to \cite[Corollary~3.13]{BD} in the
pointwise approach. 
\begin{cor}
   \label{cor:super-crit-d1-2}
Let $X$ be a BMC with kernel $\cp$ and initial distribution
$\nu$ such   that Assumptions \ref{hyp:Q} (ii) and  \ref{hyp:Q2} are in
force  with  $\alpha \in (1/\sqrt{2},
1) $ in  \reff{eq:L2-erg-crit}. Assume $\alpha$
is the only eigen-value of $\cq$ with modulus equal to $\alpha$ (and
thus $J$ is
reduced to a singleton, say $\{j_0\}$), then we have for $f\in L^{2}(\mu)$:
\begin{equation*}
  (2\alpha)^{-n}M_{\GG_{n}}(\tilde{f}) \inP M_{\infty}(f)
  \quad \text{and} \quad
  (2\alpha)^{-n}M_{\TT_{n}}(\tilde{f}) \inP \frac{2\alpha}{2\alpha -
    1}M_{\infty, j_0}(f) ,
\end{equation*}
where $M_{\infty, j_0}(f)$ is the random variable defined in Lemma
\ref{lem:martingale}. 
\end{cor}

\section{Application to the study of symmetric BAR}\label{sec:sym-bar-M}

\subsection{Symmetric BAR}\label{sec:sym-bar}
We consider a particular case from \cite{CS86} of the real-valued  bifurcating autoregressive process (BAR) from Example  \ref{ex:valide}. We keep the  same notations. Let $a\in (-1, 1)$ and assume that $a=a_{0} = a_{1}$, $b_{0} = b_{1} = 0$ and $\rho = 0$. In this particular case the BAR has symmetric kernel as:
\[
\cp(x, dy, dz)=\cq(x, dy) \cq(x, dz).
\]
We have $\cq f(x)=\E[f(ax +\sigma G)]$ and more generally $\cq^n f(x)=\E\left[f\left(a^n x + \sqrt{1- a^{2n}} \sigma_a G\right)\right]$, where $G$ is a standard $\cn(0, 1)$ Gaussian random variable and $\sigma_a=\sigma (1- a^2)^{-1/2}$. The kernel $\Qq$ admits a unique invariant probability measure $\mu$, which is $\cn(0, \sigma_a^2)$ and whose  density, still denoted by $\mu$, with respect to the Lebesgue measure  is given by:
\begin{equation*}\label{eq:mu-SBAR}
\mu(x) = \frac{\sqrt{1 - a^{2}}}{\sqrt{2\pi\sigma^2}} \exp\left(-\frac{(1 - a^{2})x^{2}}{2 \sigma^{2}}\right). 
\end{equation*}
The density $p$ (resp. $q$) of the kernel $\cp$ (resp. $\cq$) with respect to $\mu^{\otimes 2}$ (resp. $\mu$) are given by:
\[
p(x,y,z) = q(x,y)q(x,z)
\]
and
\[
q(x,y) = \frac{1}{\sqrt{1 - a^{2}}} \exp\left(-\frac{(y-ax)^{2}}{2\sigma^{2}} + \frac{(1-a^{2})y^{2}}{2\sigma^{2}}\right) = \frac{1}{\sqrt{1 - a^{2}}} \expp{- (a^2 y^2 + a^2 x^2 -2axy )/ 2\sigma^ 2}.
\]
Notice that $q$ is symmetric.  The operator $\Qq$ (in $L^2(\mu)$) is a symmetric integral Hilbert-Schmidt operator whose eigenvalues are given by $\sigma_{p}(\Qq) = (a^{n}, n \in \NN)$, their algebraic multiplicity is one and the corresponding eigen-functions $(\bar{g}_{n}(x), n \in \NN)$ are defined for $n\in \N$ by :
\begin{equation*}\label{eq:HermiteBAR}
\bar{g}_{n}(x) = g_{n}\left(\sigma_a^{-1} \, x\right),
\end{equation*}
where $g_{n}$ is the Hermite polynomial of degree $n$ ($g_0=1$ and $g_1(x)=x$). Let  $\crr$ be the  orthogonal projection on the vector space generated by $\bar{g}_{1}$, that is $\crr f= \langle \mu, f\bar g_1 \rangle\, \bar g_1$ or equivalently, for $x\in \R$:
\begin{equation}\label{eq:R-symBAR}
\crr f(x)=\sigma_a^{-1}\, x \, \E\left[G f(\sigma_aG)\right].
\end{equation}

\medskip

Recall $\bH$ defined \eqref{eq:def-bH}. It is not difficult to check that:
\[
\bH(x)=(1-a^4)^{-1/4}\exp\left(\frac{a^2(1-a^2)}{1+a^2} \, \frac{x^2}{2\sigma^2}\right) \quad\text{for $x\in \R$},
\]
and $\bH \in L^2(\mu)$ (that  is $\int_{\R^2}  q(x,y)^2  \,  \mu(x)  \mu(y)\, dx  dy<+\infty  $).   Using elementary    computations,   it    is    possible    to   check    that $\cq  \bH\in  L^4(\mu)$   if  and  only  if   $|a|<  3^{-1/4}$  (whereas $\bH\in  L^4(\mu)$  if  and  only   if  $|a|<3^{-1/2}$).   As  $\cp$  is symmetric,   we   get $\cp(\bH\otimes^2)\leq  (\cq\bH)^2$   and   thus \eqref{eq:Hil-SchP-1}  holds for  $|a|<  3^{-1/4}$. We  also get,  using Cauchy-Schwartz inequality, that $\normm{\cp(f\sot     \bH)}=\normm{(\cq    f) (\cq    \bH)}  \leq \norm{f}_{L^4(\mu)}\,\norm{\cq     (\bH)}_{L^4(\mu)}$,      and     thus \eqref{eq:Hil-SchP-3}  holds  for   $|a|<  3^{-1/4}$.   Some  elementary computations   give   that    \eqref{eq:Hil-SchP-2}   also   holds   for $|a|\leq   0.724$   (but   \eqref{eq:Hil-SchP-2}  fails   for   $|a|\geq 0.725$). (Notice that $2^{-1/2}< 0.724< 3^{-1/4}$.)  As a consequence of Remark      \ref{rem:hypP},      if      $|a|\leq      0.724$,      then \eqref{eq:P-L2-1}-\eqref{eq:P-L2-2}  are  satisfied   and  thus  (i)  of Assumption \ref{hyp:Q} holds.

Notice that $\nu \cq^k$ is the probability distribution of $a^k X_\emptyset + \sigma_a \sqrt{1- a^{2k}}\,  G$, with $G$ a $\cn(0, 1)$ random variable independent of $X_\emptyset$. So property (ii) of Assumption \ref{hyp:Q} holds in particular if $\nu$ has compact support (with $k_0=1$) or if $\nu$ has a density with respect to the Lebesgue measure, which we still denote by $\nu$, such that $\norm{\nu/\mu}_\infty $ is finite (with $k_0\in \N$). Notice that if $\nu$ is the probability distribution of $\cn(0, \rho_0^2)$, then  $\rho_0> \sigma_a$ (resp. $\rho_0\leq  \sigma_a$) implies that   (ii) of Assumption \ref{hyp:Q} fails (resp. is satisfied). 
\medskip

Using that $(\bar g_n/\sqrt{n!}, \, n\in \N)$ is an orthonormal basis of $L^2(\mu)$ and  Parseval identity, it  is easy to check  that Assumption \ref{hyp:Q2}  holds  with  $J=\{j_0\}$,  $\alpha_{j_0} =  \alpha  =  a$, $\beta_n = a^{n}$ and $\crr_{j_0}=\crr$.  

\subsection{Numerical studies: illustration of phase transitions for the fluctuations} \label{sec:simul}

We consider the symmetric BAR  model from Section \ref{sec:sym-bar} with
$a=\alpha\in (0, 1)$. Recall $\alpha$ is an eigenvalue with multiplicity
one,  and  we  denote  by   $\crr$  the  orthogonal  projection  on  the
one-dimensional  eigenspace associated  to $\alpha$.  The expression  of
$\crr$ is given in \eqref{eq:R-symBAR}.

\medskip

In order to illustrate the effects  of the geometric rate of convergence
$\alpha$ on the fluctuations, we  plot for $\A_n\in \{\G_n, \T_n\}$, the
slope,     say     $b_{\alpha,n}$,     of    the     regression     line
$\log(\Var(|\A_{n}|^{-1}M_{\A_{n}}(f)))          $         \emph{versus}
$ \log(|\A_{n}|)$  as a  function of the  geometric rate  of convergence
$\alpha$. In  the classical cases  (e.g. Markov chains), the  points are
expected to be distributed around the horizontal line $y = -1$.  For $n$
large, we have  $\log(|\A_{n}|)\simeq n \log(2)$ and,  for the symmetric
BAR model, convergences  \eqref{eq:subcriticalGT} for $\alpha<1/\sqrt{2}$,
\eqref{eq:criticalGT}     for     $\alpha=1/\sqrt{2}$,    and     Corollary
\ref{cor:super-crit-d1-2}    for    $\alpha>1/\sqrt{2}$   yields    that
$b_{\alpha,n}         \simeq         h_{1}(\alpha)        $         with
$h_1(\alpha)=\log(\alpha^2\vee 2^{-1})/\log(2)$ as  soon as the limiting
Gaussian  random  variable   in  \eqref{eq:subcriticalGT} and
\eqref{eq:criticalGT}  or     $M_\infty      (f)$     in      Corollary
\ref{cor:super-crit-d1-2} is non-zero.  \medskip

\begin{figure}
	\centering
	\begin{subfigure}{0.45\textwidth} 
		\includegraphics[width=\textwidth]{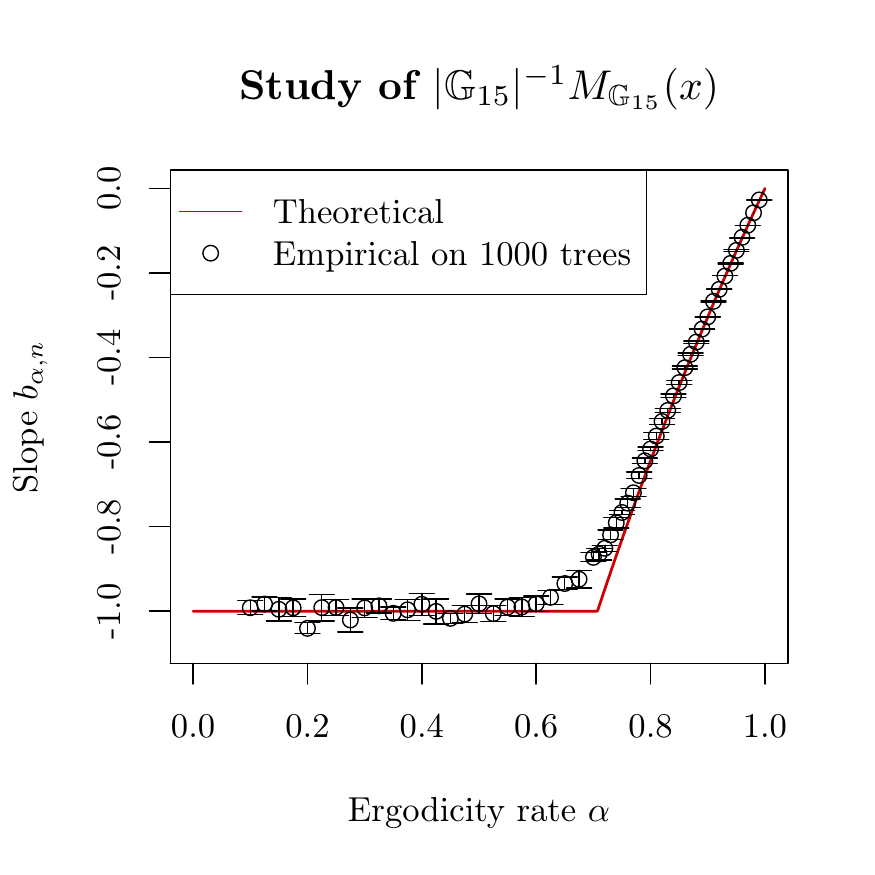}
	\end{subfigure}
	\begin{subfigure}{0.45\textwidth} 
		\includegraphics[width=\textwidth]{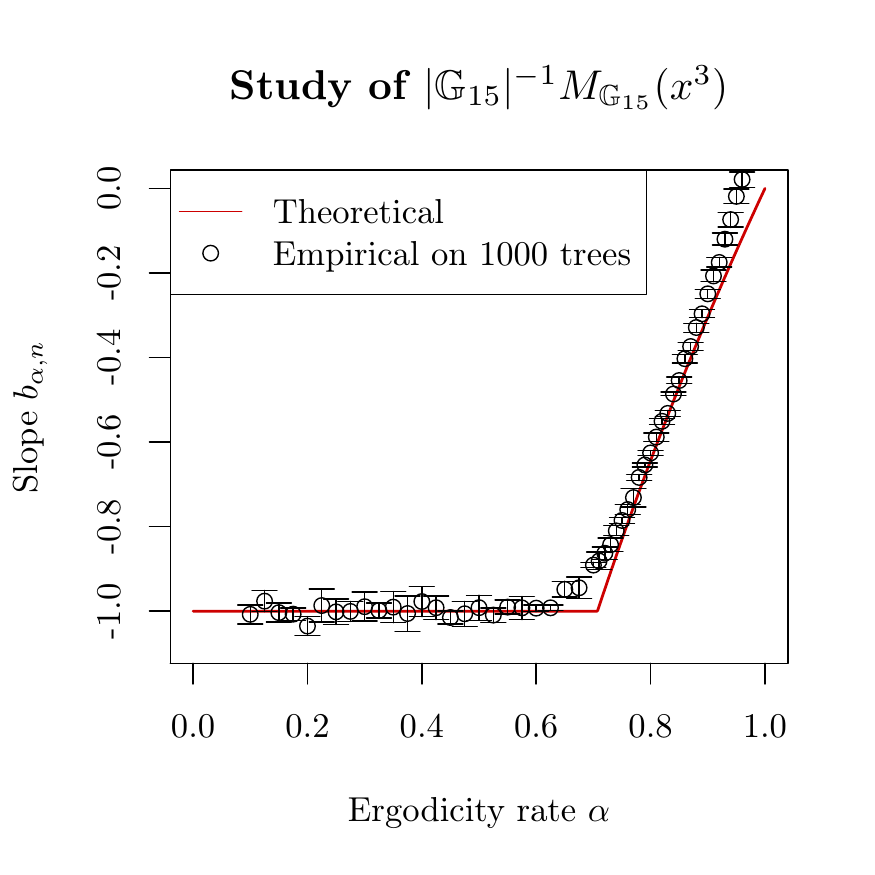}
	\end{subfigure}
	\centering
	\begin{subfigure}{0.45\textwidth} 
		\includegraphics[width=\textwidth]{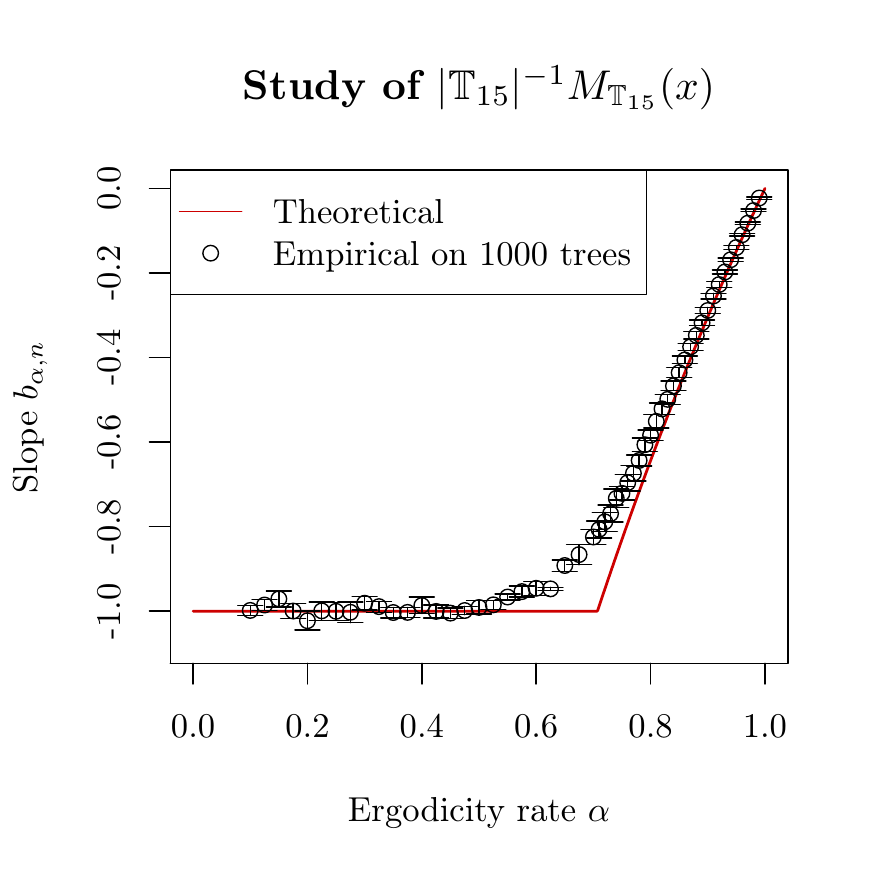}
	\end{subfigure}
	\begin{subfigure}{0.45\textwidth} 
		\includegraphics[width=\textwidth]{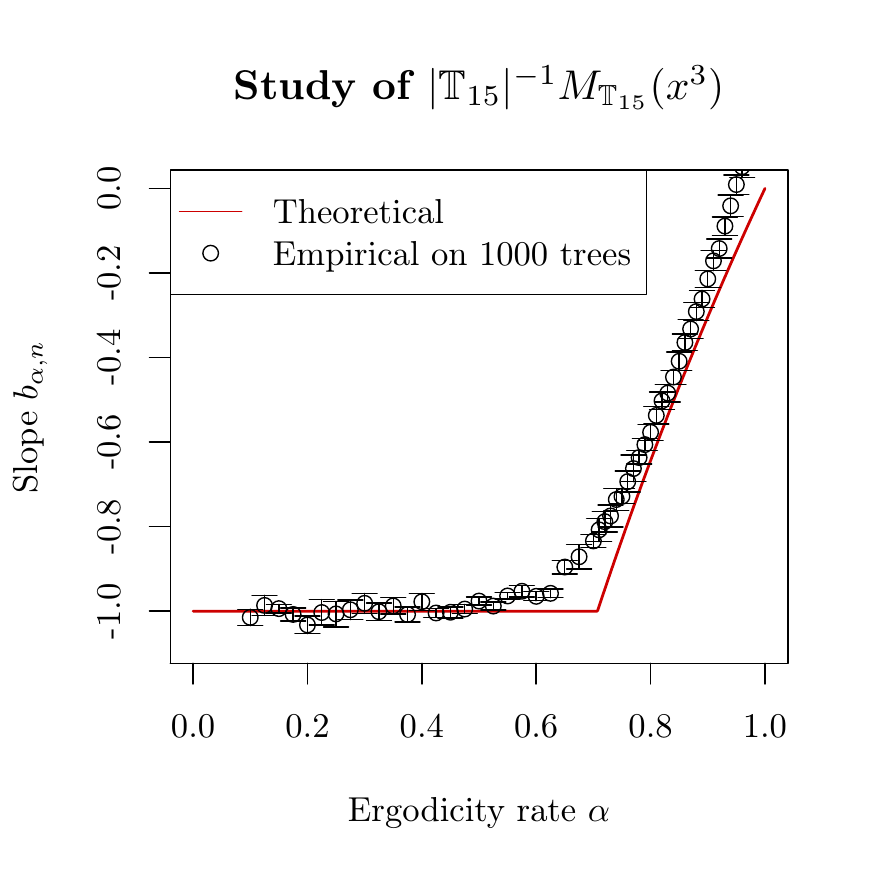}
	\end{subfigure}
        \caption{Slope  $b_{\alpha, n}$  (empirical mean  and
          confidence     interval in  black)     of    the     regression     line
          $\log(\Var(|\A_{n}|^{-1}M_{\A_{n}}(f)))     $    \emph{versus}
          $ \log(|\A_{n}|)$ as a function  of the geometric ergodic rate
          $\alpha$, for  $n=15$, $\A_{n}  \in \{\GG_{n},  \TT_{n}\}$ and
          $f(x) =  x^p$ with  $p\in \{1,  3\}$.  In  this case,  we have
          $\crr(f) \neq 0$,  where $\crr$ is the  projector defined from
          formula \eqref{eq:R-symBAR}.   One can see that  the empirical
          curve (in  black) is close to the  graph  (in red)  of the
          function  $h_1(\alpha) =  \log(\alpha^{2} \vee  2^{-1})/\log(2)$
           for $\alpha\in (0, 1)$.} \label{fig:I-phase-clt}
\end{figure}

\begin{figure}
	\centering
	\begin{subfigure}{0.45\textwidth} 
		\includegraphics[width=\textwidth]{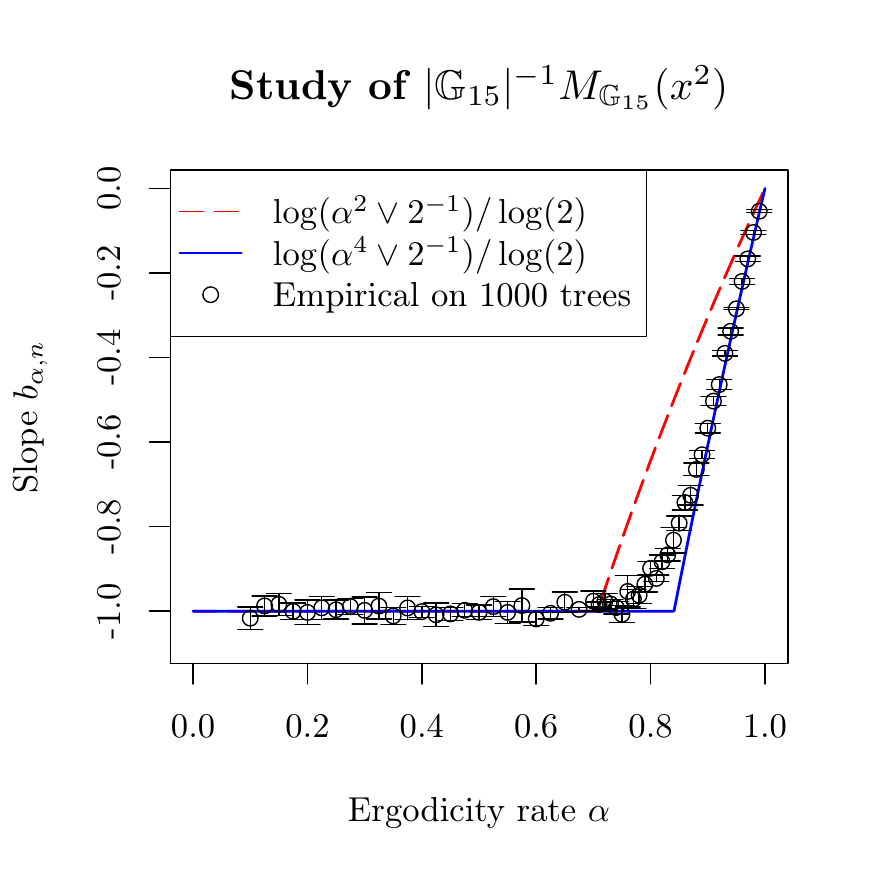}
	\end{subfigure}
	\begin{subfigure}{0.45\textwidth} 
		\includegraphics[width=\textwidth]{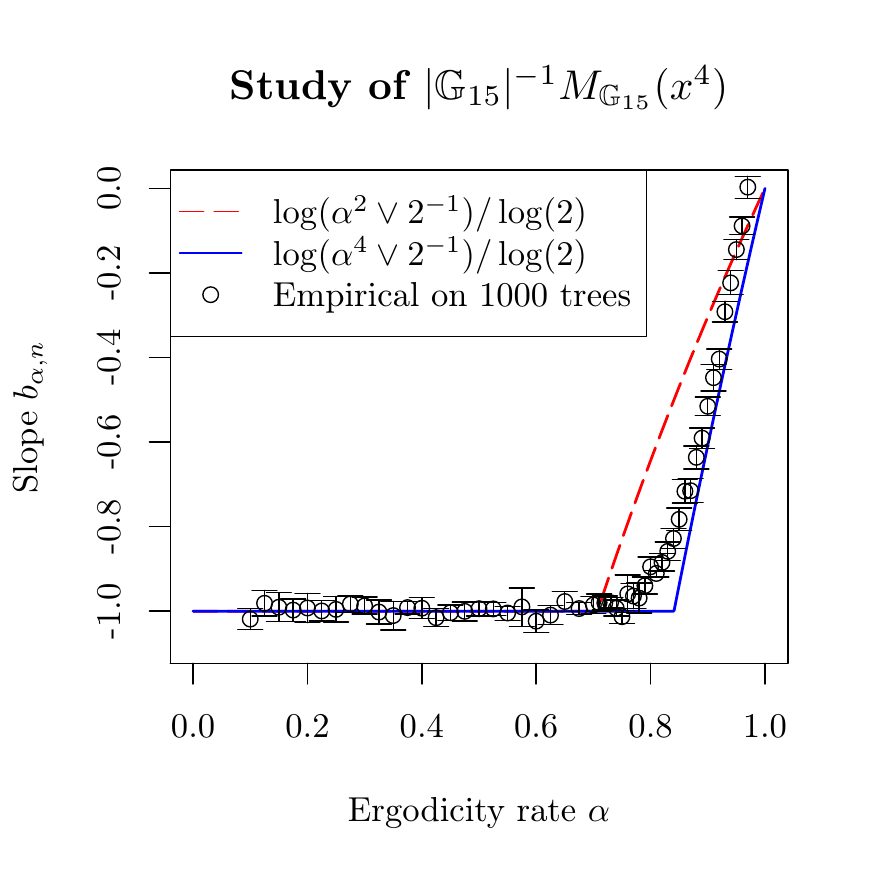}
	\end{subfigure}
	\centering
	\begin{subfigure}{0.45\textwidth} 
		\includegraphics[width=\textwidth]{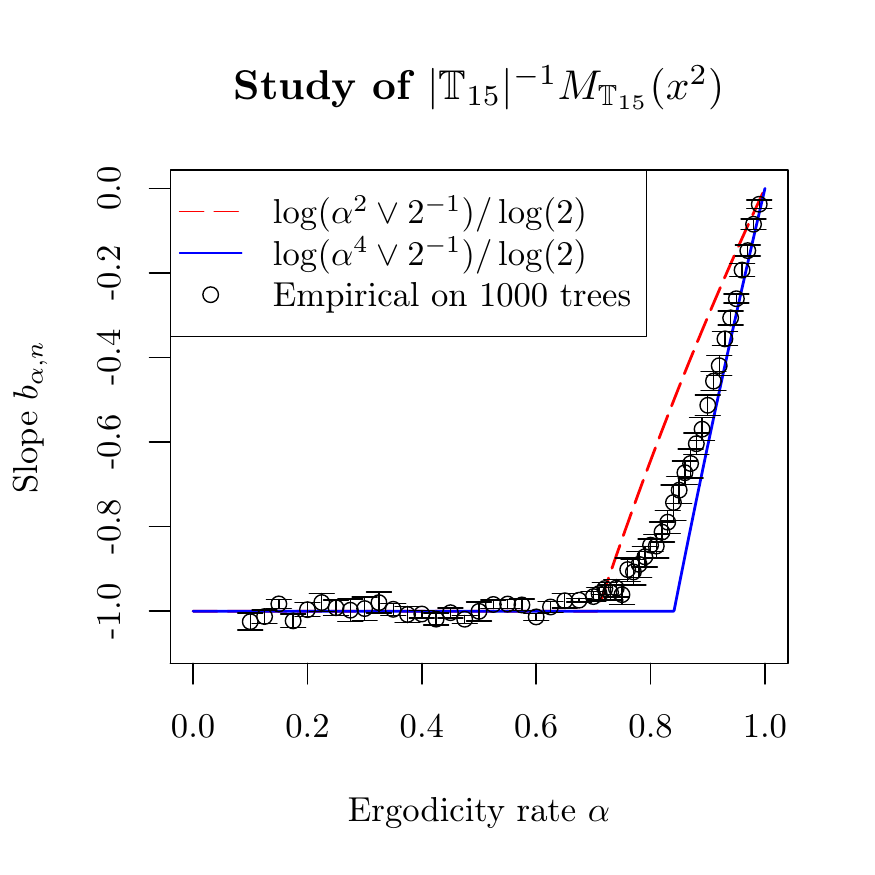}
	\end{subfigure}
	\begin{subfigure}{0.45\textwidth} 
		\includegraphics[width=\textwidth]{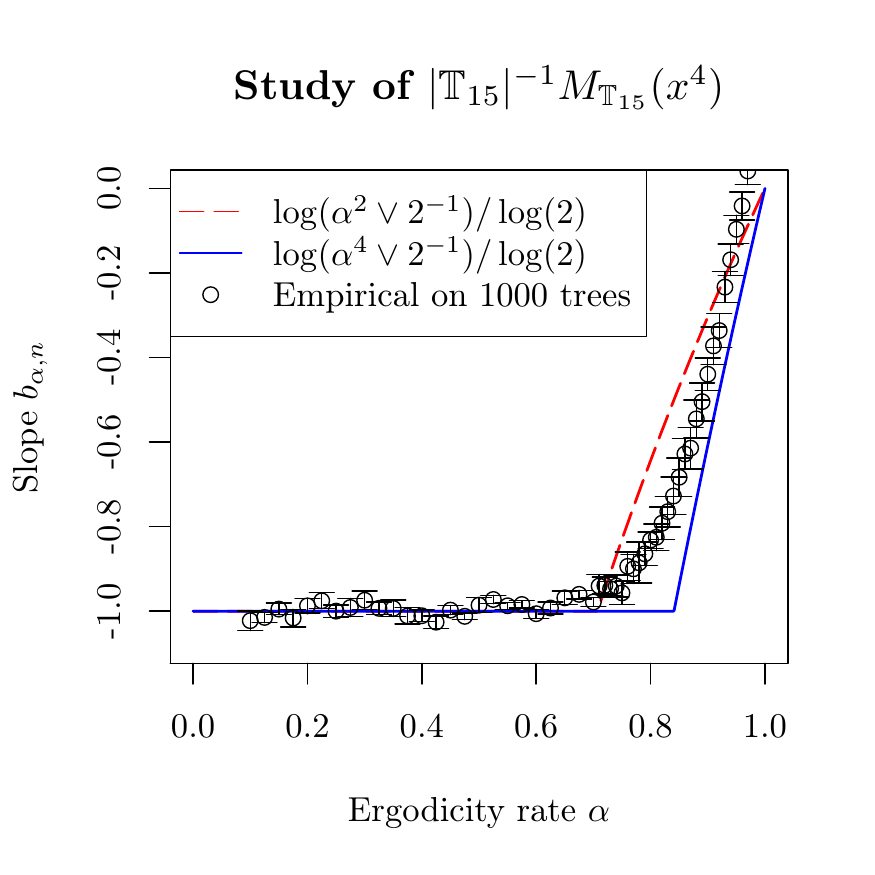}
	\end{subfigure}
\caption{Slope  $b_{\alpha, n}$  (empirical mean  and
          confidence     interval in  black)     of    the     regression     line
          $\log(\Var(|\A_{n}|^{-1}M_{\A_{n}}(f)))     $    \emph{versus}
          $ \log(|\A_{n}|)$ as a function  of the geometric ergodic rate
          $\alpha$, for  $n=15$, $\A_{n}  \in \{\GG_{n},  \TT_{n}\}$ and
          $f(x) =  x^p$ with  $p\in \{2,  4\}$.  In  this case,  we have
          $\crr(f) = 0$,  where $\crr$ is the  projector defined from
          formula \eqref{eq:R-symBAR}.   One can see that  the empirical
          curve (in  black) does not agree with the  graph  (dash line in red)  of the
          function  $h_1(\alpha) =  \log(\alpha^{2} \vee
          2^{-1})/\log(2)$ for $2\alpha^2>1$;
          but it is close to  the graph  (in blue)  of the
          function  $h_2(\alpha) =  \log(\alpha^{4} \vee
          2^{-1})/\log(2)$ for $\alpha\in (0, 1)$.} \label{fig:I-phase-cltii-bis}
\end{figure}

For  our  illustrations, we  consider  the  empirical moments  of  order
$p\in \{1, \ldots, 4\}$, that is we  use the functions $f(x) = x^p$.  As
we     can     see      in     Figures     \ref{fig:I-phase-clt}     and
\ref{fig:I-phase-cltii-bis},  these curves  present  two  trends with  a
phase   transition   around  the   rate   $\alpha   =  1/\sqrt{2}$   for
$p  \in \{1,3\}$  and  around  the rate  $\alpha^{2}  = 1/\sqrt{2}$  for
$p \in \{2,4\}$. For convergence  rates $\alpha \in (0,1/\sqrt{2})$, the
trend  is similar  to  that  of classic  cases.   For convergence  rates
$\alpha  \in  (1/\sqrt{2},1)$, the  trend  differs  to that  of  classic
cases. One can observe that  the slope $b_{\alpha,n}$ increases with the
value of geometric convergence rate  $\alpha$.  We also observe that for
$\alpha >  1/ \sqrt{2}$, the empirical  curves agrees with the  graph of
$h_1(\alpha)=\log(\alpha^2\vee 2^{-1})/\log(2)$ for  $f(x) = x^{p}$ when
$p$  is odd,  see Figure  \ref{fig:I-phase-clt}. However,  the empirical
curves does not  agree with the graph  of $h_1$ for $f(x)  = x^{p}$ when
$p$ is even, see Figure  \ref{fig:I-phase-cltii-bis}, but it agrees with
the             graph             of            the             function
$h_2(\alpha)=\log(\alpha^4\vee  2^{-1})/\log(2)$.  This  is  due to  the
fact that for $p$ even, the function $f(x)=x^p$ belongs to the kernel of
the projector $\crr$ (which  is clear from formula \eqref{eq:R-symBAR}),
and thus $M_\infty (f)=0$.  In fact, in those two cases, one should take
into account  the projection on  the eigenspace associated to  the third
eigenvalue,   which    in   this    particular   case   is    equal   to
$\alpha^2$.  Intuitively,  this  indeed  give a  rate  of  order  $h_2$.
Therefore, the normalization given for  $f(x)=x^p$ when $p$ even, is not
correct.

\section{Proof of Theorem \ref{cor:subcritical-L2}} 
\label{sec:proof-sub-cor-L2}

In the following proofs, we will  denote by  $C$ any unimportant finite  constant which may  vary from  line to line  (in particular $C$ does not  depend on $n$ nor on $\bF$).


Let $(p_n, n\in \N)$ be a non-decreasing sequence of elements of $\N^*$
such that, for all $\lambda>0$:
\begin{equation}
   \label{eq:def-pn}
p_n< n, \quad \lim_{n\rightarrow \infty } p_n/n=1
\quad\text{and}\quad \lim_{n\rightarrow \infty } n-p_n - \lambda \log(n)=+\infty .
\end{equation}
When there is no ambiguity, we write $p$ for $p_n$. 
\medskip

Let $i,j\in \T$. We write $i\preccurlyeq  j$ if $j\in i\T$. We denote by
$i\wedge j$  the most recent  common ancestor of  $i$ and $j$,  which is
defined  as   the  only   $u\in  \T$   such  that   if  $v\in   \T$  and
$ v\preccurlyeq i$, $v \preccurlyeq j$  then $v \preccurlyeq u$. We also
define the lexicographic order $i\leq j$ if either $i \preccurlyeq j$ or
$v0  \preccurlyeq i$  and $v1  \preccurlyeq j$  for $v=i\wedge  j$.  Let
$X=(X_i, i\in  \T)$ be  a $BMC$  with kernel  $\cp$ and  initial measure
$\nu$. For $i\in \T$, we define the $\sigma$-field:
\begin{equation*}\label{eq:field-Fi}
\cf_{i}=\{X_u; u\in \T \text{ such that  $u\leq i$}\}.
\end{equation*}
By construction,  the $\sigma$-fields $(\cf_{i}; \, i\in \T)$ are nested
as $\cf_{i}\subset \cf_{j} $ for $i\leq  j$.
\medskip

We define for $n\in \N$, $i\in \G_{n-p_n}$ and $\bF\in F^\N$ the
martingale increments:
\begin{equation}
   \label{eq:def-DiF}
\Delta_{n,i}(\bF)= N_{n,i}(\bF) - \E\left[N_{n,i}(\bF) |\,
  \cf_i\right]
\quad\text{and}\quad \Delta_n(\bF)=\sum_{i\in \G_{n-p_n}} \Delta_{n,i}(\bF).
\end{equation}
Thanks to \reff{eq:def-NiF}, we have:
\[
\sum_{i\in \G_{n-p_n}} N_{n, i}(\bF) 
= |\G_n|^{-1/2} \sum_{\ell=0}^{p_n}  M_{\G_{n-\ell}} (\tilde f_\ell)
= |\G_n|^{-1/2} \sum_{k=n-p_n}^{n}  M_{\G_{k}} (\tilde f_{n-k}).
\]
Using the branching Markov property, and \eqref{eq:def-NiF}, we get for
$i\in \G_{n-p_n}$:
\[
\E\left[N_{n,i}(\bF) |\,
  \cf_i\right]
=\E\left[N_{n,i}(\bF) |\,
  X_i\right]
= |\G_n|^{-1/2} \sum_{\ell=0}^{p_n}
\E_{X_i}\left[M_{\G_{p_n-\ell}}(\tilde f_\ell)\right].
\] 
We deduce from \reff{eq:nof-D} with $k=n-p_n$ that:
\begin{equation}
   \label{eq:N=D+R}
N_{n, \emptyset}(\bF)
= \Delta_n(\bF) + R_0(n)+R_1(n),
\end{equation}
with
\begin{equation}
   \label{eq:reste01}
R_0(n)= |\G_n|^{-1/2}\, \sum_{k=0}^{n-p_n-1}
  M_{\G_k}(\tilde f_{n-k})
\quad\text{and}\quad
R_1(n)= \sum_{i\in \G_{n-p_n}}\E\left[N_{n,i}(\bF) |\,
  \cf_i\right].
\end{equation}

We first state a very useful Lemma which holds in sub-critical, critical
and super-critical cases.  
\begin{lem}
   \label{lem:L2MG}
   Let $X$  be a BMC  with kernel  $\cp$ and initial  distribution $\nu$
   such that (ii)  from Assumption \ref{hyp:Q} (with $k_0\in  \N$) is in
   force. There exists a finite constant $C$, such that for all $f\in
   \cb_+(S)$ all $n\geq k_0$, we have: 
\begin{equation}
   \label{eq:EGn-nu0}
|\G_n|^{-1} \E[M_{\G_n} (f)]\leq  C  \norm{f}_{L^1(\mu)}
\quad\text{and}\quad
|\G_n|^{-1} \E\left[M_{\G_n} (f)^2\right]\leq  C  \sum_{k=0}^{n}
2^{k} \normm{\cq^kf}^2. 
\end{equation}
 \end{lem}
\begin{proof}
   Using the first moment formula \reff{eq:Q1},  (ii) from Assumption
   \ref{hyp:Q} and the fact that $\mu$ is invariant for $\cq$, we get that:
\[
|\G_n|^{-1} \E[M_{\G_n} (f)]
=\langle \nu, \cq^n f \rangle
\leq  \norm{\nu_0}_\infty\, \langle \mu, \cq^{n-k_0}f  \rangle
=  \norm{\nu_0}_\infty\, \langle \mu, f  \rangle.
\]

We also have:
\begin{align*}
   |\G_n|^{-1} \E\left[M_{\G_n} (f)^2\right]
&= \langle \nu, \cq^n (f^2) \rangle
+ \sum_{k=0}^{n-1} 2^{k} \, \langle \nu, \cq^{n-k-1} \left( \cp(\cq^k f
  \otimes^2\right) \rangle\\
&\leq   \langle \nu, \cq^n (f^2) \rangle
+ \sum_{k=0}^{n-1} 2^{k} \, \langle \nu, \cq^{n-k} \left( (\cq^k f)^2
\right) \rangle\\
&\leq   \langle \nu, \cq^n (f^2) \rangle
+ \sum_{k=0}^{n-k_0} 2^{k} \, \langle \nu, \cq^{n-k} \left( (\cq^k f)^2
\right) \rangle
+ \sum_{k=n-k_0+1}^{n-1} 2^{k} \, \langle \nu, \cq^{k_0} \left(
  (\cq^{n-k_0} f)^2 
\right) \rangle\\  
&\leq  C 
\sum_{k=0}^{n-k_0} 2^{k} \normm{\cq^kf}^2,
\end{align*}
where we used  the second moment formula \reff{eq:Q2}  for the equality,
\reff{eq:majo-pfxf} for the first  inequality, Jensen inequality for the
second, and (ii)  from Assumption \ref{hyp:Q} and the fact  that $\mu$ is
invariant for $\cq$ for the last.
\end{proof}

We set for $k\in \N^*$: 
\begin{equation}\label{eq:def-ck}
c_{k}(\bF)=\sup_{n\in \N}\norm{f_n}_{L^{k}(\mu)} \quad \text{and}\quad q_k(\bF)=\sup_{n\in \N}\norm{\cq( f_n^k)}_{\infty }^{1/k}.
\end{equation}
We will  denote by  $C$ any unimportant finite  constant which may  vary
from  line to line  (but in particular $C$ does not  depend on $n$ nor
on $\bF$, but may depends on $k_0$ and $\norm{\nu_0}_\infty   $).
\medskip

\begin{rem}\label{rem:Nn=Nnk0}
Recall $k_{0}$ given in Assumption \ref{hyp:Q} (ii). Let $\bF = (f_{\ell}, \ell \in \NN)$ be a bounded sequence in $L^{4}(\mu)$.  We have
\begin{equation}\label{eq:Nn0=Nn0k0+R}
N_{n,\emptyset}(\bF) = N_{n,\emptyset}^{[k_{0}]}(\bF) \, + \, |\GG_{n}|^{-1/2} \sum_{\ell = 0}^{k_{0} - 1} M_{\GG_{\ell}}(\tilde{f}_{n-\ell}),
\end{equation}
where we set:
\begin{equation}
   \label{eq:def-NOf2}
\boxed{ N_{n, \emptyset}^{[k_0]}(\bF) = |\G_n|^{-1/2
  }\sum_{\ell=0}^{n-k_0} M_{\G_{n-\ell}}(\tilde f_{\ell}) .} 
\end{equation}
Using the Cauchy-Schwartz inequality, we get
\begin{equation}\label{eq:reste}
|\GG_{n}|^{-1/2} |\sum_{\ell = 0}^{k_{0} - 1}
M_{\GG_{\ell}}(\tilde{f}_{n-\ell})| \leq C c_{2}(\bF) |\GG_{n}|^{-1/2}
\, + \, |\GG_{n}|^{-1/2} \sum_{\ell = 0}^{k_{0}-1}
M_{\GG_{\ell}}(|f_{n-\ell}|)  .
\end{equation}   
Since the sequence $\bF$ is bounded in $L^{4}(\mu)$ and since $k_{0}$ is finite, we have, for all $\ell \in \{0, \ldots, k_{0}-1\}$, $\lim_{n \rightarrow \infty} |\GG_{n}|^{-1/2} M_{\GG_{\ell}}(|f_{n-\ell}|) = 0$ a.s. and then that (used \eqref{eq:reste})
\begin{equation*}
\lim_{n \rightarrow \infty} |\GG_{n}|^{-1/2} |\sum_{\ell = 0}^{k_{0} - 1} M_{\GG_{\ell}}(\tilde{f}_{n-\ell})| = 0 \quad \text{a.s.}
\end{equation*}
Therefore, from \eqref{eq:Nn0=Nn0k0+R}, the study of $N_{n,\emptyset}(\bF)$ is reduced to that of $N^{[k_{0}]}_{n,\emptyset}(\bF)$.   
\end{rem}


Recall   $(p_n,  n\in   \N)$  is   such  that \reff{eq:def-pn}  holds.  Assume  that  $n$  is  large  enough  so  that
$n-p_n -1\geq k_0$. We have:
\[
N^{[k_0]}_{n, \emptyset}(\bF) = \Delta_n(\bF) + R_0^{k_0} (n)+R_1(n),
\]
where $\Delta_n(\bF)$ and $R_1(n)$ are defined in \reff{eq:def-DiF} and \reff{eq:reste01}, and :
\[
R_0^{k_0} (n)= |\G_n|^{-1/2}\, \sum_{k=k_0}^{n-p_n-1} M_{\G_k}(\tilde f_{n-k}).
\]

\begin{lem}\label{lem:cvR0-L2}
Under the assumptions of Theorem \ref{cor:subcritical-L2}, we have the following convergence: 
\[
\lim_{n\rightarrow \infty } \E[R_0^{k_0} (n)^2] =0.
\]
\end{lem}
\begin{proof}
Assume $n-p\geq k_0$. We write:
\[
R_0^{k_0} (n)= |\G_n|^{-1/2}\, \sum_{k=k_0}^{n-p-1}\sum_{i\in \G_{k_0}}  M_{i\G_{k-k_0}}(\tilde f_{n-k}).
\]
We have that $
\sum_{i\in \G_{k_0}} \E[ M_{i\G_{k-k_0}}(\tilde f_{n-k})^2] 
=
\E[M_{\G_{k_0}}(h_{k,n})]$,
where:
\[
h_{k,n}(x)=\E_x[ M_{\G_{k-k_0}}(\tilde f_{n-k})^2].
\] 
We deduce from  (ii) from Assumption \ref{hyp:Q}, see \reff{eq:EGn-nu0}, that $
\E[M_{\G_{k_0}}(h_{k,n})]\leq  C \langle \mu,  h_{k, n}\rangle$.
We have also that:
\begin{align*}
\langle \mu,  h_{k, n}\rangle 
= \E_\mu[ M_{\G_{k-k_0}}(\tilde f_{n-k})^2] 
& \leq  C \, 2^k \sum_{\ell=0}^k 2^\ell \normm{\cq^\ell \tilde
  f_{n-k}}^2
 \leq  C \, 2^k\, c_{2}^2(\bF) \sum_{\ell=0}^{k}
  2^{\ell}\alpha^{2\ell}
\leq  C 2^{k}\,c_{2}^2(\bF) ,
\end{align*}
where we used \reff{eq:EGn-nu0} for the first inequality (notice one can
take $k_0=0$ in this case as we consider the expectation $\E_\mu$), 
\reff{eq:L2-erg}  in the second, and $2\alpha^2<1$ in the last. We deduce
that:

\begin{equation}\label{eq:IRk0nL2sub}
\E[R_0^{k_0}(n)^2]^{1/2} \leq  |\G_n|^{-1/2}\, \sum_{k=k_0}^{n-p-1} \left(2^{k_0} \E\left[M_{\G_{k_0}}(h_{k,n})\right]\right)^{1/2} \leq  C\, 2^{- p/2} c_{2}(\bF),
\end{equation}
where we used that the sequence $\bF$ is bounded in $L^2(\mu)$. Use that $\lim_{n\rightarrow \infty } p=\infty $ to conclude. 
\end{proof}

We have the following lemma. 
\begin{lem}
   \label{lem:cvR1-L2}
Under the assumptions of Theorem \ref{cor:subcritical-L2}, we have the
following convergence: 
\[
\lim_{n\rightarrow \infty } \E\left[R_1(n)^2\right]=0.
\]
\end{lem}

\begin{proof}
We  set  for $p\geq  \ell \geq 0$, $n-p\geq k_0$ and $j\in \G_{k_0}$:
\begin{equation*}
   \label{eq:def-R1ln-L2}
R_{1,j} (\ell,n ) = \sum_{i\in j\G_{n-p-k_0}}\E\left[N_{n,i}^\ell(f_\ell) |\,
  \cf_i\right],
\end{equation*}
so that $R_1(n)=\sum_{\ell=0}^{p} \sum_{j\in \G_{k_0}} R_{1,j} (\ell,n )
$. 
We have for $i\in \G_{n-p}$:
\begin{equation}\label{eq:ENlni}
|\G_n|^{1/2} \E\left[N_{n,i}^\ell(f_\ell) |\,\cf_i\right] = \E\left[M_{i\G_{p-\ell}}(\tilde f_\ell) |X_i\right]=\E_{X_i} \left[M_{\G_{p-\ell}}(\tilde f_\ell)\right] = |\G_{p-\ell}|\, \cq^{p-\ell} \tilde f_\ell(X_i),
\end{equation}
where we used definition \reff{eq:def-Nnil} of $N_{n,i}^\ell$ for the first equality, the Markov property of $X$ for the second and  \reff{eq:Q1} for the third. Using \reff{eq:ENlni}, we get for $j\in \G_{k_0}$: 
\[
R_{1,j}(\ell,n)= |\G_n|^{-1/2} \,  |\G_{p-\ell}|\, M_{j\G_{n-p-k_0}}
(\cq^{p-\ell} \tilde f_\ell).
\]
We deduce from the Markov property of $X$ that $\E[R_{1,j}(\ell,n)^2|\, \cf_j]=2^{-n+2(p-\ell)} \,
h_{\ell, n} (X_j)$ with $h_{\ell,n}(x)=\E_x\left[M_{\G_{n-p-k_0}}
(\cq^{p-\ell} \tilde f_\ell)^2\right]$. 
We have, thanks to  (ii) from Assumption \ref{hyp:Q}, see
\reff{eq:EGn-nu0}, that: 
\[
\sum_{j\in \G_{k_0}} \E[R_{1,j}(\ell,n)^2] = 2^{-n+2(p-\ell)} \,
\E\left[M_{\G_{k_0}} (h_{\ell, n})\right]  
\leq C 2^{-n+2(p-\ell)} \,  \langle \mu,  h_{\ell,n}\rangle . 
\]
We have:
\begin{align}
\langle \mu,  h_{\ell, n} \rangle = \E_\mu
  \left[M_{\G_{n-p-k_0}}(\cq^{p-\ell} \tilde f_\ell)^2\right]  
& \leq   C  \, 2^{n-p} \sum_{k=0}^{n-p-k_0} 2^{k} \,  \normm{\cq^k
  \cq^{p-\ell} \tilde f_\ell} ^2 \nonumber 
  \\ 
&\leq  C \,2^{n-p} \alpha^{2(p-\ell)} \, c_{2}^2(\bF), \nonumber
\end{align}
where we used \reff{eq:EGn-nu0} for the first inequality (notice one can
take $k_0=0$ in this case as we consider the expectation $\E_\mu$), 
\reff{eq:L2-erg}  in the second, and $2\alpha^2<1$ in the last.
We deduce that:
\[
   \sum_{j\in \G_{k_0}} \E\left[R_{1,j}(\ell, n)^2\right]
\leq C \alpha^{2(p-\ell)} 2^{p-2\ell} c_{2}^2(\bF).
\]
We get that:
\begin{equation*}\label{eq:IR1nL2sub}
\E\left[R_1(n)^2\right]^{1/2} \leq  \sum_{\ell=0}^{p}  \left(2^{k_0} \sum_{j \in \G_{k_0}} \E\left[R_{1,j}(\ell, n)^2\right] \right)^{1/2} \leq C \,  c_{2}(\bF) \, a_{1, n},
\end{equation*}
with the sequence $(a_{1, n},n\in \N)$ defined by:
\begin{equation*}\label{eq:def-a1n}
a_{1, n}= (2\alpha^2)^{p/2} \sum_{\ell=0}^{p}
(2\alpha)^{-\ell}. 
\end{equation*}
The sequence $(a_{1, n}, n\in \N)$ does not depend on $\bF$ and
converges to 0 since $\lim_{n \rightarrow \infty} p = \infty$, $2\alpha^{2} < 1$ and
\begin{equation*}\label{eq:sum2alpha-}
\sum_{\ell=0}^{p} (2\alpha)^{-\ell}\leq  
\begin{cases} 
2\alpha/(2\alpha -1) & \text{if $2\alpha>1$},\\ 
p+1 & \text{if $2\alpha=1$},\\ 
(2\alpha)^{-p}/(1-2\alpha) & \text{if $2\alpha<1$}. 
\end{cases}
\end{equation*} 
Then use that  $\bF$ is bounded in $L^2(\mu)$ to conclude.
\end{proof}
\begin{rem}\label{rem:Nk0=Dnf}
From the proofs of Lemmas \ref{lem:cvR0-L2} and \ref{lem:cvR1-L2}, we have that $\EE[(N_{n,\emptyset}^{[k_{0}]}(\bF) - \Delta_{n}(\bF))^{2}] \leq a_{0,n} c_{2}(\bF)$, where the sequence $(a_{0,n}, n\in \NN)$ converges to $0$ as $n$ goes to infinity. 
\end{rem}

\medskip

We now study the bracket of $\Delta_n$:
\begin{equation*}\label{eq:def-Vn-subc}
V(n)= \sum_{i\in \G_{n-p_n}} \E\left[ \Delta_{n, i}(\bF)^2|\cf_i\right]. 
\end{equation*}
Using \reff{eq:def-NiF} and \reff{eq:def-DiF}, we write:
\begin{equation}
   \label{eq:def-V}
V(n) = |\G_n|^{-1} \sum_{i\in \G_{n-p_n}} \E_{X_i}\left[
    \left(\sum_{\ell=0}^{p_n} M_{\G_{p_n-\ell}}(\tilde f_\ell) \right)^2
  \right]-R_2(n)=V_1(n) +2V_2(n) - R_2( n),
\end{equation}
with:
\begin{align*}
V_1(n)
& =   |\G_n|^{-1} \sum_{i\in \G_{n-p_n}} 
\sum_{\ell=0}^{p_n} \E_{X_i}\left[
   M_{\G_{p_n-\ell}}(\tilde f_\ell) ^2  \right] ,\\
V_2(n)
& =  |\G_n|^{-1} \sum_{i\in \G_{n-p_n}} \sum_{0\leq \ell<k\leq p_n
  } \E_{X_i}\left[ 
   M_{\G_{p_n-\ell}}(\tilde f_\ell)  M_{\G_{p_n-k}}(\tilde f_k) 
  \right], \\
R_2( n)
&=\sum_{i\in \G_{n-p_n}} \E\left[ N_{n,i} (\bF)
  |X_i \right] ^2.
\end{align*}

\begin{lem}
   \label{lem:cvR2-L2}
Under the assumptions of Theorem \ref{cor:subcritical-L2}, we have the following convergence:
\[
\lim_{n\rightarrow \infty } \E\left[R_2(n)\right]=0.
\]
\end{lem}

\begin{proof}
We define the sequence $(a_{2, n}, n\in \N)$ for $n\in \N$ by: 
\begin{equation*}\label{eq:def-a2n}
a_{2, n}= 2^{-p} 
 \left(\sum_{\ell=0}^{p}  (2\alpha)^{\ell} \right)^2.
\end{equation*}
Notice that the sequence $(a_{2, n}, n\in \N)$   converges to 0 since $\lim_{n\rightarrow \infty }p=\infty $, $2\alpha^2<1$  and 
\begin{equation*}\label{eq:sum2a}
\sum_{\ell=0}^{p} (2\alpha)^{\ell}\leq  
\begin{cases}
(2\alpha)^{p+1}/(2\alpha-1) & \text{if $2\alpha>1$},\\
p+1 & \text{if $2\alpha=1$},\\
   1/(1-2\alpha ) & \text{if $2\alpha<1$}.
\end{cases}
\end{equation*}
We now compute $ \E_x\left[R_2(n)\right]$. 
\begin{align}
\nonumber
  \E_x\left[R_2(n)\right]
 &= |\G_n|^{-1} \, \sum_{i\in \G_{n-p}} \E_x\left[
\E_x\left[\sum_{\ell=0}^{p} M_{i\G_{p-\ell}}(\tilde f_\ell)
  |X_i\right]^2
\right]\\
\nonumber
 &= |\G_n|^{-1} \, \sum_{i\in \G_{n-p}} \E_x\left[
\left(\sum_{\ell=0}^{p} \E_{X_i} \left[M_{\G_{p-\ell}}(\tilde f_\ell)
 \right]\right)^2
\right].
\\ \nonumber
&= |\G_n|^{-1} \,  |\G_{n-p}|\, 
\cq^{n-p} \left(\Big(\sum_{\ell=0}^{p}  |\G_{p-\ell}|\, \cq^{p-\ell} \tilde
  f_\ell\Big)^2\right)(x)
\end{align}
where we used the definition of $N_{n,i}(\bF)$ for the first equality, the Markov property of $X$ for the second,  \reff{eq:Q1} for the third. From the latter equality, we have using (ii) from Assumption \ref{hyp:Q}:
\begin{align}
\nonumber
\E\left[R_2(n)\right] 
&= |\G_n|^{-1} \,  |\G_{n-p}|\, \langle \nu,\cq^{n-p}
  \left(\Big(\sum_{\ell=0}^{p}  |\G_{p-\ell}|\, \cq^{p-\ell} \tilde
  f_\ell\Big)^2\right) \rangle \\ \nonumber
&\leq C 2^{-p} \left(\sum_{\ell=0}^{p}  |\G_{p-\ell}|\, \norm{
  \cq^{p-\ell} \tilde f_\ell}_{L^2(\mu)}\right)^2. 
\end{align}
We deduce that:
\begin{equation*}\label{eq:majoR2-subL2-vit}
\E\left[R_2(n)\right] \leq  C\, c_2^2(\bF) \, a_{2, n},
\end{equation*}
Then use that  $\bF$ is bounded in $L^2(\mu)$ to conclude.
\end{proof}
\begin{rem}\label{rem:Vn=V1n+V2n}
In particular, we have obtained from the previous proof that $\EE[|V(n) - V_{1}(n) - V_{2}(n)|] \leq C c_{2}^{2}(\bF) a_{2,n}$, with the sequence $(a_{2,n}, n \in \NN)$ going to 0 as $n$ goes to infinity.
\end{rem}

\begin{lem}\label{lem:cvV2-L2}
Under the assumptions of Theorem \ref{cor:subcritical-L2},  we have that in probability $ \lim_{n\rightarrow \infty } V_2(n) =\ssub_2(\bF)$ with $\ssub_2(\bF)$ finite and defined in \reff{eq:S2}. 
\end{lem}

\begin{proof}
Using \reff{eq:Q2-bis}, we get:
\begin{equation}\label{eq:decom-V2}
V_2(n)= V_5(n)+ V_6(n),
\end{equation}
with
\begin{align*}
V_5(n) &=  |\G_n|^{-1} \sum_{i\in \G_{n-p}} \sum_{0\leq \ell<k\leq  p } 2^{p-\ell} \cq^{p-k} \left( \tilde f_k \cq^{k-\ell} \tilde f_\ell\right)(X_i),\\
V_6(n) &=     |\G_n|^{-1} \sum_{i\in \G_{n-p}} \sum_{0\leq \ell<k<  p } \sum_{r=0}^{p-k-1}  2^{p-\ell+r} \, \cq^{p-1-(r+k)}\left(\cp\left(\cq^r \tilde f_k \sot \cq  ^{k-\ell+r} \tilde f_\ell \right)\right)(X_i).
\end{align*}
We consider the term $V_6(n)$. We have:
\begin{equation*}
V_6(n)=|\G_{n-p}|^{-1} M_{\G_{n-p}} (H_{6,n}),
\end{equation*}
with:
\begin{equation}\label{eq:H6n-crit-det}
H_{6,n}=\sum_{\substack{0\leq \ell< k \\ r\geq 0}}
h_{k,\ell,r}^{(n)}\,  \ind_{\{r+k<  p\}} 
\quad\text{and} \quad 
h_{k, \ell,r}^{(n)} =  2^{r-\ell} \, \cq^{p-1-(r+k)}\left(\cp\left(\cq^r \tilde f_k \sot \cq  ^{k-\ell+r} \tilde f_\ell \right)\right).
\end{equation}
Define $H_6(\bF)=\sum_{0\leq \ell< k; r\geq 0} h_{k,\ell,r}$ with $h_{k,\ell,r}=2^{r-\ell} \langle \mu, \cp\left(\cq^r \tilde f_k \sot \cq^{k-\ell+r} \tilde f_\ell\right)\rangle=\langle \mu, h_{k, \ell,r}^{(n)}\rangle$. 
Thanks to \reff{eq:mp}  and  \reff{eq:L2-erg}, we get that:
\begin{equation}
   \label{eq:majo-H6h}
|h_{k,\ell,r}|
\leq C \, 2^{r-\ell} \normm{\cq^r \tilde f_k}\normm{\cq^{k-\ell +r}
  \tilde f_\ell} 
\leq C \, 2^{r-\ell} \alpha^{k-\ell
  +2r}\normm{f_\ell} \normm{f_k}.
\end{equation}
We deduce that $|h_{k,\ell,r}|\leq C \, 2^{r-\ell} \alpha^{k-\ell
  +2r} c_2^2(\bF)$ and, as   the sum  
$\sum_{0\leq \ell<k, \, r\geq 0}2^{r-\ell} \alpha^{k-\ell +2r}
$ is finite: 
\begin{equation}
   \label{eq:majoH6f-L2}
|H_6(\bF)|\leq C \, c_2^2(\bF).
\end{equation}
We write $H_6(\bF)=H_6^{[n]}(\bF)+ B_{6,n}(\bF)$, with
\begin{equation*}\label{eq:def-B6nsubL2}
H_6^{[n]}(\bF)=\sum_{\substack{0\leq \ell< k \\ r\geq 0}} h_{k,\ell,r} \,  \ind_{\{r+k<  p\}} \quad\text{and}\quad B_{6,n}(\bF)=\sum_{\substack{0\leq \ell< k \\ r\geq 0}} h_{k,\ell,r} \,  \ind_{\{r+k\geq   p\}} .
\end{equation*}
As  $\lim_{n\rightarrow \infty } \ind_{\{r+k\geq p\}}=0$, we get from
\reff{eq:majo-H6h}, 
\reff{eq:majoH6f-L2} 
and dominated convergence  that
$\lim_{n\rightarrow \infty } B_{6,n}(\bF)=0$ and thus:
\begin{equation}
   \label{eq:limH6n}
\lim_{n\rightarrow \infty } H_6^{[n]}(\bF)=H_6(\bF).
\end{equation}
We set $A_{6,n}(\bF)=H_{6,n}-H_6^{[n]}(\bF)=  \sum_{\substack{0\leq \ell< k \\ r\geq 0}}
(h_{k,\ell,r}^{(n)}-h_{k,\ell,r}) \,  \ind_{\{r+k<  p\}}$, so that from
the    definition     of     $V_{6}(n)$,     we    get     that:
\begin{equation*}
   \label{eq:V6-H6f}
  V_6(n)-  H_6^{[n]}(\bF)=  |\G_{n-p}|^{-1}\, M_{\G_{n-p}}  (A_{6,n}(\bF)).
\end{equation*}
We now study the second moment of $ |\G_{n-p}|^{-1}\, M_{\G_{n-p}}
(A_{6,n}(\bF))$.  Using \reff{eq:EGn-nu0}, we get for $n-p\geq k_0$:
\[
|\G_{n-p}|^{-2}\, \E\left[M_{\G_{n-p}} (A_{6,n}(\bF))^2 
\right]
\leq  C\, |\G_{n-p}|^{-1}\, \sum_{j=0}^{n-p} 2^j \normm{\cq^j (A_{6,
    n}(\bF))}^2. 
\]
Recall $c_k(\bF)$ and $q_k(\bF)$ from \reff{eq:def-ck}. We deduce that
\begin{align}
\nonumber
   \normm {\cq^j (A_{6,n}(\bF))}
&\leq \sum_{\substack{0\leq \ell< k\\ r\geq 0}}  \normm {\cq^j
  h_{k,\ell,r}^{(n)}-h_{k,\ell,r}} 
  \ind_{\{r+k<  p\}} \\ 
\nonumber
&\leq C \sum_{\substack{0\leq \ell< k\\ r\geq 0}} 
  2^{r-\ell } \, \alpha^{p-1-(r+k)+j} 
\normm{\cp\left(\cq^r \tilde f_k \sot \cq  ^{k-\ell+r} \tilde f_\ell
  \right)} \ind_{\{r+k<  p\}}\\ \nonumber
&\leq C c_2^2(\bF)\, \alpha^j \, \sum_{\substack{0\leq \ell< k\\
  r\geq 1}}  
  2^{r-\ell } \,\alpha^{p-(r+k)} \alpha^{k-\ell+2r} 
\, \ind_{\{r+k<  p\}}\\ \nonumber
&\hspace{2cm} + C  
\alpha^j \, \sum_{0\leq \ell< k}  
  2^{-\ell } \,\alpha^{p-k} \normm{\cp\left(\tilde f_k \sot \cq  ^{k-\ell} \tilde f_\ell
  \right)} 
\, \ind_{\{k<  p\}}\\
\nonumber
&\leq C c_2(\bF)c_4(\bF) \, \alpha^j \, \sum_{\substack{0\leq \ell< k\\
  r\geq 0}}  
  2^{r-\ell } \,\alpha^{p-(r+k)} \alpha^{k-\ell+2r} 
\, \ind_{\{r+k<  p\}}\\ \nonumber
&\leq C c_2(\bF)c_4(\bF) \, \alpha^j, 
\end{align}
where  we  used the  triangular  inequality  for the  first  inequality;
\reff{eq:L2-erg}   for   the   second;   \reff{eq:P-L2-1}   for   $r\geq
1$  and \reff{eq:L2-erg} again for the third;
\reff{eq:P-L2-2} for $r=0$ to get the $c_4(\bF)$ term and $c_2(\bF)\leq
c_4(\bF)$ for the fourth;
and that   $  \sum_{{0\leq \ell< k,\,  r\geq 0}} 
  2^{r-\ell} \alpha^{k
    -\ell+2r} $ is finite for the
last. As $\sum_{j=0}^{\infty } (2\alpha^2)^j$ is finite, we deduce that:
\begin{equation}
   \label{eq:majo-L2A6}
\E\left[\left( V_6(n)-  H_6^{[n]}(\bF)\right)^2\right]=
|\G_{n-p}|^{-2}\, \E\left[M_{\G_{n-p}} (A_{6,n}(\bF))^2 
\right]\leq C c_2^2(\bF)c_4^2(\bF) \, 2^{-(n-p)}.
\end{equation}
\medskip

We now consider the term $V_5(n)$ defined just after \reff{eq:decom-V2}:
\[
V_5(n)=|\G_{n-p}|^{-1} M_{\G_{n-p}} (H_{5,n}),
\]
with
\[
H_{5,n}=\sum_{0\leq \ell< k }
h_{k,\ell}^{(n)}\,  \ind_{\{k\leq   p\}} 
\quad\text{and} \quad 
h_{k, \ell}^{(n)} =  2^{-\ell} \, \cq^{p-k}\left(\tilde f_k
  \cq  ^{k-\ell} \tilde f_\ell \right). 
\]
Define $H_5(\bF)=\sum_{0\leq \ell< k} h_{k,\ell}$ with $h_{k,\ell}= 2^{-\ell} \langle \mu, \tilde f_k \cq^{k-\ell} \tilde f_\ell\rangle$.
We have using Cauchy-Schwartz inequality and  \reff{eq:L2-erg} that:
\begin{equation}
   \label{eq:majo-H5h}
|h_{k,\ell}|\leq C 
 \, 2^{-\ell} \alpha^{k-\ell}\normm{f_\ell} \normm{f_k}\leq C 
 \, 2^{-\ell} \alpha^{k-\ell} \,c_2^2(\bF).
\end{equation}
As   the sum  
$\sum_{0\leq \ell<k}2^{-\ell} \alpha^{k-\ell}
$ is finite, we deduce that: 
\begin{equation}
   \label{eq:majoH5f-L2}
|H_5(\bF)|\leq C \, c_2^2(\bF).
\end{equation}
We write $H_5(\bF)=H_5^{[n]}(\bF)+ B_{5,n}(\bF)$, with
\begin{equation}
   \label{eq:def-B5nsubL2}
H_5^{[n]}(\bF)=\!\!\sum_{0\leq \ell< k} \!\! h_{k,\ell}\ind_{\{k\leq p\}}
=\!\!\sum_{0\leq \ell< k}\!\!
2^{-\ell} \langle \mu,
\tilde
    f_k \cq^{k-\ell} \tilde f_\ell\rangle\ind_{\{k\leq p\}}
\quad\text{and}\quad
B_{5,n}(\bF)=\!\!\sum_{0\leq \ell< k }\!\!
h_{k,\ell} \,  \ind_{\{k>  p\}} .
\end{equation}
As  $\lim_{n\rightarrow  \infty  }   \ind_{\{k>  p\}}=0$,  we  deduce  from
\reff{eq:majo-H5h}  and \reff{eq:majoH5f-L2} that   $\lim_{n\rightarrow
  \infty  } B_{5,n}(\bF)=0$ by  dominated  convergence and thus:
\begin{equation}
   \label{eq:limH5n}
\lim_{n\rightarrow \infty } H_5^{[n]}(\bF)=H_5(\bF).
\end{equation}
We set $A_{5,n}(\bF)= H_{5,n} -H_5^{[n]}(\bF)= \sum_{0\leq \ell< k }
(h_{k,\ell}^{(n)}-h_{k,\ell}) \,  \ind_{\{k\leq   p\}}$, so that from     the    definition     of     $V_{5}(n)$,     we    get     that:
\begin{equation}
   \label{eq:V5-H5f}
  V_5(n)-  H_5^{[n]}(\bF)=  |\G_{n-p}|^{-1}\, M_{\G_{n-p}}  (A_{5,n}(\bF)).
\end{equation}
We now  study the
second  moment of  $  |\G_{n-p}|^{-1}\,  M_{\G_{n-p}} (A_{5,n}(\bF))  $.
Using \reff{eq:EGn-nu0}, we get for $n-p\geq k_0$:
\[
|\G_{n-p}|^{-2}\, \E\left[M_{\G_{n-p}} (A_{5,n}(\bF))^2 
\right]
\leq  C\, |\G_{n-p}|^{-1}\, \sum_{j=0}^{n-p} 2^j \normm{\cq^j (A_{5,
    n}(\bF))}^2. 
\]
We also have  that:
\begin{align}
\nonumber
   \normm {\cq^j(A_{5,n}(\bF))}
&\leq \sum_{0\leq \ell< k}  \normm {\cq^j h_{k,\ell}^{(n)}-h_{k,\ell}}
  \ind_{\{k\leq   p\}} \\ \nonumber
&\leq C \sum_{0\leq \ell< k}
  2^{-\ell } \, \alpha^{p-k+j} 
\normm{ \tilde f_k \cq  ^{k-\ell} \tilde f_\ell
} \ind_{\{k\leq   p\}}\\
\label{eq:L2A5Hf}
&\leq C c_4^2(\bF) \, \alpha^j , 
\end{align}
where we  used the  triangular inequality for  the first  inequality, 
\reff{eq:L2-erg} for the second, and Cauchy-Schwartz inequality for the
last.  As $\sum_{j=0}^{\infty }
(2\alpha^2)^j$ is finite, we deduce that: 
\begin{equation}\label{eq:majo-L2A5}
\E\left[\left(V_5(n)-  H_5^{[n]}(\bF)\right)^2\right] = |\G_{n-p}|^{-2}\, \E\left[M_{\G_{n-p}} (A_{5,n}(\bF))^2 \right] \leq  C\,c_4^4(\bF)\, 2^{-(n-p)} . 
\end{equation}
\medskip

Since $c_2(\bF) \leq  c_4(\bF)$, we deduce from \reff{eq:majo-L2A6}  and
\reff{eq:majo-L2A5}, as $V_2(n)=V_5(n)+V_6(n)$ (see \reff{eq:decom-V2}),
that: 
\begin{equation*}
   \label{eq:majo-L2A2}
\E\left[\left( V_2(n)-  H_2^{[n]}(\bF)\right)^2\right]
\leq C \, c_4^4(\bF) \, 2^{-(n-p)}
\quad\text{with}\quad
 H_2^{[n]}(\bF)= H_6^{[n]}(\bF)+ 
H_5^{[n]}(\bF).
\end{equation*}
 Since, according to \reff{eq:limH6n}
and \reff{eq:limH5n} and 
$\ssub_2(\bF)=H_6(\bF)+H_5(\bF)$ (see \reff{eq:S2}), we get 
$\lim_{n\rightarrow \infty } H_2^{[n]}(\bF)=\ssub_2(\bF)$. This implies
that   $\lim_{n\rightarrow   \infty  }
V_2(n)=\ssub_2(\bF)$  in 
  probability. 

\end{proof}

We now study the limit of $V_1(n)$. 

\begin{lem}\label{lem:cvV1-L2}
Under the assumptions of Theorem \ref{cor:subcritical-L2},  we have that in probability $\lim_{n\rightarrow \infty } V_1(n) =\ssub_1(\bF)<+\infty $ with $\ssub_1(\bF)$ finite and defined in \reff{eq:S1}. 
\end{lem}
\begin{proof}
Using \reff{eq:Q2}, we get:
\begin{equation}\label{eq:DV4nsub}
V_1(n)= V_3(n)+ V_4(n),
\end{equation}
with
\begin{align*}
V_3(n) &=  |\G_n|^{-1} \sum_{i\in \G_{n-p}} \sum_{\ell=0}^p 2^{p-\ell}\, \cq^{p-\ell} (\tilde f_\ell^2)(X_i),\\ 
V_4(n) &=     |\G_n|^{-1} \sum_{i\in \G_{n-p}} \sum_{\ell=0}^{p-1}\,\sum_{k=0}^{p-\ell -1} 2^{p-\ell+k} \,\cq^{p-1-(\ell+k)}\left(\cp\left(\cq^k \tilde f_\ell \otimes^2\right)\right)(X_i).  
\end{align*}
We first consider the term $V_4(n)$. We have:
\[
V_4(n)=|\G_{n-p}|^{-1} M_{\G_{n-p}} (H_{4,n}),
\]
with:
\[
H_{4,n} = \sum_{\ell\geq 0, \, k\geq 0} h_{\ell,
  k}^{(n)}\,  \ind_{\{\ell+k<  p\}} \quad \text{and} \quad h_{\ell,
  k}^{(n)} = 2^{k-\ell}\,  \cq^{p-1-(\ell+k)} \left(\cp\left(\cq^k
    \tilde f_\ell \otimes^2 \right)\right). 
\]
Define the constant $H_4(\bF)= \sum_{\ell\geq 0, \, k\geq 0} h_{\ell, k} $ with $h_{\ell, k} = 2^{k-\ell}\,  \langle \mu,\cp\left(\cq^k \tilde f_\ell \otimes^2 \right) \rangle$. 
Thanks to \reff{eq:majo-pfxf} and \reff{eq:L2-erg}, we have:
\begin{equation}
   \label{eq:majo-H4h}
|h_{\ell,k}|\leq  2^{k-\ell} \normm{\cq^k \tilde f_\ell}^2 
\leq  C\, 2^{k-\ell} \alpha^{2k}\normm{f_\ell}^2
\leq  C\, 2^{k-\ell} \alpha^{2k} \, c_2^2(\bF),
\end{equation}
and thus, as the sum $\sum_{\ell\geq 0, \, k\geq 0} 2^{k-\ell}
\alpha^{2k} $ is finite:
\begin{equation}
   \label{eq:majoH4f-L2}
|H_4(\bF)|\leq C \, c_2^2(\bF).
\end{equation}
We write $H_4(\bF)=H_4^{[n]}(\bF)+ B_{4,n}(\bF)$, with
\begin{equation*}\label{eq:def-B4nsubL2}
H_4^{[n]}(\bF)= \sum_{\ell\geq 0, \, k\geq 0} h_{\ell,k} \,  \ind_{\{\ell+k<  p\}} \quad\text{and}\quad B_{4,n}(\bF)=\sum_{\ell\geq 0, \, k\geq 0} h_{\ell,k} \,  \ind_{\{\ell+k\geq   p\}}  .
\end{equation*}
Using that  $\lim_{n\rightarrow \infty } \ind_{\{\ell+k\geq p\}}=0$, we deduce  from
\reff{eq:majo-H4h}, 
\reff{eq:majoH4f-L2} 
and dominated convergence  that
$\lim_{n\rightarrow \infty } B_{4,n}(\bF)=0$,
and thus:
 \begin{equation}
   \label{eq:limH4n}
\lim_{n\rightarrow \infty } H_4^{[n]}(\bF)=H_4(\bF).
\end{equation}
We set $
A_{4,n}(\bF)= H_{4,n} -H_4^{[n]}(\bF)= \sum_{\ell\geq 0, \, k\geq 0}
(h_{\ell,k}^{(n)}-h_{\ell,k}) \,  \ind_{\{\ell+k<  p\}} $, so that from
the    definition     of     $V_{4}(n)$,     we    get     that: 
\begin{equation*}
   \label{eq:V6-H4f}
  V_4(n)-  H_4^{[n]}(\bF)=  |\G_{n-p}|^{-1}\, M_{\G_{n-p}}  (A_{4,n}(\bF)).
\end{equation*}
We now study the second moment of $ |\G_{n-p}|^{-1}\, M_{\G_{n-p}}
(A_{4,n}(\bF))$.  Using \reff{eq:EGn-nu0}, we get for $n-p\geq k_0$:
\[
|\G_{n-p}|^{-2}\, \E\left[M_{\G_{n-p}} (A_{4,n}(\bF))^2 
\right]
\leq  C\, |\G_{n-p}|^{-1}\, \sum_{j=0}^{n-p} 2^j \normm{\cq^j (A_{4,
    n}(\bF))}^2. 
\]
Using \reff{eq:majo-pfxf}, we obtain  that $\normm{\cp(\tilde f_\ell
  \otimes \tilde f_\ell)} \leq  c_4^2(\bF)$. 
We deduce that:
\begin{align}
\nonumber
   \normm {\cq^j (A_{4,n}(\bF))}
&\leq \sum_{\ell\geq 0,\, k\geq 0}  \normm {\cq^j
  h_{\ell,k}^{(n)}-h_{\ell,k}} 
  \ind_{\{\ell+k<  p\}} \\ 
\nonumber
&\leq C \sum_{\ell\geq 0,\, k\geq 0} 
  2^{k-\ell} \, \alpha^{p-1-(\ell+k)+j} 
\normm{\cp\left(\cq^k \tilde f_\ell \otimes^2
  \right)} \ind_{\{\ell+k<  p\}}\\
\nonumber
&\leq C \, c_2^2(\bF) \, \alpha^j\!\! \!  \sum_{\ell\geq 0,\, k> 0} \! \!\! 
  2^{k-\ell } \,\alpha^{p-(\ell+k)} \alpha^{2k} 
\, \ind_{\{\ell+k<  p\}}\\ \nonumber
&\hspace{2cm}
 +  C \, \alpha^j \sum_{\ell\geq 0} 
  2^{-\ell} \, \alpha^{p-\ell} 
\normm{\cp\left(\tilde f_\ell \otimes ^2
  \right)} \ind_{\{\ell<  p\}}\\ \nonumber
&\leq C \, c_4^2(\bF) \, \alpha^j,
\end{align}
where  we  used the  triangular  inequality  for the  first  inequality;
\reff{eq:L2-erg}   for   the   second;   \reff{eq:P-L2-1}   for   $k\geq
1$ and \reff{eq:L2-erg} again for the third; and \reff{eq:majo-pfxf} as
well as $c_2(\bF)\leq c_4(\bF)$ for the last. 
As $\sum_{j=0}^{\infty } (2\alpha^2)^j$ is finite, we deduce that:
\begin{equation}
   \label{eq:majo-L2A4}
\E\left[\left(V_4(n)-  H_4^{[n]}(\bF)\right)^2\right]
=|\G_{n-p}|^{-2}\, \E\left[M_{\G_{n-p}} (A_{4,n}(\bF))^2 
\right]\leq C  \, c_4^4(\bF)\, 2^{-(n-p)}.
\end{equation}
\medskip

We now consider the term $V_3(n)$ defined just after \reff{eq:DV4nsub}:
\[
V_3(n)=|\G_{n-p}|^{-1} M_{\G_{n-p}} (H_{3,n}),
\]
with 
\[
H_{3,n} = \sum_{\ell\geq 0} h_{\ell}^{(n)}\,  \ind_{\{\ell\leq   p\}}
\quad \text{and} \quad h_{\ell}^{(n)} = 2^{-\ell}\,  \cq^{p-\ell}
\left(\tilde f_\ell^2\right). 
\]
Define the constant  $H_3(\bF)=
\sum_{\ell\geq 0} h_\ell$
 with
$h_\ell= 2^{-\ell} \langle
\mu,  \tilde f_\ell^2 \rangle=\langle \mu, h^{(n)}_\ell \rangle$.
As   $h_\ell  \leq   \normm{  f_\ell}^2\leq   c_2^2(\bF)$,  we   get  that
$H_3(\bF)\leq          2         c_2^2(\bF)$.      
We write $H_3(\bF)=H_3^{[n]}(\bF)+ B_{3,n}(\bF)$, with
\begin{equation*}\label{eq:def-B3nsubL2}
H_3^{[n]}(\bF)= \sum_{\ell\geq 0} h_{\ell} \,  \ind_{\{\ell\leq p\}}\quad\text{and}\quad B_{3,n}(\bF)=\sum_{\ell\geq 0} h_{\ell} \,  \ind_{\{\ell>  p\}}  .
\end{equation*}
As  $\lim_{n\rightarrow \infty } \ind_{\{\ell> p\}}=0$, we get from
 dominated convergence  that
$\lim_{n\rightarrow \infty } B_{3,n}(\bF)=0$ and thus:
\begin{equation}
   \label{eq:limH3n}
\lim_{n\rightarrow \infty } H_3^{[n]}(\bF)=H_3(\bF).
\end{equation}
We set $
A_{3,n}(\bF)= H_{3,n} -H_3^{[n]}(\bF)= \sum_{\ell\geq 0}
(h_{\ell}^{(n)}-h_{\ell}) \,  \ind_{\{\ell\leq   p\}} $, 
so that from     the    definition     of     $V_{3}(n)$,     we    get
that:
\begin{equation}\label{eq:V3-H3f}
V_3(n) -  H_3^{[n]}(\bF) =  |\G_{n-p}|^{-1}\, M_{\G_{n-p}}  (A_{3,n}(\bF)).
\end{equation}
We now study the second moment of $ |\G_{n-p}|^{-1}\,
M_{\G_{n-p}}(A_{3,n}(\bF))$.  Using \reff{eq:EGn-nu0}, we get for
$n-p\geq k_0$: 
\begin{equation*}\label{eq:majoA3}
|\G_{n-p}|^{-2}\, \E\left[M_{\G_{n-p}} (A_{3,n}(\bF))^2 \right] \leq  C\, |\G_{n-p}|^{-1}\, \sum_{j=0}^{n-p} 2^j \normm{\cq^j (A_{3,n}(\bF))}^2. 
\end{equation*}
We have  that
\begin{align}
\nonumber
   \normm {\cq^j (A_{3,n}(\bF))}
&\leq \sum_{\ell\geq 0}  \normm {\cq^j
  h_{\ell}^{(n)}-h_{\ell}} 
  \ind_{\{\ell\leq   p\}} \\ \nonumber 
&\leq C \sum_{\ell\geq 0}  2^{-\ell} \,
\normm{ \cq^{j+p-\ell} \tilde g}
 \ind_{\{\ell\leq  p\}}\quad\text{with}\quad g=\tilde f^2_\ell\\ \nonumber
&\leq C \sum_{\ell\geq 0}
  2^{-\ell} \, \alpha^{j+p-\ell} 
\normm{\tilde f_\ell^2 } \ind_{\{\ell\leq  p\}}\\
\nonumber
&\leq C \, c_4^2(\bF) \, \alpha^j,
\end{align}
where  we  used the  triangular  inequality  for the  first  inequality;
and \reff{eq:L2-erg}    for   the  third. As $\sum_{j=0}^{\infty } (2\alpha^2)^j$ is finite, we deduce that:
\begin{equation}\label{eq:majo-L2A3}
\E\left[\left( V_3(n)-  H_3^{[n]}(\bF)\right)^2\right] = |\G_{n-p}|^{-2}\, \E\left[M_{\G_{n-p}} (A_{3,n}(\bF))^2 \right]\leq C  \, c_4^4(\bF)\, 2^{-(n-p)}.
\end{equation}
\medskip

Since $c_2(\bF) \leq  c_4(\bF)$, we deduce from \reff{eq:majo-L2A4} and \reff{eq:majo-L2A3} that:
\begin{equation*}
   \label{eq:majo-L2A1}
\E\left[\left( V_1(n)-  H_1^{[n]}(\bF)\right)^2\right]
\leq C \,c_4^4(\bF) \, 2^{-(n-p)}
\quad\text{with}\quad
H_1^{[n]}(\bF)= H_4^{[n]}(\bF)+ H_3^{[n]}(\bF). 
\end{equation*}
 Since,
according to \reff{eq:limH4n} and \reff{eq:limH3n} 
$\ssub_1(\bF)=H_4(\bF)+H_3(\bF)$ (see \reff{eq:S1}), we get 
$\lim_{n\rightarrow \infty } H_1^{[n]}(\bF)=\ssub_1(\bF)$. This implies
that   $\lim_{n\rightarrow   \infty  }
V_1(n)=\ssub_1(\bF)$  in 
  probability. 

\end{proof}

 The  next Lemma is a direct consequence of \reff{eq:def-V} and
 Lemmas \ref{lem:cvR2-L2}, \ref{lem:cvV2-L2} and \ref{lem:cvV1-L2}.
\begin{lem}
   \label{lem:cvV-L2}
   Under the assumptions of  Theorem \ref{cor:subcritical-L2}, we have
   $\lim_{n\rightarrow\infty }  V(n)=\ssub(\bF)$ in  probability, where,
   with $\ssub_1(\bF)$  and $\ssub_2(\bF)$  defined by  \reff{eq:S1} and
   \reff{eq:S2}, we have:
\[
\ssub(\bF) = \ssub_1(\bF)+ 2\ssub_2(\bF). 
\]
\end{lem}
\medskip

We now check the Lindeberg condition using a fourth moment
condition. We set
\begin{equation}\label{eq:def-R3}
R_3(n)=\sum_{i\in \G_{n-p_n}} \E\left[\Delta_{n,i}(\bF)^4\right].
\end{equation} 

\begin{lem}\label{lem:cvG-L2}
Under the assumptions of Theorem \ref{cor:subcritical-L2},  we have that $\lim_{n\rightarrow\infty } R_3(n)=0$.
\end{lem}
\begin{proof}
We have:
\begin{align}
\nonumber R_3(n) &\leq   16 \sum_{i\in \G_{n-p}}
 \E\left[N_{n,i}(\bF)^4\right]\\ \nonumber
&\leq   16 (p+1)^3 \sum_{\ell=0}^p \sum_{i\in \G_{n-p}} \E\left[N_{n,i}^\ell(\tilde f_\ell)^4\right],
\end{align}
where we used that $(\sum_{k=0}^r a_k)^4 \leq  (r+1)^3 \sum_{k=0}^r a_k^4$ for the two inequalities (resp. with $r=1$ and $r=p$) and also Jensen inequality and \reff{eq:def-DiF} for the first and \reff{eq:def-NiF} for the last.  Using \reff{eq:def-Nnil}, we get:
\begin{equation*}\label{eq:def-R3hnl}
\E\left[N_{n,i}^\ell(\tilde f_\ell)^4\right]= |\G_n|^{-2} \E\left[h_{n,\ell}(X_i)\right], \quad\text{with}\quad h_{n,\ell}(x)=\E_x\left[M_{\G_{p-\ell}}(\tilde f_\ell) ^4\right], 
\end{equation*}
so  that:
\[
R_3(n) \leq   C n^3 \sum_{\ell=0}^p \sum_{i\in \G_{n-p}} |\G_n|^{-2} \E\left[h_{n,\ell}(X_i)\right].
\]
Using \reff{eq:EGn-nu0} (with $f$ and $n$ replaced by $h_{n, \ell}$ and $n-p$), we get that:
\begin{equation}
   \label{eq:R3-L2-9}
R_3(n) \leq   C\,  n^3\, 2^{-n-p}\,   \sum_{\ell=0}^p \E_\mu\left[M_{\G_{p-\ell}}(\tilde f_\ell) ^4\right].
\end{equation}
Now we give the main steps to get  an upper bound of
$\E_\mu\left[M_{\G_{p-\ell}}(\tilde f_\ell) ^4\right]$. Recall that:
\[
\norm{\tilde f_\ell}_{L^4(\mu)}\leq  C\, c_4(\bF).
\]
We have:
\begin{equation}
   \label{eq:L2R3-l=p}
\E_\mu\left[M_{\G_{p-\ell}}(\tilde f_\ell) ^4\right]
\leq C\,   c_4^4(\bF) \quad \text{for $\ell\in \{p-2, p-1, p\}$.}
\end{equation}
Now  we consider  the case  $0\leq \ell  \leq p-3$.   Let the  functions
$\psi_{j,p-\ell}$, with $1\leq  j\leq 9$, from Lemma  \ref{lem:M4}, with 
 $f$  replaced by  $\tilde  f_\ell$ so that for $\ell\in \{0, \ldots,
 p-3\} $
\begin{equation}
   \label{eq:R3L2-Em4}
\E_\mu\left[M_{\G_{p-\ell}}(\tilde f_\ell) ^4\right]
=\sum_{j=1}^9 \langle \mu, \psi_{j,p-\ell}\rangle.
\end{equation}
We now assume that $p-\ell-1\geq 2$.  We shall  give 
bounds on  $\langle \mu, \psi_{j,p-\ell}\rangle$
based on computations  similar to those in the second  step in the proof
of Theorem 2.1 in \cite{BDG14}. We set $h_k=\cq^{k-1} \tilde f_\ell$ so
that for $k\in \N^{*}$:
\begin{equation}
   \label{eq:hk}
\normm{h_k}\leq  C \, \alpha^k c_2(\bF)
\quad\text{and} \quad
\norm{h_k}_{L^4(\mu)} \leq  C \,  c_4(\bF).
\end{equation}
We recall the notation $f\otimes f=f\otimes^2$. We deduce for $k\geq 2$ from
\reff{eq:P-L2-1} applied with $h_k=\cq h_{k-1}$ and for $k=1$ from \reff{eq:mp2}
and \reff{eq:hk} that: 
\begin{equation}
   \label{eq:phk}
\normm{\cp(h_k\otimes ^2)} \leq 
\begin{cases}
 C\, \alpha^{2k} c_2^2(\bF)&
\text{for $k\geq 2$,}   \\
 C\,  c_4^2(\bF)&
\text{for $k=1$.}   
\end{cases}
\end{equation}

\medskip
\textbf{Upper bound of $\langle \mu, |\psi_{1,p-\ell}| \rangle$}. We have:
\begin{equation*}
   \label{eq:L2R3p1}
 \langle \mu, |\psi_{1,p-\ell}| \rangle\leq C\,  2^{p-\ell}\, \langle \mu,
 \cq^{p-\ell} (\tilde f_\ell^4) \rangle
\leq C\,  2^{p-\ell}\, c_4^4(\bF).
\end{equation*}

\medskip
\textbf{Upper bound of $|\langle \mu, \psi_{2,p-\ell} \rangle|$}. Using Lemma
\ref{lem:majo-f3} for the second inequality and \reff{eq:hk} for the third, 
we get:  
\begin{align}
 \nonumber
 |\langle \mu, \psi_{2, p-\ell} \rangle|
&\leq  C    2^{2(p-\ell)}\,
 \sum_{k=0}^{p-\ell-1} 2^{-k} \langle \mu, \cq^k
 \cp  \left( \cq^{p-\ell-k - 1}( |\tilde f_\ell|^3) \sot |h_{p-\ell- k}|
  \right) \rangle \\
\nonumber
&\leq  C    2^{2(p-\ell)}\,
 \sum_{k=0}^{p-\ell-1} 2^{-k} \, c_4^3(\bF)\, \norm{ h_{p-\ell
  -k}}_{L^4(\mu)} \\ 
  \nonumber
&\leq   C\, 2^{2(p-\ell)}\, c_4^4(\bF).
\end{align}

\medskip
\textbf{Upper bound of $\langle \mu, |\psi_{3,p-\ell} \rangle|$}. Using
\reff{eq:mp}, we
easily get:
\begin{equation*}
\label{eq:L2R3p3}
\langle \mu,| \psi_{3, p-\ell}| \rangle
\leq  C \, 2^{2(p-\ell)} \sum_{k=0}^{p-\ell-1} 2^{-k}\, \langle \mu, 
  \cq^k\cp \left( \cq^{p-\ell-k - 1} (\tilde f_\ell^2) \otimes^2
  \right) \rangle
\leq   C\, 2^{2(p-\ell)}\, c_4^4(\bF).
\end{equation*}

\medskip
\textbf{Upper bound of $\langle \mu, |\psi_{4,p-\ell} \rangle|$}. Using \reff{eq:mp} and then \reff{eq:phk} with  $p-\ell -1\geq 2$, we get:
\begin{align}
\nonumber
\langle \mu,| \psi_{4, p-\ell}| \rangle
& \leq C \, 2^{4(p-\ell)} \, 
\langle \mu, \cp \left( |\cp(h_{p-\ell-1}\otimes^2)\otimes^2|\right) \rangle
  \\ \nonumber
& \leq C \, 2^{4(p-\ell)} \, \normm{\cp(h_{p-\ell-1}\otimes^2)}^2\\ 
\nonumber
& \leq C \, 2^{4(p-\ell)} \, \alpha^{4(p-\ell)} \, c_2^4(\bF)\\ \nonumber
& \leq C \, 2^{2(p-\ell)} \, c_2^4(\bF).
\end{align}

\medskip
\textbf{Upper bound of $\langle \mu, |\psi_{5,p-\ell} \rangle|$}. 
We have:
\[
\langle \mu,| \psi_{5, p-\ell}| \rangle
 \leq C \, 2^{4(p-\ell)} \, 
 \sum_{k=2}^{p-\ell-1} \sum_{r=0}^{k -1} 2^{-r}\,  \Gamma_{k,r}^{[5]},
\]
with 
\[
\Gamma_{k,r}^{[5]}= 2^{-2k }  \langle \mu, 
  \cp \left( \cq^{k -r- 1} |\cp (h_{p-\ell- k} \otimes^2)|
  \otimes^2 \right) \rangle.
\]
Using
\reff{eq:mp}  and then  \reff{eq:phk}, we 
 get:
\begin{align*}
 \Gamma_{k,r}^{[5]}
  & \leq C \,  2^{-2k}  \normm{\cp
  (h_{p-\ell- k} \otimes^2)}^2\\ 
& \leq C \, 2^{-2(p-\ell)} \, c_4^4(\bF) \, \ind_{\{k=p-\ell-1\}} + C \, 
  2^{-2k} \alpha^{4(p-\ell -k)} c_2^4(\bF)\,  \ind_{\{k\leq p-\ell-2\}}. 
\end{align*}
We deduce that $\langle \mu,| \psi_{5, p-\ell}| \rangle
  \leq C \, 2^{2(p-\ell)} \, c_4^4(\bF)$.
  
\medskip
\textbf{Upper bound of $\langle \mu, |\psi_{6,p-\ell} |\rangle$}.
We have:
\[
 \langle \mu,| \psi_{6, p-\ell}| \rangle
\leq C\, 2^{3(p-\ell)} \, 
  \sum_{k=1}^{p-\ell-1} \sum_{r=0}^{k -1} 
 2^{-r}\, \Gamma_{k,r}^{[6]},
\]
with 
\[
\Gamma_{k,r}^{[6]}=  2^{-k } \, \langle \mu,  \cq^r \cp \left(
\cq^{k -r-1}|\cp \left( h_{p-\ell-k} \otimes^2 \right)|\sot
\cq^{p-\ell-r-1}(\tilde f^2_\ell) \right) \rangle.
\]
 Using
\reff{eq:mp} and then  \reff{eq:phk}, we 
 get:
\begin{align*}
\Gamma_{k,r}^{[6]}
&\leq  C\,
  2^{-k } \, \normm{\cp \left( h_{p-\ell-k} \otimes^2 \right)} 
\normm{\cq^{p-\ell-r-1}(\tilde f^2_\ell)}\\
&\leq   C\, 2^{-(p-\ell)} \, 
  c_4^4(\bF)\,  \ind_{\{k=p-\ell-1\}}
+   C\,
  2^{-k } \, \alpha^{2(p-\ell-k)}\, c_2^2(\bF)\,
  c_4^2(\bF)\,\ind_{\{k \leq  p-\ell - 2\}}
\end{align*}
We deduce that $\langle \mu,| \psi_{6, p-\ell}| \rangle\leq 
    C\, 2^{2(p-\ell)} \, \, c_4^4(\bF) $.

\medskip
\textbf{Upper bound of $|\langle \mu, \psi_{7,p-\ell} \rangle|$}. 
We have:
\begin{equation*}
   \label{eq:y7}
 |\langle\mu, \psi_{7, p-\ell} \rangle|
\leq C\, 2^{3(p-\ell)} \, 
  \sum_{k=1}^{p-\ell-1} \sum_{r=0}^{k -1} 
 2^{-r}\, \Gamma_{k,r}^{[7]},
\end{equation*}
with 
\begin{equation*}\label{eq:L2R3p7-G}
\Gamma_{k,r}^{[7]}= 2^{-k }  |\langle \mu , \cq^r \cp \left( \cq^{k -r-1}\cp \left( h_{p-\ell-k} \sot  \cq^{p-\ell-k -1} (\tilde f_\ell^2) \right)\sot h_{p-\ell -r} \right)\rangle|.
\end{equation*}
For $k\leq  p-\ell -2$, we have:
\begin{align*}
 \Gamma_{k,r}^{[7]}
&\leq C\,     2^{-k }  \normm{
\cp \left( h_{p-\ell-k} \sot  \cq^{p-\ell-k -1} (\tilde
  f_\ell^2) \right)}\normm{
h_{p-\ell -r} }\\
\nonumber
&\leq  C\, 
  2^{-k }  \normm{ h_{p-\ell-k-1}}\normm{\cq^{p-\ell-k -2} (\tilde
  f_\ell^2)} \, \alpha^{p-\ell-r} c_2(\bF) \ind_{\{k \leq p - \ell - 2\}}\\
\nonumber
&\leq  C\,
  2^{-k }  \alpha ^{2(p-\ell-k)}\, c_2^2(\bF)\, c_4^2(\bF)\,
  \ind_{\{k\leq p-\ell -2\}}, 
\end{align*}
where we used  \reff{eq:mp} for the first  inequality;   \reff{eq:P-L2-1}    for
 the  second; and \reff{eq:hk}  for
the third. We now consider the case  $k=  p-\ell -1$. 
Let $g\in \cb_+(S)$. As $2ba^2\leq b^3+ a^3$ for $a,b$ non-negative, we
get that
$g\otimes  g^2 \leq   g^3\sot \ind$ and thus:
\begin{equation}
   \label{eq:Pgg2}
\cp(g\sot g^2)\leq 2\cq(g^3).
 \end{equation}  
 Writing  $A_r=\Gamma_{p-\ell-1,r}^{[7]}$, we  get using  \reff{eq:Pgg2}
 for the first inequality and Lemma \ref{lem:majo-f3} for the second:
\begin{align*}
   A_r
&= 2^{-p-\ell -1}  |\langle \mu ,  \cp \left(
\cq^{p-\ell -2-r}\cp \left( \tilde f_\ell \sot  \tilde
  f_\ell^2 \right)\sot 
h_{p-\ell -r} \right)\rangle|\\
&\leq C\, 2^{-(p-\ell)}  \langle \mu ,  \cp \left(
\cq^{p-\ell -1-r}|\tilde f_\ell^3| \sot |\cq^{p-\ell -1-r}\tilde f_\ell|
 \right)\rangle\\
&\leq C\, 2^{-(p-\ell)} c_4^4(\bF).
\end{align*}
 Since  $c_2(\bF) \leq c_4(\bF)$, we  deduce that
$ |\langle\mu, \psi_{7, p-\ell} \rangle|\leq  C\, 
2^{2(p-\ell)} \, c_4^4(\bF)$. 

\medskip
\textbf{Upper bound of $\langle \mu, |\psi_{8,p-\ell} |\rangle$}. 
We have:
\begin{equation*}\label{eq:y8}
\langle \mu,| \psi_{8, p-\ell}| \rangle \leq   C\,  2^{4(p-\ell)} \,
\sum_{k=2}^{p-\ell-1} \sum_{r=1}^{k -1} \sum_{j=0}^{r-1} \,  2^{-j }\, \Gamma_{k,r,j}^{[8]},
\end{equation*}
with 
\[
\Gamma_{k,r,j}^{[8]}
\leq  2^{-k-r } \langle \mu,\cq^j\cp \left(|
\cq^{r-j-1}\cp \left( h_{p-\ell-r} \otimes^2 \right)|
\sot
|\cq^{k-j-1}\cp \left( h_{p-\ell-k} \otimes^2 \right)|
\right) \rangle.
\]
Using \reff{eq:mp} and then  \reff{eq:phk}
(twice and noticing that $p-\ell-r\geq 2$), we get:
\begin{align*}
\Gamma_{k,r,j}^{[8]}
&\leq C\, 
 2^{-k-r } \normm{\cp \left( h_{p-\ell-r} \otimes^2 \right)}
\normm{\cp \left( h_{p-\ell-k} \otimes^2 \right)}\\
&\leq C\, 
 2^{-k-r }  \alpha^{2(p-\ell -r)} \, c_2^2(\bF)\, 
\left(\alpha^{2(p-\ell-k)} c_2^2(\bF) + c_4^2(\bF)
  \ind_{\{k=p-\ell-1\}} \right).
\end{align*}
We deduce that $   \langle \mu, |\psi_{8, p-\ell} |\rangle\leq C\,
2^{2(p-\ell)} \,c_4^4(\bF)$.  

\medskip
\textbf{Upper bound of $\langle \mu, |\psi_{9,p-\ell} |\rangle$}. 
We have:
\begin{equation*}\label{eq:L2R3p9}
\langle \mu,| \psi_{9, p-\ell}| \rangle \leq   C\,  2^{4(p-\ell)} \,
\sum_{k=2}^{p-\ell-1} \sum_{r=1}^{k -1} \sum_{j=0}^{r-1} \,  2^{-j } \, \Gamma_{k,r,j}^{[9]},
\end{equation*}
with 
\[\Gamma_{k,r,j}^{[9]}
\leq 
 2^{-k-r } \,
 \langle \mu, \cq^j\cp \left(
\cq^{r-j-1}|\cp \left( h_{p-\ell-r} \sot \cq^{k -r -1}
  \cp\left(h_{p-\ell-k}\otimes^2 \right)\right)|
 \sot |h_{p-\ell-j}| \right)  \rangle .
\]
For $r\leq k-2$, we have:
\begin{align}
   \nonumber
   \Gamma_{k,r,j}^{[9]}
&\leq  C\,  2^{-k-r } \,
 \normm{\cp \left( h_{p-\ell-r} \sot \cq^{k -r -1}
  \cp\left(h_{p-\ell-k}\otimes^2 \right)\right)}
\normm{h_{p-\ell-j}}\\ \nonumber
&\leq   C\,  2^{-k-r } \,
 \normm{ h_{p-\ell-r-1}}\normm{
  \cp\left(h_{p-\ell-k}\otimes^2 \right)}
\normm{h_{p-\ell-j}}\\ \nonumber
&
\leq   C\,  2^{-k-r } \, \alpha^{2(p-\ell-r)}\,  c_2^2(\bF) 
\left(\alpha^{2(p-\ell -k)}\, c_2^2(\bF)\, \ind_{\{ k\leq  p-\ell -2\}}
+ c_4^2(\bF)\, \ind_{\{ k=p-\ell -1\}}
\right),
\end{align}
where we used \reff{eq:mp} for the first inequality; \reff{eq:P-L2-1} as
$p-\ell -r \geq  2$ and $k-r-1\geq 1$ for the  second; and \reff{eq:hk}
(two times) and \eqref{eq:phk} (one time) for the last. For $r= k-1$ and $k\leq p-\ell -2$, we have:
\begin{align}
   \nonumber
   \Gamma_{k,r,j}^{[9]}
&\leq   C\,  2^{-2k } \,
 \normm{\cp \left( h_{p-\ell-k+1} \sot
  \cp\left(h_{p-\ell-k}\otimes^2 \right)\right)}
\normm{h_{p-\ell-j}}\\ \nonumber
&\leq  C\,   2^{-2k } \,
 \normm{ h_{p-\ell-k}}\normm{h_{p-\ell-k-1}}^2
\normm{h_{p-\ell-j}}\\ \nonumber
&\leq  C\,   2^{-2k } \, \alpha^{4(p-\ell -k)}\,  c_2^4(\bF),
\end{align}
where we used \reff{eq:mp} for the first inequality;
\reff{eq:P-L2-3}\footnote{Notice this is the only place in the proof of
 Corollary \ref{cor:subcritical-L2} where we use \reff{eq:P-L2-3}.}  as
$p-\ell -k \geq  2$  for the  second; and \reff{eq:hk}
(three times) for the last. 
For $r= k-1= p-\ell -2$, we have:
\begin{align}
   \nonumber
   \Gamma_{k,r,j}^{[9]}
&\leq   C\,  2^{-2(p-\ell)} \,
 \normm{\cp \left(\cq \tilde f_\ell \sot
  \cp\left( \tilde f_\ell\otimes^2 \right)\right)}
\normm{h_{p-\ell-j}}\\ \nonumber
&\leq   C\,  2^{-2(p-\ell) } \,
 \normm{\cp \left( \cq \tilde f_\ell \sot \cq(\tilde f_\ell^2)\right)}
\normm{h_{p-\ell-j}}\\ \nonumber
&\leq   C\,  2^{-2(p-\ell) } \, c_4^2(\bF)\, \alpha^{p-\ell -j}\,
  c_2^2(\bF), 
\end{align}
where we used  \reff{eq:mp} for the first inequality,
\reff{eq:majo-pfxf} (with $f$ replaced by $f_\ell$) for the
second and \reff{eq:P-L2-1} as well as \reff{eq:phk} (with $p-\ell -j
\geq 2$) for the last. Taking all together, we deduce that $
\langle \mu,| \psi_{9, p-\ell}| \rangle
\leq   C\, 2^{2(p-\ell)}\, c_4^2(\bF)\, 
  c_2^2(\bF)$. 

\medskip
Wrapping all the upper bounds with \reff{eq:R3L2-Em4} we deduce 
that for $\ell\in \{0, \ldots,
 p-3\} $
\[
\E_\mu\left[M_{\G_{p-\ell}}(\tilde f_\ell) ^4\right]
\leq  C\, 2^{2(p-\ell)} c_4^4(\bF).
\]
Thanks to \reff{eq:L2R3-l=p}, this equality holds  for $\ell\in \{0,
\ldots, p\}$. We deduce from \reff{eq:R3-L2-9} that:
\begin{equation}
   \label{eq:R3-31}
R_3(n) \leq   C\,
 n^3\, 2^{-(n-p)}\,c_4^4(\bF) . 
 \end{equation} 
This proves that $\lim_{n \rightarrow \infty} R_{3}(n) = 0$.
\end{proof}

We can  now use  Theorem 3.2 and  Corollary 3.1, p.~58, and  the Remark
p.~59  from  \cite{hh:ml}  to  deduce  from  Lemmas \ref{lem:cvV-L2}
and \ref{lem:cvG-L2} that  $\Delta_n(\bF)$ converges in distribution
towards a Gaussian  real-valued   random  variable  with   deterministic
variance $\ssub(\bF)$ given  by \reff{eq:ssub}.  Using  \reff{eq:N=D+R},
Remark \ref{rem:Nn=Nnk0} and Lemmas  \ref{lem:cvR0-L2}  and  \ref{lem:cvR1-L2}, we  then  deduce  Theorem \ref{cor:subcritical-L2}.

\section{Moments formula for BMC}
\label{sec:moment}
Let $X=(X_i,  i\in \T)$ be a  BMC on $(S, \cs)$  with probability kernel
$\cp$.          Recall           that          $|\G_n|=2^n$          and
$M_{\G_n}(f)=\sum_{i\in \G_n} f(X_i)$. We also recall that
$2\cq(x,A)=\cp(x, A\times S) + \cp(x, S\times A)$ for $A\in \cs$. 
We use the convention that $\sum_\emptyset=0$.
\medskip

We recall the following well known and easy to establish  many-to-one
formulas  for BMC.

\begin{lem}
   \label{lem:Qi}
Let $f,g\in \cb(S)$, $x\in S$ and $n\geq m\geq 0$. Assuming that all the
quantities below are well defined,  we have:
\begin{align}
   \label{eq:Q1}
\E_x\left[M_{\G_n}(f)\right]
&=|\G_n|\, \cq^n f(x)= 2^n\, \cq^n f(x) ,\\
   \label{eq:Q2}
\E_x\left[M_{\G_n}(f)^2\right]
&=2^n\, \cq^n (f^2) (x) + 
 \sum_{k=0}^{n-1} 2^{n+k}\,   \cq^{n-k-1}\left( \cp
  \left(\cq^{k}f\otimes 
    \cq^k f \right)\right) (x),\\
   \label{eq:Q2-bis}
\E_x\left[M_{\G_n}(f)M_{\G_m}(g)\right]
&=2^{n} \cq^{m} \left(g \cq^{n-m} f\right)(x)\\
  \nonumber &\hspace{2cm} + \sum_{k=0}^{m-1} 2^{n+k}\, \cq^{m-k-1}
  \left(\cp\left(\cq^k g \sot \cq^{n-m+k} f\right) \right)(x). 
\end{align}
\end{lem}

We also give some bounds on $\E_x\left[M_{\G_{n}}(f) ^4\right]$,  
see the proof of Theorem 2.1 in \cite{BDG14}.  We will use the notation:
\[
g\otimes^2=g\otimes g.
\]
\begin{lem}
   \label{lem:M4}
   There exists a  finite constant $C$ such that for  all $f\in \cb(S)$,
   $n\in \N$ and  $\nu$ a probability measure on $S$,  assuming that all
   the quantities below are well defined, there exist functions $\psi_{j,
     n}$ for $1\leq j\leq 9$ such that:
\[
\E_\nu\left[M_{\G_n}(f)^4\right]= \sum_{j=1}^9 \langle \nu, \psi_{j, n}
\rangle,
\]
and,   with
$h_{k}= \cq^{k - 1} (f) $ and (notice that either $|\psi_j|$ or
$|\langle \nu, 
\psi_j \rangle|$ is bounded), writing  $\nu   g=\langle
\nu  ,   g   \rangle$:
\begin{align*}
 | \psi_{1, n}|
&\leq C \,2^n \cq^n(f^4),\\
 | \nu \psi_{2, n}|
&\leq C\,  2^{2n}\,
 \sum_{k=0}^{n-1} 2^{-k} |\nu \cq^k
 \cp  \left( \cq^{n-k - 1}( f^3) \sot h_{n- k} \right)|,\\
|\psi_{3, n}|
&\leq C 2^{2n} \sum_{k=0}^{n-1} 2^{-k}\,  \cq^k
  \cp \left( \cq^{n-k - 1} (f^2) \otimes^2
  \right),\\
|\psi_{4, n}|
&\leq C \, 2^{4n} \, 
\cp \left( |\cp(h_{n-1}\otimes^2)\otimes^2|\right), \\
|\psi_{5, n}|
&\leq C\,  2^{4n} \, 
 \sum_{k=2}^{n-1} \sum_{r=0}^{k -1}  2^{-2k-r }  \cq^r
  \cp \left( \cq^{k -r- 1} |\cp (h_{n- k} \otimes^2)|
  \otimes^2 \right),\\
|\psi_{6, n}|
&\leq C\, 2^{3n} \, 
  \sum_{k=1}^{n-1} \sum_{r=0}^{k -1} 
  2^{-k-r }  \cq^r| \cp \left(
\cq^{k -r-1}\cp \left( h_{n-k} \otimes^2 \right)\sot
\cq^{n-r-1}(f^2) \right)|,\\
|\nu \psi_{7, n}|
&\leq  C\, 2^{3n} \,   \sum_{k=1}^{n-1} \sum_{r=0}^{k -1} 
  2^{-k-r } |\nu \cq^r \cp \left(
\cq^{k -r-1}\cp \left( h_{n-k} \sot  \cq^{n-k -1} (f^2) \right)\sot
h_{n-r} \right)|,\\
|\psi_{8, n}|
&\leq C\,  2^{4n} \, 
  \sum_{k=2}^{n-1} \sum_{r=1}^{k -1} \sum_{j=0}^{r-1}
  2^{-k-r-j } \cq^j \cp \left(|
\cq^{r-j-1}\cp \left( h_{n-r} \otimes^2 \right)|\sot
|\cq^{k-j-1}\cp \left( h_{n-k} \otimes^2 \right)|
\right),\\
|\psi_{9, n}|
& \leq C\,  2^{4n} \, 
  \sum_{k=2}^{n-1} \sum_{r=1}^{k -1} \sum_{j=0}^{r-1}
  2^{-k-r-j } \cq^j |\cp \left(
\cq^{r-j-1}|\cp \left( h_{n-r} \sot \cq^{k -r -1}
  \cp\left(h_{n-k}\otimes^2 \right)\right) \sot h_{n-j} 
\right)|.
\end{align*}
\end{lem}
We shall use the following lemma in order to bound the term $ | \nu \psi_{2, n}|$. 

\begin{lem}
   \label{lem:majo-f3}
Let $\mu$ be an invariant probability measure on $S$ for $\cq$. Let $f,
g\in L^4(\mu)$. Then we have for all $r \in \N$:
\begin{equation*}
\langle \mu , \cp( \cq^r |f|^3 \otimes  |g|) \rangle
\leq  2 \norm{f}_{L^4(\mu)}^3 \, \norm{g}_{L^4(\mu)}.
\end{equation*}
\end{lem}
\begin{proof}
We have
\begin{align*}
 \langle \mu , \cp( \cq^r |f|^3 \otimes  |g|) \rangle
&\leq    \langle \mu , \cp( (\cq^r |f|^3)^{4/3} \otimes 1) \rangle^{3/4} \, 
\langle\mu , \cp( 1 \otimes  g^4) \rangle^{1/4} \\
&\leq   2 \langle \mu , \cq ( (\cq^r |f|^3)^{4/3} ) \rangle^{3/4} \, 
\langle\mu , \cq (  g ^4) )\rangle^{1/4} \\
&\leq   2 \langle \mu ,  |f|^4  \rangle^{3/4} \, 
\langle\mu , |g|)^4 \rangle^{1/4},
\end{align*}
where we used H\"older inequality and that  $v\otimes w=(v\otimes 1) \,
(1\otimes w)$   for the first inequality, that 
$\cp(v\otimes 1) \leq  2 \cq v$ and 
$\cp(1\otimes v) \leq  2 \cq v$ if $v$ is non-negative for the second
inequality, Jensen's inequality and that $\mu$ is invariant for $\cq$
for the last. 
\end{proof}

\bibliographystyle{abbrv}
\bibliography{biblio}

\newpage

\section{Supplementary material  to Section  \ref{sec:critical} on
  the critical case}
We give a proof to Theorem \ref{cor:critical-L2}. 
We keep notations from Section \ref{sec:proof-sub-cor-L2} on the sub-critical case, and adapt very closely the arguments of this section. We recall that $c_{k}(\bF)=\sup\{\norm {f_n}_{L^{k}(\mu)}, \, n\in \N\}$ for all $k \in \NN$. We recall that $C$  denotes any unimportant finite constant which may vary from line to line, which does not depend on $n$ or  $\bF$.

\begin{lem}\label{lem:cvRk0L2crit}
Under the assumptions of Theorem \ref{cor:critical-L2}, we have that $\lim_{n \rightarrow \infty} \EE[n^{-1}R_{0}^{k_{0}}(n)^{2}] = 0.$
\end{lem}

\begin{proof}
Mimicking  the proof of Lemma \ref{lem:cvR0-L2}, we get:
\begin{equation*}
\lim_{n \rightarrow \infty} \EE[R_{0}^{k_{0}}(n)^{2}]^{1/2} \leq \lim_{n
  \rightarrow \infty} C c_{2}(\bF)\, \sqrt{n}  2^{-p/2} = 0.
\end{equation*}
This trivially implies the result. 
\end{proof}

\begin{lem}\label{lem:cvR1L2crit}
Under the assumptions of Theorem \ref{cor:critical-L2}, we have that
$\lim_{n \rightarrow \infty} \EE[n^{-1} R_{1}(n)^{2}] = 0.$  
\end{lem}

\begin{proof}
Mimicking the proof of Lemma \ref{lem:cvR1-L2}, we get 
$\EE[R_{1}(n)^{2}]^{1/2} \leq C c_{2}(\bF)\sqrt{n-p}$. As
$\lim_{n\rightarrow\infty } p/n=1$, this 
implies that $\lim_{n \rightarrow \infty}\EE[n^{-1}R_{1}(n)^{2}] = 0$.
\end{proof}

Similarly to Lemma \ref{lem:cvR2-L2}, we get the following result on
$R_2(n)$. 

\begin{lem}\label{lem:cvR2L2crit}
Under the assumptions of Theorem \ref{cor:critical-L2}, we have that $\lim_{n \rightarrow \infty} \EE[n^{-1/2} R_{2}(n)] = 0.$
\end{lem}

We now consider the asymptotics of $V_2(n)$. 
\begin{lem}\label{lem:cvL2V2-crit}
Under the assumptions of Theorem \ref{cor:critical-L2}, we have that
 $\lim_{n \rightarrow \infty} n^{-1}V_{2}(n) = \scrit
_2(\bF)$ in probability, where 
 $\scrit _2(\bF)$, defined in \eqref{eq:S2-crit},  is well defined
 and finite. 
\end{lem}
In the proof, we shall use the analogue of \reff{eq:P-L2-2} with $f$ replaced by
$\hat f$ in the left hand-side, whereas $f\in L^4(\mu)$ does imply that
$\tilde f \in L^4(\mu)$ but does not imply
that $\hat f \in L^4(\mu)$. 
 Thanks to \eqref{eq:P-L2-2}, we get  for $f\in L^4(\mu)$ and $g\in
L^2(\mu)$,  
as $\crr_j f= \alpha_j^{-1} \cq \crr_j f$ and $|\alpha_j|=\alpha$, that:
\begin{align}
\nonumber
\normm{ \cp \left(\hat f \sot \cq g\right)}
&\leq  \normm{\cp \left( \tilde f \sot \cq g\right)}+ \alpha^{-1} 
\sum_{j\in J} \normm{ \cp \left(\cq (\crr_j f) \sot \cq g\right)}\\
&\nonumber 
\leq  C \,\left( \norm{f}_{L^4(\mu)} +  \normm{f}\right)\normm{g}\\
&\label{eq:majocp-hatf-qg} 
\leq  C \,\norm{f}_{L^4(\mu)}\normm{g}.
\end{align}

\begin{proof} 
We keep the decomposition \reff{eq:decom-V2} of $V_2(n)=V_5(n)+V_6(n)$ given in the proof of Lemma \ref{lem:cvV2-L2}. 
We recall $V_{6}(n)= |\G_{n-p}|^{-1} M_{\G_{n-p}} (H_{6, n})$ with $H_{6, n}$ defined in \reff{eq:H6n-crit-det}. We set
\begin{equation*}
\bar{H}_{6,n} = \sum_{0\leq \ell<k\leq  p; \, r\geq 0}  \bar h_{k,\ell,r}^{(n)}\,  \ind_{\{r+k<  p\}} \quad \text{and} \quad \bar{V}_{6}(n) = |\GG_{n-p}|^{-1} M_{\GG_{n-p}}(\bar{H}_{6,n}),
\end{equation*}
where  for $0\leq \ell<k\leq p$ and $0\leq r<p-k$: 
\[
\bar  h_{k,\ell,r}^{(n)} =  2^{r-\ell} \,  \alpha^{k-\ell +2r} \,
\cq^{p-1-(r+k)}(\cp f_{k, \ell,r})
=  2^{-(k+\ell)/2}\, \cq^{p-1-(r+k)}(\cp f_{k, \ell,r}), 
\]
where  we used  that  $2\alpha^2=1$. For $f \in L^{2}(\mu)$, we recall $\hat{f}$ defined in \eqref{eq:fhatcrit-S}. We set: 
\begin{align*}
h_{k,\ell,r}^{(n,1)} 
&= 2^{r - \ell }\Qq^{p-1-(r+k)}(\Pp(\Qq^{r}(\hat{f}_{k}) \sot
  \Qq^{k-\ell+r}(\hat{f}_{\ell}))), \\ 
h_{k,\ell,r}^{(n,2)}
& = 2^{r - \ell}\Qq^{p-1-(r+k)}(\Pp(\Qq^{r}(\hat{f}_{k}) \sot \Qq^{k
  - \ell + r}(\sum_{j \in J} \crr_{j}(f_{\ell})))), \\ 
h_{k,\ell,r}^{(n,3)} 
&= 2^{r - \ell }\Qq^{p-1-(r+k)}(\Pp(\Qq^{r}(\sum_{j \in J}
  \crr_{j}(f_{k})) \sot \Qq^{k-\ell+r}(\hat{f}_{\ell}))),
\end{align*}
so that $h^{(n)}_{k, \ell, r}= \bar h^{(n)}_{k, \ell, r} + \sum_{i=1}^3 
h^{(n,i)}_{k, \ell, r}$. 
Thanks to \eqref{eq:P-L2-1} for $r \geq 1$ and
\reff{eq:majocp-hatf-qg}  for $r = 0$, we have using Jensen's inequality, \eqref{eq:L2-erg-crit} and the fact that the sequence $(\beta_{r}, r \in \NN)$ is nonincreasing:
\begin{equation*}
\normm{h_{k,\ell,r}^{(n,1)}} \leq C 2^{-(k + \ell)/2} \beta_{ r} \normm{f_{\ell}} 
\begin{cases} \norm{f_{k}}_{L^{2}(\mu)} & \text{for $r \geq 1$} ,\\ 
\norm{f_{k}}_{L^{4}(\mu)}  & \text{for $r = 0.$}\end{cases}
\end{equation*}
Using the same arguments,  that $\langle \mu, \crr_j(g)
\rangle=0$ for $g\in L^2(\mu)$ (as $\crr_j(g)$ is an eigen-vector of
$\cq$ associated to $\alpha_j$)   and that  $\normm{\sum_{j \in J} \crr_{j}(f_{\ell})} \leq
C\, \normm{f_{\ell}}$ (as $\crr_j$ are bounded operators on $L^2(\mu)$), we get:
\begin{equation*}
\normm{h_{k,\ell,r}^{(n,2)}} + \normm{h_{k,\ell,r}^{(n,3)}} \leq C 2^{-(k + \ell)/2}
\beta_{r}\, \normm{f_\ell} 
\begin{cases} \normm{f_{k}} 
& \text{for $r \geq 1$}, \\ 
\norm{f_{k}}_{L^{4}(\mu)} 
& \text{for $r = 0.$} 
\end{cases}
\end{equation*}
We deduce that
\begin{equation}\label{eq:I-h1nklr}
\sum_{i=1} ^3 \normm{h_{k,\ell,r}^{(n,i)}} \leq C c_{2}(\bF) c_{4}(\bF)
2^{-(k + \ell)/2} \beta_{r}. 
\end{equation}

Using \eqref{eq:EGn-nu0}  for the first inequality,  Jensen's inequality
for  the second  inequality,  the triangular  inequality  for the  third
inequality and \eqref{eq:I-h1nklr} for the last inequality, we get:
\begin{align}
\nonumber \E\left[\left(V_6(n) - \bar V_6(n)\right)^2\right]
&=|\GG_{n-p}|^{-2}\EE[M_{\GG_{n-p}}(H_6(n) - \bar H_6(n))^{2}] \\
&\leq C |\GG_{n-p}|^{-1} \sum_{m = 0}^{n-p} 2^{m}
 \normm{\cq^{m} (H_6(n) - \bar H_6(n)) }^2 \nonumber   \\ 
&\leq C \normm{H_6(n) - \bar H_6(n)} ^{2} \nonumber \\ 
&\leq C \Big(\sum_{0 \leq \ell < k < p} \sum_{r = 0}^{p - k - 1}
\sum_{i=1}^3 \normm{h_{n,k,\ell,r}^{(n,i)}}\Big)^{2} \nonumber  
\\ 
&\leq C c_{2}(\bF)^{2}c_{4}(\bF)^{2} \Big(\sum_{r=0}^{p} \beta_{r}\Big)^{2}. \nonumber
\end{align} 
We deduce that
\begin{equation*}
\EE[(V_{6}(n) - \bar{V}_{6}(n))^{2}] \leq C c_{2}(\bF)^{2}c_{4}(\bF)^{2}
\Big(\sum_{r=0}^{p} \beta_{r}\Big) ^2 , 
\end{equation*}
and then that
\begin{equation}
\label{eq:vbar6-crit-L2}
\lim_{n\rightarrow\infty} \EE[n^{-2}(V_{6}(n) - \bar{V}_{6}(n))^{2}] = 0.
\end{equation}
We set $  H_6^{[n]} = \sum_{0\leq \ell<k\leq  p; \, r\geq 0} h_{k,\ell,r}\,  \ind_{\{r+k<  p\}} $ with  for $0\leq \ell<k\leq p$ and $0\leq r<p-k$: 
\begin{equation*}\label{eq:hat-crit}
 h_{k, \ell,r} = 2^{-(k+\ell)/2} \langle \mu, \cp f_{k, \ell,r} \rangle = \langle \mu, \bar h_{k, \ell, r}^{(n)}  \rangle.
\end{equation*}
We have that
\begin{equation*}\label{eq:hatH6L2}
{H}_{6}^{[n]} = \sum_{0 \leq \ell < k < p} \sum_{r = 0}^{p-k-1} 
h_{k, \ell,r}
= \langle \mu, \bar H_{6, n} \rangle. 
\end{equation*}
We have:
\begin{align}
\EE[(\bar{V}_{6}(n) - {H}_{6}^{[n]})^{2}] 
&\leq C |\GG_{n-p}|^{-1}  \sum_{m = 0}^{n-p} 2^{m} \normm{\Qq^{m}
  (\bar{H}_{6,n} - {H}_{6}^{[n]})}^{2} \nonumber\\ 
&\leq C |\GG_{n-p}|^{-1}  \sum_{m = 0}^{n-p} 2^{m} 
\left(\sum_{0 \leq \ell < k \leq p} \sum_{r =
  0}^{p-k-1}  \alpha^{m+p- r -k } 2^{-(k +\ell)/2} \normm{\cp f_{k, \ell, r}}
\right)^2\nonumber\\
& \leq C (n-p) |\GG_{n-p}|^{-1} \left(\sum_{0 \leq \ell < k \leq p} \sum_{r =
  0}^{p-k-1} 2^{-(p+\ell -r)/2} 
  \normm{\Pp(f_{k,\ell,r})}
\right)^{2} \nonumber 
\\ 
& \leq C (n-p) |\GG_{n-p}|^{-1}  \left(\sum_{0 \leq \ell < k < p}
2^{-(\ell + k)/2}  \|\sum_{j \in J}
  \crr_{j}(f_{k})\|_{L^{2}(\mu)} \|\sum_{j \in J}
  \crr_{j}(f_{\ell})\|_{L^{2}(\mu)}
\right)^{2} \nonumber\\ 
&\leq C(n-p)|\GG_{n-p}|^{-1} \,  c_2^4(\bF), \nonumber
\end{align}
where we used \eqref{eq:EGn-nu0} for the first inequality,
\eqref{eq:L2-erg}  for the second, $\alpha = 1/\sqrt{2}$ for
the third,  
\eqref{eq:P-L2-1} and the fact that $\Qq(\sum_{j \in J} \crr_{j} f) =
\sum_{j \in J}\alpha_{j} \crr_{j}(f)$, with $|\alpha_{j}| = 1/\sqrt{2}$,
for the fourth, $\|\sum_{j \in J} \crr_{j}(f)\|_{L^{2}(\mu)}
\leq \normm{f}$ for the 
last.  From the latter inequality we conclude that: 
\begin{equation}
\label{eq:V6bar-crit-L2}
\lim_{n \rightarrow \infty} \EE[n^{-2}(\bar{V}_{6}(n) - {H}_{6}^{[n]})^{2}]  = 0.
\end{equation}
We set for $k, \ell\in \N$: $h_{k,\ell}^*= 2^{-(k+\ell)/2} \langle \mu, \cp (f^*_{k, \ell} ) \rangle$ and we consider the sums
\begin{equation*}\label{eq:HstarV6}
H_0^*= \sum_{0\leq \ell <k} (k+1) |h_{k,\ell}^*|  \quad \text{and} \quad H^{*}_{6}(\bF) =  \sum_{0\leq \ell <k} h_{k,\ell}^* = \scrit_2(\bF). 
\end{equation*}
Using \eqref{eq:mp}, we have:
\begin{equation*}
|h_{k,\ell}^{*}| \leq C 2^{-(k+\ell)/2} \sum_{j \in J} \normm{\crr_{j}(f_{k})} \normm{\crr_{j}(f_{\ell})} \leq C2^{-(k+\ell)/2}  c_{2}^2(\bF).
\end{equation*}
This implies that $H_{0}^{*} \leq C c_{2}^2(\bF)$, $H_{6}^{*}(\bF) \leq C c_{2}^2(\bF)$ and then that $H_{0}^{*}$ and
$H_{6}^{*}(\bF)$ are well defined.
We write:
\begin{equation*}
 h_{ k,\ell, r}= h_{k,\ell}^* +  h_{k,\ell, r}^{\circ}, \quad \text{with} \quad h_{k,\ell,r}^{\circ} = 2^{-(k+\ell)/2} \langle \mu, \cp f^\circ_{k, \ell,r} \rangle,
\end{equation*}
where we recall that $ f^\circ_{k, \ell, r}= f_{k, \ell, r}-
f^*_{k, \ell}$, 
and 
\begin{equation}\label{eq:decomphatHn}
H_{6}^{[n]}  = H_{6}^{[n],*} +H_{6}^{[n],\circ} 
\end{equation}
with  
\begin{equation*}\label{eq:Hnhat-H0nhat}
H_{6}^{[n],*} = \sum_{0\leq \ell<k\leq  p} (p-k) h^*_{k, \ell} \quad
\text{and} \quad 
H_{6}^{[n],\circ} = \sum_{0\leq \ell<k\leq  p; \, r\geq 0}  
h_{k,\ell,r}^\circ \,  \ind_{\{r+k<  p\}}.
\end{equation*} 
Recall $\lim_{n\rightarrow \infty } p/n=1$. We have:
\[
|n^{-1} H_{6}^{[n],*} - H^{*}_{6}(\bF)|\leq   |n^{-1} p - 1| |H^{*}_{6}(\bF)| + n^{-1} H_0^* +\sum_{\substack{0\leq  \ell<k \\ k>p}} |h_{k, \ell}^*| ,
\]
so that $\lim_{n\rightarrow \infty } |n^{-1}  H_{6}^{[n],*}- H^{*}_{6}(\bF)| = 0$ and
thus:
\begin{equation}\label{eq:limtH-hH}
\lim_{n\rightarrow \infty } n^{-1}  H_{6}^{[n], *} = H^{*}_{6}(\bF).
\end{equation}
We now prove that  $n^{-1}  H_{6}^{[n], \circ}$ converges towards 0. We have:
\begin{equation}\label{eq:detail-f0}
f^\circ_{k, \ell, r}= \sum_{j, j'\in J, \, \theta_j \theta_{j'}\neq 1} (\theta_{j'}\theta_{j})^r \theta_{j'}^{k-\ell} \,\, \crr_{j} f_k \sot \crr_{j'} f_\ell.
\end{equation}
This gives:
\begin{align}
| H^{[n],\circ}_{6}| 
& =\Big|  \sum_{0\leq \ell<k\leq  p, \, r\geq 0}  2^{ -(k+\ell)/2}
  \langle \mu, \cp f^\circ_{k, \ell, r}  \rangle\,  \ind_{\{r+k<
  p\}}\Big| \nonumber \\ 
& \leq  \sum_{0\leq \ell<k\leq  p}  2^{ -(k+\ell)/2}
  \sum_{j, j'\in J, \, \theta_j \theta_{j'}\neq 1} \Big| \langle \mu,
  \cp( \crr_{j} f_k \sot \crr_{j'} f_\ell) \rangle \Big| \,
  \Big|\sum_{r=0}^{p-k-1} (\theta_{j'}\theta_{j})^r
  \Big|, \label{eq:2IhatH0n} 
\end{align}
where we used \reff{eq:detail-f0} for the inequality. 
Using  \reff{eq:mp} in the upper bound \reff{eq:2IhatH0n}, we get
\begin{equation*}
\Big| \langle \mu, \cp( \crr_{j'} f_k \sot \crr_j f_\ell) \rangle \Big| \leq 2\normm{\crr_{j'}(f_k)}\normm{\crr_j(f_\ell)} \leq C\normm{f_k}\normm{f_\ell}.
\end{equation*}
This implies that $| H^{[n],\circ}_{6}| \leq c$, with
\begin{equation*}
c = C \, c_{2}(\bF)^{2} \, \sum_{0\leq \ell<k\leq  p}  2^{ -(k+\ell)/2} \sum_{j, j'\in J, \, \theta_j \theta_{j'}\neq 1}  |1- \theta_{j'}\theta_{j}|^{-1}.
\end{equation*}
Since  $J$ is
finite, we deduce that $c$ is finite. This gives that
$\lim_{n\rightarrow \infty } n^{-1}  H^{[n],\circ} _{6} = 0$. Recall that
$ H_{6}^{[n]}$ and $H^{*}_{6}(\bF)$ are complex numbers (\textit{i.e.} constant
functions). Use
\reff{eq:decomphatHn} and \reff{eq:limtH-hH} to get that: 
\begin{equation}\label{eq:limhatHnV6}
\lim_{n\rightarrow \infty } n^{-1}  H_{6}^{[n]} = H^{*}_{6}(\bF) 
\end{equation}
It follows  from \reff{eq:vbar6-crit-L2}, \reff{eq:V6bar-crit-L2}
and \reff{eq:limhatHnV6} 
that:
\begin{equation}\label{eq:LV6ncrit}
\lim_{n \rightarrow \infty} \EE[(n^{-1}V_{6}(n) - H_{6}^{*}(\bF))^{2}] = 0.
\end{equation} 
\medskip

We recall $H_{5}^{[n]}(\bF)$ defined in \eqref{eq:def-B5nsubL2}. From
\eqref{eq:V5-H5f}, we have:  
\begin{equation*}
\EE[n^{-2}V_{5}(n)^{2}] \leq 2 n^{-2} |\G_{n-p}|^{-2}\, \E\left[M_{\G_{n-p}} (A_{5,n}(\bF))^2 \right] + 2 n^{-2} H_{5}^{[n]}(\bF)^{2}.
\end{equation*}
Using \eqref{eq:majo-H5h} with $\alpha = 1/\sqrt{2}$, we get $|H_5^{[n]}
(\bF)|\leq  C \, c_2^2(\bF)$ and thus:
\begin{equation*}
\lim_{n \rightarrow \infty} n^{-2} H_{5}^{[n]}(\bF)^{2} = 0.
\end{equation*} 
Next, as \reff{eq:L2A5Hf} holds for $\alpha=1/\sqrt{2}$, we get  \reff{eq:majo-L2A5}
with the
right hand-side replaced by $ C\,c_4^4(\bF)\, (n-p)2^{-(n-p)} $, and
thus: 
\begin{equation*}
\lim_{n \rightarrow \infty} n^{-2} |\G_{n-p}|^{-2}\, \E\left[M_{\G_{n-p}} (A_{5,n}(\bF))^2 \right] = 0.
\end{equation*}
It then follows that:
\begin{equation*}
\lim_{n \rightarrow \infty} \EE[n^{-2}V_{5}(n)^{2}] = 0.
\end{equation*}
Finally, since $V_{2}(n) = V_{5}(n) + V_{6}(n)$, we get thanks to \reff{eq:HstarV6} 
that in probability $\lim_{n\rightarrow \infty } n^{-1}
V_2(n)=H^*_6(\bF)=\scrit_2(\bF)$. 
\end{proof}

\begin{lem}\label{lem:cvV1L2-crit}
Under the assumptions of Theorem \ref{cor:critical-L2},  we have
that in probability $\lim_{n\rightarrow \infty } V_1(n)
=\scrit_1(\bF)$, where $\scrit_1(\bF)$, defined in \eqref{eq:S1-crit},
 is well defined and finite. 
\end{lem}

\begin{proof}
We recall the decomposition \eqref{eq:DV4nsub}: $V_{1}(n) = V_{3}(n) +
V_{4}(n).$ First, following the proof of \eqref{eq:LV6ncrit} in the
spirit of the proof of \reff{eq:majo-L2A4}, 
we get:
\begin{equation*}
\lim_{n \rightarrow \infty}\EE[(n^{-1}V_{4}(n) - H^{*}_{4}(\bF))^{2}] = 0 \quad \text{with} \quad H_{4}^{*}(\bF) = \sum_{\ell\geq 0}2^{-\ell} \langle \mu, \cp( \sum_{j\in J}\,\,  \crr_j (f_{\ell})\sot \overline \crr_j (f_{\ell}))\rangle=\scrit_1(\bF).
\end{equation*}
Let us stress that the proof requires to use \reff{eq:mp2}. Since
$\sum_{\ell\geq 0}2^{-\ell} |\langle \mu, \cp( \sum_{j\in J}\,\,  \crr_j
(f_{\ell})\sot \overline \crr_j (f_{\ell}))\rangle|
\leq  \sum_{\ell\geq 0}2^{-\ell}  c_2^2(\bF)$, we deduce that $\scrit_1(\bF)$
is well defined and finite. 
Next, from \eqref{eq:V3-H3f} we have
\begin{equation*}
\EE[n^{-2}V_{3}(n)^{2}] \leq 2 n^{-2} |\G_{n-p}|^{-2}\, \E\left[M_{\G_{n-p}} (A_{3,n}(\bF))^2 \right] + 2 n^{-2} H_{3}^{[n]}(\bF)^{2}.
\end{equation*} 
It follows from \eqref{eq:majo-L2A3} (with an extra term $n-p$ as
$2\alpha^2=1$ in the right hand side) and \eqref{eq:limH3n} that
$\lim_{n \rightarrow \infty}\EE[n^{-2} V_{3}(n)^{2}] = 0.$ Finally the
result of the lemma follows as $V_1=V_3+V_4$. 
\end{proof}

We now check the Lindeberg condition using a fourth moment
condition. Recall $R_3(n)=\sum_{i\in \G_{n-p_n}}
\E\left[\Delta_{n,i}(\bF)^4\right]$ defined in \reff{eq:def-R3}.

\begin{lem}\label{lem:cvG-L2-b}
Under the assumptions of Theorem \ref{cor:critical-L2},  we have that $\lim_{n\rightarrow\infty } n^{-2}R_3(n)=0.$
\end{lem}
\begin{proof}
  Following line  by line the  proof of Lemma \ref{lem:cvG-L2}  with the
  same  notations  and  taking  $\alpha  = 1/  \sqrt{2}$,  we  get  that
  concerning    $|\langle    \mu,     \psi_{i,p-\ell}    \rangle|$    or
  $\langle   \mu,    |\psi_{i,p-\ell}|   \rangle$,   the    bounds   for
  $i\in \{1, 2, 3,  4\}$ are the same; the bounds  for $i\in \{5, 6,7\}$
  have an extra $(p-\ell)$ term, the  bounds for $i\in \{8, 9\}$ have an
  extra $(p-\ell)^2$ term. This leads to (compare with \reff{eq:R3-31}):
\[
R_3(n) \leq   C\, n^5\, 2^{-(n-p)}\,c_4^4(\bF) 
\]
which implies that $\lim_{n\rightarrow\infty } n^{-2}R_3(n)=0$.
\end{proof}

The proof of  Theorem \ref{cor:critical-L2} then follows the proof of
Theorem \ref{cor:subcritical-L2}.

\section{Supplementary material  to Section  \ref{sec:main-res-super} on
  the supercritical case}
\subsection{Complementary results and proof of Corollary \ref{cor:super-crit-d1-2}}

Now, we state the main result of this section, whose proof is given in
Section \ref{sec:proof-super-thm}. Recall that
$\theta_j=\alpha_j/\alpha$ and $|\theta_j|=1$ and $M_{\infty , j}$ is
defined in Lemma \ref{lem:martingale}. 
\begin{theo}\label{theo:super-critical}
Let $X$ be a BMC with kernel $\cp$ and initial distribution
$\nu$ such  that Assumptions \ref{hyp:Q} (ii) and \ref{hyp:Q2} are in
force with
$\alpha\in (1/\sqrt{2}, 1)$ in
\reff{eq:L2-erg-crit}. We  have the  following convergence for all sequence $\bF=(f_\ell,   \ell\in   \N)$ uniformly bounded  in   $L^2(\mu)$   (that is $\sup_{\ell\in \N} \norm{f_\ell}_{L^2(\mu)}<+\infty $):
\[
(2\alpha^{2})^{-n/2} N_{n, \emptyset}(\bF)
- \sum_{\ell\in \N} (2\alpha)^{-\ell} \sum_{j\in J} \theta_j^{n-\ell}
M_{\infty  ,j}(f_\ell)  
\; \xrightarrow[n\rightarrow \infty ]{\P}  \; 0. 
\]
\end{theo}
\begin{rem}\label{rem:theo-Scrit}
We stress that if for all $\ell \in \NN$, the orthogonal projection of $f_{\ell}$ on the eigen-spaces corresponding to the eigenvalues $1$ and $\alpha_{j}$, $j \in J$, equal 0, then $M_{\infty,j}(f_{\ell}) = 0$ for all $j\in J$ and in this case, we have
\begin{equation*}
(2\alpha^{2})^{-n/2} N_{n, \emptyset}(\bF) \; \xrightarrow[n\rightarrow \infty ]{\P}  \; 0. 
\end{equation*}
\end{rem}
 
As a direct consequence  of Theorem \ref{theo:super-critical} and Remark
\ref{rem:simpleN0n},  we  deduce  the  following  results.  Recall  that
$\tilde f= f -\langle \mu, f \rangle$.
\begin{cor}\label{cor:CtheoSC}
Under the assumptions of Theorem \ref{theo:super-critical}, we have for all $f \in L^{2}(\mu)$:
\begin{align*}
(2\alpha)^{-n} M_{\TT_{n}}(\tilde{f}) - \sum_{j\in J} \theta_{j}^{n}(1 -
  (2\alpha\theta_{j})^{-1})^{-1}M_{\infty,j}(f) 
& \inP 0 \\
(2\alpha)^{-n} M_{\GG_{n}}(\tilde{f}) - \sum_{j\in J} \theta_{j}^{n}
  M_{\infty,j}(f) 
& \inP 0.
\end{align*}
\end{cor}
\begin{proof}
  We first  take $\bF = (f,f,\ldots)$  and next $\bF =  (f,0,\ldots)$ in
  Theorem \ref{theo:super-critical}, and then use \reff{eq:def-NOf1}.
\end{proof}

We directly deduce the following Corollary.
\begin{cor}
   \label{cor:super-crit-d1}
Under the hypothesis of Theorem \ref{theo:super-critical}, if $\alpha$
is the only eigen-value of $\cq$ with modulus equal to $\alpha$ (and
thus $J$ is
reduced to a singleton), then we have:
\[
(2\alpha^{2})^{-n/2} N_{n, \emptyset}(\bF)
\; \xrightarrow[n\rightarrow \infty ]{\P}  \;
 \sum_{\ell\in \N} (2\alpha)^{-\ell}  M_{\infty }(f_\ell) ,
\]
where, for $f\in F$,  $M_{\infty }(f)=\lim_{n\rightarrow \infty } (2\alpha)^{-n}
M_{\G_n} (\crr (f))$, and $\crr$ is the projection on the eigen-space
associated to the eigen-value $\alpha$. 
\end{cor}

The  Corollary \ref{cor:super-crit-d1-2} is then a direct consequence of Corollary
\ref{cor:super-crit-d1}. 



\subsection{Proof of Lemma \ref{lem:martingale}}
\label{sec:martingaleSuC}

Let $f \in L^{2}(\mu)$ and $j\in J$. Use that $\crr_j(L^{2}(\mu))\subset \C L^{2}(\mu)$ to deduce that
$\E\left[|M_{n,j}(f)|^2\right]$ is finite. 
We have for $n\in \N^*$:
\begin{align*}
\E[M_{n,j}(f)|\Hh_{n-1}] 
&= (2\alpha_j)^{-n} \sum_{i\in\GG_{n-1}} \E[\crr_j f(X_{i0}) +
  \crr_j f(X_{i1})|\Hh_{n-1}]\\ 
&= (2\alpha_j)^{-n} \sum_{i\in\GG_{n-1}} 2\, \cq \crr_j f(X_{i}) \\ 
&= (2\alpha_j)^{-(n-1)}\sum_{i\in\GG_{n-1}} \crr_j f (X_{i}) \\
& = M_{n-1,j}(f),
\end{align*}
where the second equality follows from branching Markov property and the
third  follows from the fact that $\crr_j$ is the projection on the
eigen-space associated to the eigen-value $\alpha_j$ of $\cq$. This
gives that $M_{j}(f)$ is a $\ch$-martingale. 
We also  have, writing $f_j$ for $\crr_j(f)$:
\begin{align}
\nonumber
\E\left[|M_{n,j}(f)|^2\right]
&= (2\alpha)^{-2n} \, \E\left[M_{\G_n} (f_j) M_{\G_n} (\overline
  f_j)\right]\\
\nonumber
&= (2\alpha^2)^{-n} \,\langle \nu, \cq^n (|f_j|^2)  \rangle
+ (2\alpha)^{-2n} \sum_{k=0}^{n-1} 2^{n+k}\langle \nu, \cq^{n-k-1}
  \cp\left(\cq^k f_j  \sot \cq^{k} \overline f_j\right) \rangle\\
\nonumber
& \leq \, C \, (2\alpha^{2})^{-n} \,\langle \mu, \cq^{n-k_{0}} (|f_j|^2)  \rangle + (2\alpha)^{-2n} \sum_{k=0}^{n-1} 2^{n+k}\langle \nu, \cq^{n-k-1} \cp\left(|\cq^k f_j|  \otimes^{2} \right) \rangle\\
\label{eq:mart-majo-L2} &\leq C (2\alpha^{2})^{-n} \|f_{j}\|^{2}_{L^{2}(\mu)} + C \, (2\alpha^{2})^{-n} \sum_{k = 0}^{n-k_{0}} 2^{k} \, \|\Qq^{k} f_{j}\|^{2}_{L^{2}(\mu)}
\end{align}
where  we used  the  definition  of $M_{n,j}$  for  the first  equality,
\eqref{eq:Q2-bis} with $m=n$ for the  second equality, Assumption \ref{hyp:Q} (ii) for the first term of the first inequality, the fact that $\cq^k f_j  \sot \cq^{k} \overline f_j \leq |\Qq^{k} f_{j}| \otimes^{2}$ for the second term of the first inequality and for the last inequality, we followed the lines of the proof of Lemma \ref{lem:L2MG}.  Finally, using that $|\Qq^{k} f_{j}| = \alpha^{k} |f_{j}|$, this implies     that $\sup_{n\in\N}\E\left[|M_{n,j}(f)|^2\right] < + \infty.$ Thus the martingale  $M_{j}(f)$ converges a.s. and  in  $L^{2}$  towards a limit.

\subsection{Proof of Theorem \ref{theo:super-critical}}
\label{sec:proof-super-thm}
Recall the sequence $(\beta_{n},n\in\NN)$ defined in Assumption
\ref{hyp:Q2} and  the $\sigma$-field $\Hh_{n} = \sigma\{X_{u},u\in\TT_{n}\}$.
Let  $(\p_{n},n\in\NN)$ be a sequence of integers such that $\p_n$ is even
and (for $n\geq 3$):
\begin{equation}\label{eq:hypSCrit-p}
\frac{5n}{6} <\p_n<n, \quad  \lim_{n\rightarrow \infty } (n-\p_n)=\infty
\quad  \text{and} \quad \lim_{n\rightarrow \infty} 
\alpha^{-(n-\p_n)} \beta_{\p_{n}/2} = 0. 
\end{equation}
Notice such sequences exist. 
When there is no ambiguity, we shall write $\p$  for $\p_n$. Using Remark \ref{rem:Nn=Nnk0}, it suffices to do the proof with $N_{n,\emptyset}^{[k_{0}]}(\bF)$ instead of $N_{n,\emptyset}(\bF).$
We deduce from \eqref{eq:nof-D} that:
\begin{equation}\label{eq:decomSCrit}
N^{[k_0]}_{n, \emptyset}(\bF) = R_0^{k_0} (n) + R_{4}(n) + T_{n}(\bF),
\end{equation}
with notations from \reff{eq:N=D+R} and   \reff{eq:reste01}:
\begin{align*}
R_{0}^{k_{0}}(n) &= |\GG_{n}|^{-1/2} \sum_{k=k_{0}}^{n-\p_n-1}
  M_{\GG_{k}}(\tilde{f}_{n-k}),\\
 T_{n}(\bF)=R_1(n) &= \sum_{i\in\GG_{n-\p_n}} \EE[N_{n,i}(\bF)|\Hh_{n-\p_n}],\\
 R_{4}(n)=\Delta_n &= \sum_{i\in \GG_{n-\p_n}} \left(N_{n,i}(\bF) -
            \EE[N_{n,i}(\bF)|\Hh_{n-\p_n}]\right). 
\end{align*}

Furthermore, using the branching Markov property, we get for all $i\in \G_{n-\p_n}$:
\begin{equation}\label{eq:bmpSCrit}
\EE[N_{n,i}(\bF)|\Hh_{n-\p_n}] = \EE[N_{n,i}(\bF)|X_{i}].
\end{equation}
We have the following elementary lemma.
\begin{lem}\label{lem:cvR0SCrit}
Under the assumptions of Theorem \ref{theo:super-critical}, we have the
following convergence: 
\begin{equation*}
\lim_{n\rightarrow \infty} (2\alpha^{2})^{-n}\, \EE\left[R^{[k_{0}]}_{0}(n)^{2}\right] = 0.
\end{equation*}
\end{lem}
\begin{proof}
We follow the proof of Lemma \ref{lem:cvR0-L2}. As $2\alpha^{2} > 1$ and using the first inequality of \eqref{eq:IRk0nL2sub} we
get that for some constant $C$ which does not depend on $n$ or $\p$:
\begin{equation*}
\EE\left[R_{0}^{k_{0}}(n)^{2}\right]^{1/2} \leq C\, 
2^{-\p/2} (2\alpha^2)^{(n-\p)/2} .
\end{equation*}
It follows from the previous inequality that $
(2\alpha^{2})^{-n} \E\left[R_{0}(n)^{2}\right] \leq 
C  (2\alpha)^{-2\p}$. 
Then use $2\alpha > 1$ and $\lim_{n\rightarrow \infty} \p = \infty$ to conclude.
\end{proof}
Next, we have the following lemma.
\begin{lem}\label{lem:cvR4SCrit}
Under the assumptions of Theorem \ref{theo:super-critical}, we have the following convergence:
\begin{equation*}
\lim_{n\rightarrow \infty} (2\alpha^{2})^{-n}\EE\left[R_{4}(n)^{2}\right] = 0.
\end{equation*}
\end{lem}
\begin{proof}
First, we have:
\begin{align}
\EE[R_{4}(n)^{2}] 
&= \EE\left[\left(\sum_{i\in \GG_{n-\p}} (N_{n,i}(\bF) -
  \EE[N_{n,i}(\bF)|X_{i}])\right)^{2}\right] \nonumber\\
 &= \EE\left[\sum_{i\in \GG_{n-\p}}\EE[(N_{n,i}(\bF) -
   \EE[N_{n,i}(\bF)|X_{i}])^{2}|\Hh_{n-\p}]\right]
 \nonumber \\ 
&\leq  \EE\left[\sum_{i\in \GG_{n-\p}}\EE[N_{n,i}(\bF)^2|X_i]\right],
\label{eq:R4SCrit1} 
\end{align}
where we used \eqref{eq:bmpSCrit} for the first equality and the branching Markov chain property for the second and the last inequality.  Note that for all $i\in\GG_{n-\p}$ we have 
\begin{equation*}
\EE\left[\EE[N_{n,i}(\bF)^2|X_i]\right] = |\GG_{n}|^{-1}  \EE\left[\EE\left[\left(\sum_{\ell=0}^{\p} M_{i\GG_{\p-k}}(\tilde{f}_{\ell})\right)^{2}|X_{i}\right]\right], 
\end{equation*}
where we used the definition of $N_{n,i}(\bF)$. Putting the latter equality in \eqref{eq:R4SCrit1} and using the first inequality of \eqref{eq:EGn-nu0}, we get
\begin{equation*}
\EE[R_{4}(n)^{2}] \leq |\GG_{n}|^{-1} \, \EE[M_{\GG_{n-p}}(h_{\hat{p}})]
\leq C \, 2^{-\hat{p}} \langle \mu,h_{\p} \, \rangle, \quad \text{with}
\quad h_{\p}(x) = \EE_{x}[(\sum_{\ell = 0}^{\p} M_{\GG_{\p -
    \ell}}(\tilde{f}))^{2}]. 
\end{equation*}
Using the second inequality of \eqref{eq:EGn-nu0} and \eqref{eq:L2-erg}, we get
\begin{equation*}
\langle \mu,h_{\p} \rangle = \EE_{\mu}[(\sum_{\ell = 0}^{\p} M_{\GG_{\p
    - \ell}}(\tilde{f}))^{2}]
\leq \left(\sum_{\ell = 0}^{p}
  \EE_{\mu}[(M_{\GG_{p-\ell}}(\tilde{f}_{\ell}))^{2}]^{1/2} \right)^2
\leq C \,  (2 \alpha)^{2\p}.
\end{equation*}
This implies that 
\begin{equation*}
(2\alpha^{2})^{-n}\EE\left[R_{4}(n)^{2}\right] \leq C \,  (2\alpha^{2})^{-n} \, (2 \alpha^{2})^{\p}  \,  =  \, C \,  (2\alpha^{2})^{\p - n}.
\end{equation*}
We then conclude using $2\alpha^{2} > 1$ and \eqref{eq:hypSCrit-p}. 
\end{proof}

Now, we study the third term of the right hand side of
\eqref{eq:decomSCrit}. First, note that: 
\begin{align*}
T_{n}(\bF) 
&= \sum_{i\in\GG_{n-\p}}  \EE[N_{n,i}(\bF)|X_{i}] \\ 
&= \sum_{i\in\GG_{n-\p}}|\GG_{n}|^{-1/2} \sum_{\ell = 0}^{\p}
  \EE_{X_{i}}[M_{\GG_{\p - \ell}}(\tilde{f}_{\ell})] \\ 
&= |\GG_{n}|^{-1/2} \sum_{i\in\GG_{n-\p}} \sum_{\ell = 0}^{\p} 2^{\p -
  \ell} \Qq^{\p - \ell}(\tilde{f}_{\ell})(X_{i}), 
\end{align*}
where we used \eqref{eq:bmpSCrit} for the first equality, the definition
\reff{eq:def-NiF}  of  $N_{n}(\bF)$  for  the  second  equality  and
\eqref{eq:Q1} for the last equality. Next, projecting in the eigen-space
associated to the eigenvalue $\alpha_{j}$, we get
\begin{equation*}\label{eq:decomTSCrit}
T_{n}(\bF) = T^{(1)}_{n}(\bF) + T^{(2)}_{n}(\bF),
\end{equation*}
where, with $\hat f= f - \langle \mu, f \rangle - \sum_{j\in J}
\crr_j(f)$ defined in \reff{eq:fhatcrit-S}:
\begin{align*}
T^{(1)}_{n}(\bF) 
&= |\GG_{n}|^{-1/2} \sum_{i\in \GG_{n-\p}} \sum_{\ell=0}^{\p} 2^{\p-\ell}
  \left(\Qq^{\p-\ell}(\hat{f}_{\ell}) \right)(X_{i}),\\ 
T^{(2)}_{n}(\bF) 
& = |\GG_{n}|^{-1/2} \sum_{i\in \GG_{n-\p}}
  \sum_{\ell=0}^{\p} 2^{\p-\ell} 
 \alpha^{\p-\ell}\sum _{j\in J} 
\theta_j^{\p-\ell}\crr_j (f_\ell)
(X_{i}). 
\end{align*}
We have the following lemma.
\begin{lem}\label{lem:cvT1SCrit}
Under the assumptions of Theorem \ref{theo:super-critical}, we have the following convergence:
\begin{equation*}
\lim_{n\rightarrow \infty}(2\alpha^{2})^{-n/2} \EE[|T^{(1)}_{n}(\bF)|] = 0.
\end{equation*}
\end{lem}
\begin{proof}
Recall $\p$ is even. We set $h_{\p} = \sum_{\ell=0}^{\p} 2^{\p-\ell} \Qq^{\p-\ell}(\hat {f}_{\ell}).$ We have:
\begin{align*}
(2\alpha^{2})^{-n/2}\EE[|T^{(1)}_{n}(\bF)|] 
  &\leq (2\alpha)^{-n} \EE[M_{\GG_{n-\p}}(|h_{\p}|)] \\
  &\leq C \, (2\alpha)^{-n} 2^{n-\p} \langle \mu, |h_{\p}| \rangle \\ 
  & \leq C\, (2\alpha)^{-n} 2^{n-\p} \|h_{\p}\|_{L^{2}(\mu)} \\
  &\leq C \, (2\alpha)^{-n} 2^{n-\p} \sum_{\ell=0}^{\p} 2^{\p-\ell}
    \alpha^{\p-\ell}  \beta_{\p-\ell} \norm{f_\ell}_{L^{2}(\mu)}   \\  
&= C \, \sum_{\ell=0}^\p 2^{-\ell}\alpha^{-(n-\p+\ell)}  \beta_{\p-\ell} \,,
\end{align*}
where we used the definition of $T^{(1)}_{n}(\bF)$ for the first
inequality, the first equation of \eqref{eq:EGn-nu0} for the second,
Cauchy-Schwartz inequality for the third  and \eqref{eq:L2-erg-crit} for
the last inequality.  We have:
\[
\sum_{\ell=0}^{\p/2} 2^{-\ell}\alpha^{-(n-\p+\ell)}  \beta_{\p-\ell} 
\leq \alpha^{-(n-\p)}  \beta_{\p/2}  \sum_{\ell=0}^{\p/2} (2\alpha)^{-\ell}. 
\]
Using the third condition in \reff{eq:hypSCrit-p} and that $2\alpha>1$, we deduce the right
hand-side converges to $0$ as $n$ goes to infinity. Without loss of
generality, we can assume that the sequence $(\beta_n, n\in \N^*)$
is bounded by 1. Since $\alpha>1/\sqrt{2}$, we also have:
\[
  \sum_{\ell=\p/2}^\p 2^{-\ell}\alpha^{-(n-\p+\ell)}  \beta_{\p-\ell}
  \leq  (1- 2\alpha)^{-1}  \, 2^{-\p/2} \alpha^{ -n + \p/2} \leq 
  (1- 2\alpha)^{-1}  \, 2^{ n/2 - 3\p/4}.
\]
Using that $n/2 - 3\p/4<
-n/8$, thanks to  the first condition in \reff{eq:hypSCrit-p},  we deduce the right
hand-side converges to $0$ as $n$ goes to infinity.
  Thus, we get  that $ \lim_{n \rightarrow \infty}(2\alpha^{2})^{-n/2}
  \EE[|T^{(1)}_{n}(\bF)|] = 0$. 
\end{proof}

Now, we deal with the term $T^{(2)}_{n}(\bF)$ in the following
result. Recall $M_{\infty , j}$ defined in Lemma \ref{lem:martingale}. 
\begin{lem}\label{lem:cvT2SCrit}
Under the assumptions of Theorem \ref{theo:super-critical}, we have the
following convergence: 
\begin{equation*}
(2\alpha^{2})^{-n/2}T^{(2)}_{n}(\bF) - \sum_{\ell\in \N} (2\alpha)^{-\ell} \sum_{j\in J} \theta_j^{n-\ell}M_{\infty  ,j}(f_\ell)  \; \xrightarrow[n\rightarrow \infty ]{\P}  \; 0.
\end{equation*}
\end{lem}
\begin{proof}
By definition of $T^{2}_{n}(\bF)$, we have
$T^{2}_{n}(\bF)=2^{-n/2} \sum_{\ell = 0}^{\p} (2\alpha)^{n-\ell} \sum_{j\in J}
\theta_{j}^{n-\ell} M_{n,j}(f_{\ell})$ and thus:
\begin{multline}\label{eq:T2-0}
(2\alpha^{2})^{-n/2}T^{(2)}_{n}(\bF) - \sum_{\ell\in \N}
(2\alpha)^{-\ell} \sum_{j\in J} \theta_j^{n-\ell}M_{\infty  ,j}(f_\ell)
\\ 
= \sum_{\ell = 0}^{\p} (2\alpha)^{-\ell} \sum_{j\in J}
\theta_{j}^{n-\ell} (M_{n,j}(f_{\ell}) - M_{\infty,j}(f_{\ell})) -
\sum_{\ell = \p+1}^{\infty} (2\alpha)^{-\ell} \sum_{j\in J}
\theta_{j}^{n-\ell} M_{\infty,j}(f_{\ell}). 
\end{multline}
Using that $|\theta_j|=1$, we get:
\[
\EE[|\sum_{\ell = 0}^{\p} (2\alpha)^{-\ell} \sum_{j\in J}
\theta_{j}^{n-\ell} (M_{n,j}(f_{\ell}) - M_{\infty,j}(f_{\ell}))| ] 
\leq \sum_{\ell = 0}^{\p} (2\alpha)^{-\ell} \sum_{j\in J}
\EE[|M_{n,j}(f_{\ell}) - M_{\infty,j}(f_{\ell})|].
\]
 Now, using that $(f_{\ell}, \ell\in \N)$ is uniformly bounded in $L^{2}(\mu)$, a close inspection of the proof of Lemma
 \ref{lem:martingale}, see \reff{eq:mart-majo-L2}, 
 reveals us that there exists a finite constant $C$ (depending on $\bF$)
 such that  for all $j\in J$, we have:
\begin{equation*}
\label{eq:supMn-bd}
\sup_{\ell\in\NN}\sup_{n\in\NN}\EE[|M_{n,j}(f_{\ell})|^{2}] \leq C.
\end{equation*}
The $L^2(\nu)$ convergence in Lemma \ref{lem:martingale} yields that:
\begin{equation}
\label{eq:supMi-bd}
\sup_{\ell\in\NN}\EE[|M_{\infty ,j}(f_{\ell})|^{2}] \leq C
\quad\text{and}\quad
\sup_{\ell \in \NN} \sup_{n \in \NN} \sum_{j\in J}\EE[|M_{n,j}(f_{\ell}) - M_{\infty,j}(f_{\ell})|] < 2|J|\sqrt{C}.
\end{equation}
Since Lemma \ref{lem:martingale} implies that $\lim_{n\rightarrow \infty
} \EE[|M_{n,j}(f_{\ell}) - M_{\infty,j}(f_{\ell})|]=0$, we deduce, as $2\alpha>1$ by the
dominated convergence theorem that:
\begin{equation}\label{eq:T2-1}
\lim_{n\rightarrow +\infty}\EE[|\sum_{\ell = 0}^{\p} (2\alpha)^{-\ell} \sum_{j\in J} \theta_{j}^{n-\ell} (M_{n,j}(f_{\ell}) - M_{\infty,j}(f_{\ell}))|] = 0.
\end{equation}

On the other hand, we have
\begin{equation}\label{eq:T2-2}
\EE[|\sum_{\ell = \p+1}^{\infty} (2\alpha)^{-\ell} \sum_{j\in J}
  \theta_{j}^{n-\ell} M_{\infty,j}(f_{\ell})|] 
\leq \sum_{\ell = \p+1}^{\infty} (2\alpha)^{-\ell} \sum_{j\in J}
  \EE[|M_{\infty,j}(f_{\ell})|] 
 \leq |J|\sqrt{C}\, \sum_{\ell = \p+1}^{\infty} (2\alpha)^{-\ell} , 
\end{equation}
where we used $|\theta_{j}| = 1$ for the first inequality and the
Cauchy-Schwarz inequality  and \reff{eq:supMi-bd}
for the second inequality. 
Finally, from \eqref{eq:T2-0}, \eqref{eq:T2-1} and \eqref{eq:T2-2} (with
$\lim_{n\rightarrow \infty} \sum_{\ell = \p+1}^{\infty} (2\alpha)^{-\ell}
= 0$) , we get the result of the lemma.
\end{proof}

\end{document}